\documentclass[final]{siamltex}

\usepackage{latexsym,amsmath,amssymb,amsfonts,amscd,multirow,dsfont}
\usepackage{epsfig,verbatim,epstopdf,graphics}
\usepackage{bm,color}
\usepackage{subfigure}

\newtheorem{thm}{Theorem}[section]

\newtheorem{lem}[thm]{Lemma}
\newtheorem{exa}[thm]{Example}
\newtheorem{rem}[thm]{Remark}
\newtheorem{prop}[thm]{Property}

\newcommand{\bk}{\mathbf{k}}
\newcommand{\bx}{\mathbf{x}}

\newcommand{\bE}{\mathbf{E}}

\newfont{\iams}{msbm9}

\newcommand{\commentbis}[1]{}
\newcommand{\be}{\begin{eqnarray}}
\newcommand{\ee}{\end{eqnarray}}
\newcommand{\beno}{\begin{eqnarray*}}
\newcommand{\eeno}{\end{eqnarray*}}

\newcommand{\barr}[1]{\begin{array}{#1}}
\newcommand{\earr}{\end{array}}

\newcommand{\beq}{\begin{equation}}
\newcommand{\eeq}{\end{equation}}
\newcommand{\beqa}{\begin{eqnarray}}
\newcommand{\eeqa}{\end{eqnarray}}

\newcommand{\bv}{{\bf v}}

\newcommand{\bJ}{{\bf J}}

\newcommand{\bV}{{\bf V}}
\newcommand{\bn}{{\bf n}}
\newcommand{\bzer}{{\mathbf{0}}}

\newcommand{\bP}{{\bf P}}

\newcommand{\bzero}{\mathbf{0}}
\newcommand{\bone}{\mathbf{1}}
\newcommand{\bl}{\mathbf{l}}
\newcommand{\bi}{\mathbf{i}}

\newcommand{\bj}{\mathbf{j}}
\newcommand{\bW}{\mathbf{W}}

\newcommand{\ba}{{\bm{a}}}
\newcommand{\bal}{{\bm{\alpha}}}
\newcommand{\bb}{{\bm{\beta}}}
\newcommand{\bq}{\mathbf{q}}

\title
{A Sparse Grid Discontinuous Galerkin Method for High-Dimensional
Transport Equations and Its Application to Kinetic Simulations}

\author{
 Wei Guo
\thanks{Department of Mathematics, Michigan State University,
East Lansing, MI 48824 U.S.A.
 {\tt wguo@math.msu.edu}}
\and
 Yingda Cheng
\thanks{Department of Mathematics, Michigan State University,
East Lansing, MI 48824 U.S.A.
 {\tt ycheng@math.msu.edu}. Research is supported by NSF grants DMS-1318186 and DMS-1453661.}
}

\date{\today}

\begin{document}

\maketitle

\begin{abstract}
In this paper, we develop  a sparse grid discontinuous Galerkin (DG) scheme for  transport equations  and applied it to  kinetic simulations. The method uses the weak formulations of traditional Runge-Kutta DG (RKDG) schemes for hyperbolic problems and 
is proven to be $L^2$ stable and convergent. A major advantage of the scheme lies in its low computational and storage cost  due to    the employed sparse finite element approximation space. This attractive feature is explored in simulating Vlasov and Boltzmann transport equations. Good performance in accuracy and conservation is  verified by numerical tests in up to four dimensions.
\end{abstract}

\begin{keywords}
discontinuous Galerkin methods; sparse grid; high-dimensional transport equations; Vlasov equation; Boltzmann equation.
\end{keywords}

\section{Introduction}

In this paper, we develop a sparse grid DG method for high-dimensional transport equations. 
High-dimensional transport problems are ubiquitous in science and engineering, and most evidently in kinetic  simulations where  it is necessary to track the evolution of  probability density functions of particles. Deterministic kinetic simulations are very demanding     due to the large computational and storage cost. To make the schemes more attractive comparing with the alternative probabilistic methods,  an appealing approach is to  explore the sparse grid techniques \cite{bungartz2004sparse, garcke2013sparse} with the aim of breaking the curse of dimensionality \cite{bellman1961adaptive}.
 In the context of wavelets or sparse grid methods for kinetic transport equations,  we  mention the work of using wavelet-MRA methods for Vlasov equations \cite{besse2008wavelet}, the combination technique for linear gyrokinetics \cite{kowitz2013combination},  sparse adaptive finite element method \cite{widmer2008sparse}, sparse discrete ordinates method
\cite{grella2011sparseo} and sparse tensor spherical harmonics \cite{grella2011sparse} for radiative transfer, among many others. 

This paper focuses on the DG method \cite{Cockburn_2000_history}, which is a class
of finite element methods using  discontinuous approximation space for the numerical solution and the test functions. The RKDG scheme \cite{Cockburn_2001_RK_DG} developed in a series of papers for hyperbolic equations became very popular  due to its provable convergence, excellent conservation properties and accommodation for adaptivity and parallel implementations. Recent years have seen great growth in the interest of applying DG methods to kinetic systems (see for example  \cite{Ayuso2009, heath2012discontinuous, qiu2011positivity, cheng_vp, cheng_gamba_proft_2010}) because of the  conservation properties and   long time performance of the resulting simulations. However, the DG method is still deemed too costly in a realistic setting, often requiring more  degrees of freedom than other high order numerical calculations. 

 Recently, we developed a sparse grid DG method for high-dimensional elliptic problems \cite{sparsedgelliptic}. A sparse DG finite element space has been constructed, reducing the degrees of freedom  from the standard {$O(h^{-d})$ to $O(h^{-1}|\log_2 h|^{d-1})$} for $d$-dimensional problems, where $h$ is the uniform mesh size in each dimension. The resulting scheme retains main properties of standard DG methods while making the computational cost tangebile for high-dimensional simulations. This  motivates the current work for the transport equations, and we use kinetic problems   as a test bed for the new algorithm.
The scheme in this paper uses weak formulation to guarantee many nice properties such as stability and conservation, while a detailed study of the approximation results is performed to obtain $L^2$ convergence rate of $O( (\log h)^d h^{k+1/2})$  for smooth enough solutions, where $k$ is the degree of polynomial. The method is demonstrated to be advantageous for kinetic simulations, because it can capture the main features of the solution with manageable cost and can conserve key macroscopic quantities in the mean time.  
The rest of this paper is organized as follows:
in Section \ref{sec:method},  we  construct the sparse grid DG formulations for linear transport equation with variable coefficient.  In Section \ref{sec:analysis}, we prove $L^2$ stability and error estimate for constant coefficient equations. The numerical performance is validated in Section \ref{sec:numerical} by   several benchmark tests. Section \ref{sec:kinetic} discusses the application of the scheme to Vlasov and Boltzmann equations, and  we conclude the paper with some remarks and future work  in Section  \ref{sec:Con}. The detail of the proof of a key lemma is gathered in the Appendix.

\section{Numerical method}
\label{sec:method}

In this section, we   construct the sparse grid DG method for
the following $d$-dimensional linear transport equation with variable coefficients on a box-shaped domain 
\begin{equation}
\label{eq:model}
\left\{\begin{array}{l}
u_t + \nabla\cdot(\ba(t,\bx) \,u) =0,\quad \bx\in\Omega=[0,1]^d,\\[2mm]
u(0,\bx) = u_0(\bx),
\end{array}\right.
\end{equation}
subject to suitable boundary conditions. 
We first review the DG finite element space on sparse grid introduced in \cite{sparsedgelliptic}, and then provide the formulation of scheme along with   implementation details. 

\subsection{DG finite element space on sparse grid}

In this subsection, we   prescribe the grid and the associated finite element space. Many of the discussions follow from our previous work for elliptic equations \cite{sparsedgelliptic}. 
  First, we introduce the hierarchical decomposition of piecewise polynomial space in one dimension on the interval $[0,1]$. We define a set of nested grids, where the $n$-th level grid $\Omega_n$ consists of $2^n$ uniform cells $I_{n}^j=(2^{-n}j, 2^{-n}(j+1)]$, $j=0, \ldots, 2^n-1,$ for any $n \ge 0.$
The nested grids result in the nested piecewise polynomial spaces. In particular, let
$$V_n^k:=\{v: v \in P^k(I_{n}^j),\, \forall \,j=0, \ldots, 2^n-1\}$$
  be the usual piecewise polynomials of degree at most $k$ on the $n$-th level grid $\Omega_n$. Then, we have $$V_0^k \subset V_1^k \subset V_2^k \subset V_3^k \subset  \cdots$$
We can now define the multiwavelet subspace $W_n^k$, $n=1, 2, \ldots $ as the orthogonal complement of $V_{n-1}^k$ in $V_{n}^k$ with respect to the $L^2$ inner product on $[0,1]$, i.e.,
\begin{equation*}
V_{n-1}^k \oplus W_n^k=V_{n}^k, \quad W_n^k \perp V_{n-1}^k.
\end{equation*}
For notational convenience, we let
 $W_0^k:=V_0^k$, which is standard piecewise polynomial space of degree $k$ on $[0,1]$.  The dimension of $W_n^k$ is $2^{n-1}(k+1)$ when $n \geq 1$, and $k+1$ when $n=0$.
In summary, we have found a hierarchical representation of  the standard piecewise polynomial space  $V_n^k$ on $\Omega_n$ as $V_n^k=\bigoplus_{0 \leq j\leq n} W_j^k$. 
 
 \medskip
Now we are ready to review the construction in multi-dimensions. First we recall some basic notations about multi-indices. For a multi-index $\mathbf{\alpha}=(\alpha_1,\cdots,\alpha_d)\in\mathbb{N}_0^d$, where $\mathbb{N}_0$  denotes the set of nonnegative integers, the $l^1$ and $l^\infty$ norms are defined as 
$$
|\bal|_1:=\sum_{m=1}^d \alpha_m, \qquad   |\bal|_\infty:=\max_{1\leq m \leq d} \alpha_m.
$$
The component-wise arithmetic operations and relational operations are defined as
$$
\bal \cdot \bb :=(\alpha_1 \beta_1, \ldots, \alpha_d \beta_d), \qquad c \cdot \bal:=(c \alpha_1, \ldots, c \alpha_d), \qquad 2^\bal:=(2^{\alpha_1}, \ldots, 2^{\alpha_d}),
$$
$$
\bal \leq \bb \Leftrightarrow \alpha_m \leq \beta_m, \, \forall m,\quad
\bal<\bb \Leftrightarrow \bal \leq \bb \textrm{  and  } \bal \neq \bb.
$$

By making use of the multi-index notation, we denote by $\bl=(l_1,\cdots,l_d)\in\mathbb{N}_0^d$ the mesh level in a multivariate sense. We  define the tensor-product mesh grid $\Omega_\bl=\Omega_{l_1}\otimes\cdots\otimes\Omega_{l_d}$ and the corresponding mesh size $h_\bl=(h_{l_1},\cdots,h_{l_d}).$ Based on the grid $\Omega_\bl$, we denote by $I_\bl^\bj=\{\bx:x_m\in(h_mj_m,h_m(j_{m}+1)),m=1,\cdots,d\}$ an elementary cell, and 
$$\bV_\bl^k:=\{\bv: \bv(\bx) \in P^k(I^{\bj}_{\bl}), \,\,  \bzero \leq \bj  \leq 2^{\bl}-\bone \}= V_{l_1,x_1}^k\times\cdots\times  V_{l_1,x_d}^k$$
the tensor-product piecewise polynomial space, where $P^k(I^{\bj}_{\bl})$ denotes the collection of polynomials of degree up to $k$ in each dimension on cell $I^{\bj}_{\bl}$. 
If we use equal mesh refinement of size $h_N=2^{-N}$ in each coordinate direction, the  grid and space will be denoted by $\Omega_N$ and $\bV_N^k$, respectively.  

Based on a tensor-product construction, the multi-dimensional increment space can be  defined as
$$\bW_\bl^k=W_{l_1,x_1}^k\times\cdots\times  W_{l_1,x_d}^k.$$ 
Therefore, space $\bV_\bl^k$ can be represented by
$$
\bV_\bl^k=\bigoplus_{0\leq j_1 \leq l_1, \ldots, 0 \leq j_d \leq l_d} \bW_\bj^k.
$$
In particular, we have the standard tensor-product polynomial space on $\Omega_N$ as
$$
\bV_N^k=\bigoplus_{\substack{ |\bl|_\infty \leq N\\\bl \in \mathbb{N}_0^d}} \bW_\bl^k.
$$
The sparse finite element approximation space on $\Omega_N$ we use in this paper, on the other hand,  is  defined by
$$ \hat{\bV}_N^k:=\bigoplus_{\substack{ |\bl|_1 \leq N\\\bl \in \mathbb{N}_0^d}}\bW_\bl^k.$$
This is a subset of   $\bV_N^k$, and   its number of degrees of freedom  scales as $O((k+1)^d2^NN^{d-1})$ \cite{sparsedgelliptic}, which is significantly less than that of $\bV_N^k$ with exponential dependence on $Nd$. This is the key for computational savings in high dimensions. 

\subsection{Formulation of the scheme}
\label{subsec:formulation}

In this subsection, we formulate a DG scheme with the sparse finite element space $\hat{\bV}_N^k$ for solving the model problem \eqref{eq:model}. For simplicity of discussion, we assume periodic boundary conditions but note that the discussion can be easily generalized to Dirichlet boundary conditions as well.
 First, we review some basic notations about jumps and averages for   piecewise functions defined on the grid $\Omega_N$. Let $T_h$ be the collection of all elementary cell $I^{\bj}_{N}, \quad 0 \leq j_m  \leq 2^{N}-1, \forall \,m=1, \ldots, d$.
 $\Gamma:=\bigcup_{T \in \Omega_N} \partial_T$ be the union of the interfaces for all the elements in $\Omega_N$ (here we have taken into account the periodic boundary condition when defining $\Gamma$) and $S(\Gamma):=\Pi_{T\in \Omega_N} L^2(\partial T)$ be the set of $L^2$ functions defined on $\Gamma$. For any $q \in S(\Gamma)$ and $\bq \in [S(\Gamma)]^d$,  we define their   averages $\{q\}, \{\bq\}$ and jumps $[q], [\bq]$ on the interior edges as follows. Suppose
$e$ is an interior edge shared by elements $T_+$ and $T_-$, we define the unit normal vectors $\bm{n}^+$ and  $\bm{n}^-$ on $e$ pointing exterior of $T_+$ and $T_-$, respectively, then
\begin{flalign*}
[ q] \  =\  \, q^- \bm{n}^- \, +  q^+ \bm{n}^+, & \quad \{q\} = \frac{1}{2}( q^- + 	q^+), \\
[ \bq] \  =\  \, \bq^- \cdot \bm{n}^- \, +  \bq^+ \cdot \bm{n}^+, & \quad \{\bq\} = \frac{1}{2}( \bq^- + \bq^+).
\end{flalign*}

The semi-discrete DG formulation for   \eqref{eq:model} is defined as follows: find $u_h\in\hat{\bV}_N^k$, such that
\begin{align}  
\label{eq:DGformulation}
\int_{\Omega}(u_h)_t\,v_h\,d\bx =& \int_{\Omega} u_h\ba\cdot\nabla v_h\,d\bx - \sum_{\substack{e \in \Gamma}}\int_{e} \widehat{\ba u_h} \cdot [v_h]\,ds,\quad  	\\
:= & A(u_h,v_h) \notag
\end{align}
for $\forall \,v_h \in \hat{\bV}_N^k,$ where $\widehat{\ba u_h}$ is defined on the element interface denotes a monotone numerical flux to ensure the $L^2$ stability of the scheme. In this paper, we use the upwind flux
\begin{equation}
\widehat{\ba u_h} = \ba\{u_h\} + \frac{|\ba\cdot \bn|}{2}[u_h],
\end{equation}
with $\bn = \bn^+$ or $\bn^-$ for the constant coefficient case.
More generally, for variable coefficients problems, we adopt the global Lax-Friedrichs flux
\begin{equation}
	\widehat{\ba u_h} = \{\ba u_h\} + \frac{\alpha}{2}[u_h],
\end{equation}
where $\alpha=\max_{\bx}{|\ba(\bx,t)\cdot \bn|}$, the maximum is taken for all possible $\bx$ at time $t$ in the computational domain. 
%

	
When implementing the scheme, we   need a set of bases  to represent the DG solution in the sparse approximation space $\hat{\bV}_N^k$. In \cite{sparsedgelliptic}, we used  the orthonormal basis functions of  $\hat{\bV}_N^k$ for the sparse IPDG method. Such bases are constructed based on the one-dimensional orthonormal multiwavelet bases first introduced in \cite{alpert1993class}. For completeness of the paper, we brief review the process of constructing the orthonormal bases. We refer readers to \cite{alpert1993class} and \cite{sparsedgelliptic} for more details. We start with the one-dimensional case. The case of mesh level $l=0$ is trivial. By using the scaled Legrendre polynomials, we can easily obtain a set of orthonormal bases in $W_0^k$ which are denoted by
$v^0_{i,0}(x),\quad i=1,\ldots,k+1.$ For the case of $l>0$, the orthonormal bases in $W_l^k$ are constructed in \cite{alpert1993class} and denoted by 
$$v^j_{i,l}(x),\quad i=1,\ldots,k+1,\quad j=0,\ldots,2^{l-1}-1.$$
Note that such multiwavelet bases retain the orthonormal property of wavelet bases for different mesh levels, i.e.,
$$\int_0^1 v^j_{i,l}(x)v^{j'}_{i',l'}(x)\,dx=\delta_{ii'}\delta_{ii'}\delta_{jj'}.$$
For the multi-dimensional cases, the basis functions for $\bW_{\bl}^k$ can be defined by a tensor-product construction
$$v^\bj_{\bi,\bl}(\bx)\doteq\prod_{m=1}^d v^{j_m}_{i_m,l_m}(x_m),\quad i_m=1,\ldots,k+1,\,j_m=0,\ldots,\max(0,2^{l_m-1}-1).$$
Therefore, a DG solution in the sparse approximation $\hat{\bV}_N^k$ can be written as 
$$u_h(\bx)=\sum_{\substack{ |\bl|_1\leq N\\ 
		\mathbf{0}\leq\bj\leq\max(2^{\bl-\mathbf{1}}-\mathbf{1},\mathbf{0})\\ \mathbf{1}\leq\bi\leq \bk+\mathbf{1}}} u^\bj_{\bi,\bl}v^\bj_{\bi,\bl}(\bx),$$
where $u^\bj_{\bi,\bl}$ denotes the corresponding degree of freedom. 

We use the  total variation diminishing (TVD) Runge-Kutta methods \cite{Shu_1988_JCP_NonOscill} to solve the ordinary differential  equations resulting from the sparse DG spatial discretization, $(u_h)_t = R(u_h).$ A commonly used third-order TVD Runge-Kutta method is given by
\begin{align*}
u_h^{(1)} &= u^{n} + \Delta t R(u^n_h),\\
u_h^{(2)} &= \frac{3}{4}u^{n} + \frac14 u_h^{(1)} +\frac14 \Delta t R(u_h^{(1)}),\\
u_h^{n+1} &= \frac{1}{3}u^{n} + \frac23 u_h^{(1)} +\frac23 \Delta t R(u_h^{(2)}),
\end{align*}
where $u_h^{n}$ denotes the numerical solution at time level $t=t^n$. 

Finally, we would like to make some remarks on the implementation issues. Unlike the traditional piecewise polynomial space, for which one element can only interact with itself and its immediate neighbors, the basis functions in the sparse space $\hat{\bV}^k_N$ are no longer locally defined due to the hierarchical structure, leading to additional challenges in implementation. In fact, it is crucial to  take full advantage of such a hierarchical (tree-like) structure when implementing the scheme to save computational cost. 
As for the numerical flux, the global Lax-Friedrichs flux is adopted since we are able efficiently compute the interface integral in \eqref{eq:DGformulation} by using the \emph {unidirectional principle}. Such an idea  has been used in the sparse IPDG method for solving variable coefficient elliptic problems. In particular,  we first project $\ba$
into space $\hat{\bV}^k_N$ and denote the resulting projection by $\ba_h$. Since $\ba_h$ is a separable function, the multi-dimensional interface integral in \eqref{eq:DGformulation} can be computed by evaluating multiplication of one-dimensional integrals. An advantage of this procedure is that we do not rely on numerical quadratures to compute the interface integrals, which can become quite complicated in the sparse grid setting.

\section{Stability and error estimate}
\label{sec:analysis}

In this  section, we provide an analysis of stability and error estimate for the DG scheme \eqref{eq:DGformulation} when $\ba$ is a constant vector.

\begin{thm}[$L^2$ stability] 
The DG scheme \eqref{eq:DGformulation} for \eqref{eq:model} is $L^2$ stable when $\ba$ is a constant vector, i.e.
\begin{equation}
\label{eq:stable}
\frac{d}{dt}\int_{\Omega} (u_h)^2\,d\bx =- \sum_{e\in\Gamma}\int_e\frac{|\ba\cdot\bn|}{2}|[u_h]|^2d s \leq 0.
\end{equation}
\end{thm}
\noindent{\it Proof}:  The proof follows the standard argument in showing $L^2$ stability for DG schemes. Let $v_h=u_h$ in the bilinear form, we have
\begin{align}
A(u_h, u_h)&=\int_{\Omega} u_h\ba\cdot\nabla u_h\,d\bx - \sum_{\substack{e \in \Gamma}}\int_{e} \widehat{\ba u_h} \cdot [u_h]\,ds \notag\\
&=\int_{\Omega}  \ba\cdot\nabla \left (\frac{u_h^2}{2} \right)\,d\bx - \sum_{\substack{e \in \Gamma}}\int_{e} \widehat{\ba u_h} \cdot [u_h]\,ds \notag\\
&=\sum_{\substack{e \in \Gamma}}\int_{e} \left[\frac{u_h^2}{2} \ba \right]\,ds  - \sum_{\substack{e \in \Gamma}}\int_{e} \left (\ba\{u_h\} + \frac{|\ba\cdot \bn|}{2}[u_h] \right)\cdot [u_h]\,ds \notag\\
&=- \sum_{e\in\Gamma}\int_e\frac{|\ba\cdot\bn|}{2}|[u_h]|^2d s \leq 0, \label{eq:astab}
\end{align}
and \eqref{eq:stable} immediately follows.
\hfill\ensuremath{\blacksquare}

\medskip
Next, we will establish $L^2$ error estimate of the sparse grid DG solution. 
Below we introduce some notations about norms and semi-norms. On the grid $\Omega_N$,
we use $\|\cdot\|_{H^s(\Omega_N)}$  to denote the standard broken Sobolev norm, i.e. $\|v\|^2_{H^s(\Omega_N)}=\sum_{\bzero \le \bj \le 2^\mathbf{N}-\mathbf{1}} \|v\|^2_{H^s(I_N^\bj)},$ where $\|v\|_{H^s(I_N^\bj)}$ is the standard Sobolev norm on $I_N^\bj,$ (and $s=0$ is used to denote the $L^2$ norm). Similarly,  we use $|\cdot|_{H^s(\Omega_N)}$ to  denote the broken Sobolev semi-norm, and $\|\cdot\|_{H^s(\Omega_\bl)}, |\cdot|_{H^s(\Omega_\bl)}$ to denote the broken Sobolev norm and semi-norm that are supported on a general grid $\Omega_\bl$.

For any set $L=\{i_1, \ldots i_r \} \subset \{1, \ldots d\}$, we define $L^c$ to be the complement set of $L$ in $\{1, \ldots d\}.$ For a non-negative integer $\alpha$ and set $L$,  we define the semi-norm on any domain denoted by $\Omega'$
%
\begin{flalign*}
|v|_{H^{\alpha,L}(\Omega')} :=  \left \| \left ( \frac{\partial^{\alpha}}{\partial x_{i_1}^{\alpha}} \cdots \frac{\partial^{\alpha}}{\partial x_{i_r}^{\alpha}}  \right ) v \right \|_{L^2(\Omega')},
\end{flalign*}
and 
$$
|v|_{\mathcal{H}^{q+1}(\Omega')} :=\max_{1 \leq r \leq d} \left ( \max_{\substack{L\subset\{1,2,\cdots,d\} \\|L|=r}} |v|_{H^{t+1, L}(\Omega')} \right ),$$
which is the norm for the mixed derivative of $v$ of at most degree $q+1$ in each direction.

The error estimate in Theorem \ref{thm:conv} relies on  the following  approximation properties of the  $L^2$  projection  onto   $\hat{\bV}_N^k$.

\begin{lem}
\label{thm:appx}
Let $\mathbf{P}$  be   the standard $L^2$ projection  onto the space $\hat{\bV}_N^k$, then for  $k \geq 1$, any $1 \leq q \leq \min \{p, k\}$, there exist constants $\bar{\bar{c}}_{k,s, q},\, B_{s}(k,q,d), \,\kappa_s(k, q, N)>0$, such that for any
$v \in \mathcal{H}^{p+1}(\Omega)$, $N\geq 1$, $d \geq 2$,  we have
$$
  | \bP v- v |_{H^s(\Omega_N)}\le  \left (\bar{\bar{c}}_{k,s, q}  +B_{s}(k,q,d) \kappa_s(k,q,N)^d \right )  2^{-N(q+1-s)} |v|_{\mathcal{H}^{q+1}(\Omega)},
$$
for $s=0, 1$,
where
\[\kappa_s(k,q,N) =
 \left\{
  \begin{array}{ll}
(N+1)C_{k,q},  & s=0,\\
2 C_{k,q}, & s=1,
  \end{array} \right.\]
  
\[B_{s}(k,q,d) =
 \left\{
  \begin{array}{ll}
2^{-(q+1)},  & s=0,\\
d^{3/2}\sqrt{\bar{c}_{k,q}} C_{k,q}^{-2}/2, & s=1,
  \end{array} \right.\]
  and the constants $C_{k,q}=\max(\tilde{c}_{k,0,q}, \hat{c}_{k,0}), \bar{c}_{k,q}=\max\left( \tilde{c}^2_{k,1,q}\hat{c}^{2}_{k,0}, \tilde{c}^2_{k,0,q}\hat{c}^2_{k,1}\right)$. $\tilde{c}_{k,s,q}, \hat{c}_{k,s}, \bar{\bar{c}}_{k,s, q}$  are constants defined in \eqref{eq:proj2}, \eqref{eq:proj3} and \eqref{eq:tensorproj}.
\end{lem}

 The proof of the lemma is provided in the Appendix. This lemma shows that the $L^2$ norm and $H^1$ semi-norm of the projection error scale  like $O(N^d 2^{-N(k+1)})$  and $O(2^{-Nk})$ with respect to $N$ when the function $v$ has bounded mixed derivatives up to  enough degrees. 
 \begin{rem}
The approximation properties for a particular projector onto the $C^0$ subset of sparse grid space have been established in \cite{schwab2008sparse} and later used in \cite{sparsedgelliptic} for showing convergence of the sparse grid DG methods for elliptic equations. However,   for hyperbolic equations, the error estimate depends on the specific property of the projector, and the projection in \cite{sparsedgelliptic} thus does not apply.
 \end{rem}
 
 \medskip
Now we are ready to establish the error estimate of the sparse grid DG scheme.
\begin{thm}[$L^2$ error estimate]
\label{thm:conv}
Let $u$ be the exact solution to \eqref{eq:linear_adv}, and $u_h$ be the numerical solution to the semi-discrete scheme \eqref{eq:DGformulation} with numerical initial condition $u_h(0)=\bP u_0$. For  $k \geq 1$,  $u_0 \in  \mathcal{H}^{p+1}(\Omega)$,  $1 \leq q \leq \min \{p, k\}$,
 $N\geq 1$, $d \geq 2$,  we have for all $t \geq 0,$
$$
\|u_h-u\|_{L^2(\Omega_N)}    \leq \left (2\sqrt{ C_d  ||\ba||_2 t} \,C_\star(k,q,d,N)+(\bar{\bar{c}}_{k,0, q}  +B_{0}(k,q,d) \kappa_0(k,q,N)^d)2^{-N/2}  \right) 2^{-N(q+1/2)} |u_0|_{\mathcal{H}^{q+1}(\Omega)},
$$
where $C_d$ is a generic constant with dependence only on $d$, $C_\star(k,q,d,N)=\max_{s=0,1} \left (\bar{\bar{c}}_{k,s, q}  +B_{s}(k,q,d) \kappa_s(k,q,N)^d \right).$
The constants $\bar{\bar{c}}_{k,s, q},\, B_{s}(k,q,d), \,\kappa_s(k, q, N)$ are defined in Lemma \ref{thm:appx}.
\end{thm}

\noindent{\it Proof}:  Denote
$$e=u_h-u=\xi-\eta,\quad\text{where}\quad \xi=u_h- \mathbf{P}u,\quad \eta = u- \mathbf{P}u,$$
and $\mathbf{P}$ is the standard $L^2$ projection of u onto the space $\hat{\bV}_N^k$. We can plug  $\xi \in \hat{\bV}_N^k$ in the semi-discrete error equation and obtain 
\begin{equation}
\label{eq:xitxi}
\int_{\Omega}\xi_t\xi d\bx = \int_{\Omega}\eta_t\xi d\bx +A(\xi,\xi) - A(\eta,\xi).
\end{equation}
Due to the definition of $L^2$ projection, $\int_{\Omega}\eta_t\xi d\bx=0.$
From \eqref{eq:astab}, we have
\begin{equation}
\label{eq:axixi}
A(\xi,\xi) =
- \sum_{e\in\Gamma}\int_e\frac{|\ba\cdot\bn|}{2}|[\xi]|^2d s.
\end{equation}
Next, we consider $A(\eta,\xi)$. Again, the definition of $A$ gives
\begin{align*}
A(\eta,\xi) =& \int_\Omega \eta\ba\cdot \nabla\xi d\bx -\sum_{e\in\Gamma}\int_e\ba\{\eta\}\cdot[\xi]d s - \sum_{e\in\Gamma}\int_e\frac{|\ba\cdot\bn|}{2}[\eta]\cdot[\xi]d s\notag.
\end{align*}
The first term on the right hand side is 0, since $\eta$ is orthogonal to space $\hat{\mathbf{V}}^k_N$ and each component of $\ba \cdot \nabla\xi$ belongs to $\hat{\mathbf{V}}^k_N$.
We can bound the other  two terms as follows.
\begin{align}
  -\sum_{e\in\Gamma}\int_e\ba\{\eta\}\cdot[\xi]d s =&  \sum_{e\in\Gamma}\int_e\{\eta\}\ba\cdot\bn^+(\xi^+ - \xi^-) d s\notag\\
  \leq &\sum_{e\in\Gamma}\int_e |\ba\cdot\bn| \{\eta\}^2 ds + \frac12\sum_{e\in\Gamma}\int_e \frac{|\ba\cdot\bn|}{2} [\xi]^2 ds \notag \\
\leq& \frac{\|\ba\|_2}{2} \sum_{T\in T_h} \|\eta\|^2_{L^2(\partial T)} + \frac12\sum_{e\in\Gamma}\int_e \frac{|\ba\cdot\bn|}{2} [\xi]^2 ds. \notag
\end{align}
and
\begin{align}
  -\sum_{e\in\Gamma}\int_e\frac{|\ba\cdot\bn|}{2}[\eta]\cdot[\xi]d s \le&  \frac12\sum_{e\in\Gamma}\int_e \frac{|\ba\cdot\bn|}{2} [\eta]^2 ds + \frac12\sum_{e\in\Gamma}\int_e \frac{|\ba\cdot\bn|}{2} [\xi]^2 ds \notag \\
  \leq& \frac{\|\ba\|_2}{2} \sum_{T\in T_h} \|\eta\|^2_{L^2(\partial T)} + \frac12\sum_{e\in\Gamma}\int_e \frac{|\ba\cdot\bn|}{2} [\xi]^2 ds. \notag
\end{align}
Hence, we get
\begin{align}
A(\eta,\xi) \leq \|\ba\|_2\sum_{T\in T_h} \|\eta\|^2_{L^2(\partial T)} + \sum_{e\in\Gamma}\int_e \frac{|\ba\cdot\bn|}{2} [\xi]^2 ds.\label{eq:aetaxi}
\end{align}
Combining \eqref{eq:xitxi}, \eqref{eq:axixi} and \eqref{eq:aetaxi} gives
\begin{align}
\frac{d}{dt}\|\xi\|^2_{L^2(\Omega_N)} \leq &  2\|\ba\|_2\sum_{T\in T_h} \|\eta\|^2_{L^2(\partial T)}.\label{eq:xil2}
\end{align}
To bound the last term, we use  the trace inequality \cite{arnold1982interior}:
$$\|\phi\|^2_{L^2(\partial T)}\leq C_d\left( \frac{1}{h_N}\|\phi\|^2_{L^2(T)}+h_N|\phi|^2_{H^1(T)}\right), \quad \forall \phi\in H^1(T),$$
where $C_d$ is a generic constant with dependence only on $d$.
Hence, by Lemma \ref{thm:appx} and also noting that $h_N=2^{-N}$, we have 
\begin{align*}
&\frac{d}{dt}\|\xi\|^2_{L^2(\Omega_N)} \leq  2 C_d \|\ba\|_2\left(\frac{1}{h_N}\|\eta\|^2_{L^2(\Omega_N)}+h_N|\eta|^2_{H^1(\Omega_N)}\right) \\
& \leq 4 C_d \|\ba\|_2 C_\star^2 2^{-2N(q+1/2)} |u_0|^2_{\mathcal{H}^{q+1}(\Omega)},
\end{align*}
where $C_\star(k,q,d,N)=\max_{s=0,1} \left (\bar{\bar{c}}_{k,s, q}  +B_{s}(k,q,d) \kappa_s(k,q,N)^d \right).$  If we take $u_h(0)=\bP u_0$, then $\xi(0)=0$, and we have
$$
\|\xi\|_{L^2(\Omega_N)} \leq 2\sqrt{ C_d  \|\ba\|_2 t} \,C_\star 2^{-N(q+1/2)} |u_0|_{\mathcal{H}^{q+1}(\Omega)}.
$$
Combining with the estimate for $\eta$ from Lemma \ref{thm:appx}, we are done.
\hfill\ensuremath{\blacksquare}

This theorem proves $L^2$ convergence rate of $O(N^d 2^{-N(k+1/2)})$ or $O( (\log h_N)^d h_N^{k+1/2})$ of the sparse grid DG solution when $u_0$ has enough smoothness measured in mixed derivatives. Compared with traditional DG schemes on Cartesian meshes \cite{Cockburn_01_SINUM_cartesian}, the convergence rate is suboptimal, partly due to the use of sparse finite element space that contributes the logarithmic factor and partly due to the use of $L^2$ projection.  For  linear hyperbolic equation, it is well known that the tensor-product of one-dimensional Gauss-Radau projection can be used to raise the convergence order. We leave detailed investigation of such error estimates to future study.

\section{Numerical tests}
\label{sec:numerical}

In this section, we present several numerical tests to validate the efficiency and efficacy of the proposed scheme  for solving the model equation \eqref{eq:model} in multi-dimensions. We use the third-order TVD-RK temporal discretization   and choose the time step $\Delta t$   as
\begin{align*}
\displaystyle\Delta t &= \frac{\text{CFL}}{\displaystyle\sum_{m=1}^d \frac{c_m}{h_N}},\quad \text{for}\quad k=1,\,2,\\[2mm]
\displaystyle\Delta t &= \frac{\text{CFL}}{\displaystyle\sum_{m=1}^d \frac{c_m}{h_N^{4/3}}},\quad \text{for}\quad k=3,
\end{align*}
for the purpose of accuracy test, where $c_m$ is the maximum wave propagation speed in $x_m$-direction and $CFL=0.1$.


\begin{exa}[Linear advection with constant coefficient]
\label{ex:linear}
We consider  
\begin{equation}
\label{eq:linear_adv}
\left\{\begin{array}{l} \displaystyle u_t + \sum_{m=1}^d u_{x_m} = 0,\quad \bx\in[0,1]^d,\\[2mm]
\displaystyle u(0,\bx) = \sin\left(2\pi\sum_{m=1}^d x_m\right),
\end{array}\right.
\end{equation}
with periodic boundary conditions.
\end{exa}

The exact solution   at $t=T$ is a smooth function, 
$$u(T,\bx) =  \sin\left(2\pi\left(\sum_{m=1}^d x_m - d\,T \right)\right).$$
In the simulation, we compute the numerical solutions up to two periods in time, meaning that we let final time $T=1$ for $d=2$, $T=2/3$ for $d=3$, and $T=0.5$ for $d=4$.
In Table \ref{table:linear_d2}, we report the degrees of freedom of the associated space, $L^2$ errors and  orders of accuracy for $k=1, 2, 3$ and up to dimension four. The degrees of freedom of the computational method are significantly reduced when compared with the traditional DG space. As for accuracy, we observe half order reduction from the optimal $(k+1)$-th order for high-dimensional computations (d=4). The order is slightly better for lower dimensions. The conclusions from this example agree well with the error estimate in Theorem \ref{thm:conv}.

\begin{table}
\caption{$L^2$ errors and orders of accuracy for Example \ref{ex:linear} at $T=1$ when $d=2$, $T=2/3$ when $d=3$, and $T=0.5$ when $d=4$. $N$ is the number of mesh levels, $h_N$ is the size of the smallest mesh in each direction, $k$ is the polynomial order, $d$ is the dimension. DOF denotes the degrees of freedom of the   sparse approximation space $\hat{V}^k_N$.  $L^2$ order is calculated with respect to  $h_N$.}
\centering
\begin{tabular}{|c|c |c c c|   c c c|  c c c|}
\hline
$N$ & $h_N$& DOF&$L^2$ error & order&  DOF&$L^2$ error & order &  DOF&$L^2$ error & order \\
\hline

  & &        \multicolumn{3}{c|}{$ k=1,  \,d=2$}          &     \multicolumn{3}{c|}{$ k=1, \, d=3$}   &       \multicolumn{3}{c|}{$ k=1, \, d=4$}     \\
\hline
3&$1/8$	&	80 	 & 3.62E-01	&	--   &  304 	&	6.58E-01	 &--	& 1008&  6.56E-01 	 &	 	-- 	 	\\
4 &$1/16$	&	192 &	 9.17E-02	&	1.98   & 	832	& 	3.72E-01	&	0.82	&	3072	& 	4.99E-01	&0.39 	\\
5 &$1/32	$	&	448& 	 1.90E-02	&	2.27 & 	2176	& 	1.19E-01	&	1.64	&	8832	& 	2.40E-01	&	1.06	\\
6 &$1/64	$	&	1024&	 	4.81E-03	&	1.98  	&	5504	& 	2.96E-02	&	2.01	&	24320	& 	9.84E-02	&	1.28	\\
7 &$1/128$		&	2304&	 	1.27E-03	&	1.92  	&	13568	&	8.85E-03	&	1.74	&	64768	& 	3.21E-02	&	1.62	\\

\hline
 &     &   \multicolumn{3}{c|}{$ k=2,  \,d=2$}          &     \multicolumn{3}{c|}{$ k=2, \, d=3$}   &       \multicolumn{3}{c|}{$ k=2, \, d=4$}     \\
 \hline
3&$1/8$	&	180& 1.48E-02	&	-- & 	1026	& 	5.17E-02	&	--	&   5103	 &	8.97E-02 	&	-- 	\\
4 &$1/16$		&	432&	 	2.13E-03	&	2.80 & 	2808	& 1.10E-02	&	2.23	&	15552	& 	2.80E-02	&	1.68	\\
5 &$1/32$		&	1008&	 	4.39E-04	&	2.28 & 	7344	& 1.79E-03	&	2.63&	44712	& 	5.82E-03	&	2.27	\\
6 &$1/64	$	&	2304&	 	4.45E-05	&	3.30 & 	18576	& 	3.97E-04	&	2.17	&	123120	& 1.37E-03	&	2.09	\\
7 &$1/128$		&	5184&	 	7.68E-06	&	2.54 & 	45792	& 	5.14E-05	&	2.95 &	327888	& 	2.58E-04	&	2.41	\\
\hline
  &   &     \multicolumn{3}{c|}{$ k=3,  \,d=2$}          &     \multicolumn{3}{c|}{$ k=3, \, d=3$}   &       \multicolumn{3}{c|}{$ k=3, \, d=4$}     \\
 \hline
3&$1/8$	&	320 &	6.36E-04	&	-- & 	2432&	 	2.10E-03	&	--	&	16128	& 	4.09E-03	&	--	\\
4 &$1/16$		&	768 &	8.93E-05	&	2.83 & 	6656&	 	2.37E-04	&	3.14	&	49152	& 	6.06E-04	&	2.75	\\
5 &$1/32	$	&	1792&	 	4.07E-06	&	4.46  	&	17408&	 	2.49E-05	&	3.25	&	141312	& 	6.85E-05	&	3.14	\\
6 &$1/64	$	&	4096&	 	3.47E-07	&	3.55  	&	44032&	 	1.83E-06	&	3.76	&	389120	& 7.19E-06	&	3.25	\\
7 &$1/128$		&	9216&	 	1.97E-08	&	4.14  	&	108544&	 	2.03E-07	&	3.18	&	1036288	& 	6.36E-07	&	3.50	\\
\hline

\end{tabular}
\label{table:linear_d2}
\end{table}

\begin{exa}[Solid body rotation]
\label{ex:rotation}
We consider solid-body-rotation problems, which are in the form of \eqref{eq:model} with
$$\ba = \left(-x_2+\frac12, x_1-\frac12\right),\quad\text{when}\quad d=2,$$
$$\ba = \left(-\frac{\sqrt{2}}{2}\left(x_2-\frac12\right), \frac{\sqrt{2}}{2}\left(x_1-\frac12\right)+ \frac{\sqrt{2}}{2}\left(x_3-\frac12\right),-\frac{\sqrt{2}}{2}\left(x_2-\frac12\right)\right),\quad\text{when}\quad d=3,$$
subject to periodic boundary conditions.
\end{exa}

Such benchmark tests are commonly used in the literature to assess performance of transport schemes. Here,  the initial profile traverses along circular trajectories centered at $(1/2,1/2)$ for $d=2$ and about the axis $\{x_1=x_3\}\cap \{x_2=1/2\}$ for $d=3$ without deformation, and it goes back to the initial state after $2\pi$ evolution.
The initial conditions are set to be the following smooth cosine bells (with $C^5$ smoothness),
\begin{equation}\label{eq:cosine} u(0,\bx)=\left\{\begin{array}{ll}b^{d-1}\cos^6\left(\frac{\pi r}{2b}\right),& \text{if}\quad r\leq b,\\
0,&\text{otherwise},
\end{array}\right.\end{equation}
where $b=0.23$ when $d=2$ and $b=0.45$ when $d=3$, and $r=|\bx-\bx_c|$ denotes the distance between $\bx$ and the center of the cosine bell with $\bx_c=(0.75,0.5)$ for $d=2$ and  $\bx_c=(0.5,0.55,0.5)$ for $d=3$. We use the global Lax-Friedrichs flux in computation. The implementation with the upwind flux is also performed for this example, and the unidirectional principle can  be applied to save cost. However,   little difference is observed, and the numerical result by the upwind flux is hence omitted for brevity. In Table \ref{table:rot}, we summarize the convergence study of the numerical solutions computed by the sparse DG method with space $\hat{\bV}_N^k$, $k=1,2,3$, including the $L^2$ errors and   orders of accuracy. For this variable coefficient equation, we observe at least $k$-th order convergence for all cases. The convergence rate for three dimensions are about half order lower than their two dimensional counterpart.

\begin{table}
	\caption{$L^2$ errors and orders of accuracy for Example \ref{ex:rotation} at $T=2\pi$. $N$ is the number of mesh levels, $h_N$ is the size of the smallest mesh in each direction, $k$ is the polynomial order, $d$ is the dimension. DOF denotes the degrees of freedom of the   sparse approximation space $\hat{V}^k_N$.  $L^2$ order is calculated with respect to  $h_N$.}
	\centering
	\begin{tabular}{|c|c |c c c|   c c c|  }
		\hline
		$N$ & $h_N$& DOF&$L^2$ error & order&  DOF&$L^2$ error & order \\
		\hline
		
		& &        \multicolumn{3}{c|}{$ k=1,  \,d=2$}          &     \multicolumn{3}{c|}{$ k=1, \, d=3$}  \\
		\hline
5   &   1/32 &   448    &  1.30E-02 & --   &  		2176& 3.47E-03 &     --   \\	
6	&   1/64 &	1024&    8.03E-03&	0.70   &	5504&	  1.62E-03& 1.10   \\
7	&	1/128 &  2304&  3.59E-03&	1.16  & 13568   &  6.27E-04
& 1.37

\\
8	&	1/256 &  5120   & 9.89E-04	&1.86 & 	32768& 2.15E-04& 1.55 
\\
9	&	1/512 & 11264& 2.04E-04	&2.28 & 77824& 6.34E-05	&1.76 \\	
\hline 
&     &   \multicolumn{3}{c|}{$ k=2,  \,d=2$}          &     \multicolumn{3}{c|}{$ k=2, \, d=3$}   \\
 \hline
5   &  1/32 & 1008  &  4.21E-03	& --       &	7344 &   4.20E-04  & --
\\ 
6	& 1/64 &	2304&	  1.03E-03 &   2.03 &	18576&  9.97E-05 &  2.08
\\
7	& 1/128 &	5184&	   1.40E-04	&	2.88 &	45792 &  2.83E-05  & 1.82
\\
8	& 1/256 &	11520&   1.78E-05	&	2.98 &	110592 & 6.53E-06	& 2.12
\\
9	& 1/512 &	25344&   2.48E-06	&	2.84 & 262656 & 	1.28E-06 &	2.36

\\
\hline 
&     &   \multicolumn{3}{c|}{$ k=3,  \,d=2$}          &     \multicolumn{3}{c|}{$ k=3, \, d=3$}   \\
\hline
4   & 1/16 &  768      &  4.26E-03 & -- & 6656 & 4.05E-04 & -- 
\\
5	& 1/32 &	1792&	7.80E-04 & 2.45 & 17408&  6.48E-05&2.64
\\
6	& 1/64 &	4096&   7.64E-05 & 	3.35 &44032&	7.15E-06	& 3.42
\\
7	& 1/128 &	9216&	7.15E-06	& 3.42 & 108544&	1.12E-06 &3.03
 \\
8	& 1/256 &	20480 & 6.61E-07& 3.44  & 262144 &1.51E-07 &	2.89	
 \\
\hline
	\end{tabular}
	\label{table:rot}
\end{table}

\begin{exa}[Deformational flow]
\label{ex:deformational}
We consider the two-dimensional deformational flow with velocity field
$$\ba=(\sin^2(\pi x_1)\sin(2\pi x_2)g(t),-\sin^2(\pi x_2)\sin(2\pi x_1)g(t)),$$
where $g(t)=\cos(\pi t/T)$ with $T=1.5$. 
 \end{exa}
We still adopt the cosine bell \eqref{eq:cosine} as the initial condition for this test, but with $\bx_c=(0.65,0.5)$ and $b=0.35$.
Note that the deformational test is more challenging than the solid body rotation due to the space and time dependent flow field. In particular, along the direction of the flow, the cosine bell deforms into a crescent shape at $t=T/2$ , then goes back to its initial state at $t=T$ as the flow reverses. In the simulations, we compute the solution up to $t=T$. The convergence study is summarized in Table \ref{table:deform_d2}.   {Similar order reduction is observed compared with Example \ref{ex:rotation}.  In Figure \ref{fig:defor}, we plot the contour plots of the numerical solutions at $t=T/2$ when the shape of the bell is greatly deformed, and $t=T$ when the solution is recovered into its initial state. It is observed that the sparse DG scheme with higher degree $k$ could better resolve the highly deformed solution structure. 

\begin{table}[htp]
	\caption{$L^2$ errors and orders of accuracy for Example \ref{ex:deformational} at $T=1.5$. $N$ is the number of mesh levels, $h_N$ is the size of the smallest mesh in each direction, $k$ is the polynomial order, $d$ is the dimension. DOF denotes the degrees of freedom of the   sparse approximation space $\hat{V}^k_N$.  $L^2$ order is calculated with respect to  $h_N$. $d=2$.
}
\centering
	\begin{tabular}{|c|c |c c c|}
		\hline
		$N$ & $h_N$& DOF&$L^2$ error & order\\
		\hline
		
		& &        \multicolumn{3}{c|}{$ k=1$}   \\
		\hline
5	& 1/32 &	448&	  1.40E-02	&	-- \\
6	& 1/64 &	1024&	 6.88E-03	&1.03 \\
7	& 1/128 &	2304& 2.65E-03	&	1.38\\
8	& 1/256 &	5120&	  8.48E-04&	1.65 \\
\hline
& &        \multicolumn{3}{c|}{$ k=2$}   \\
\hline
5	&1/64 &	1008&	  3.32E-03	&	-- \\
6	&1/128 &	2304& 1.31E-03	&	1.71\\
7	& 1/256 &	5184& 3.82E-04	&	2.33\\
8	&1/512 & 11520& 5.92E-05	&	2.61\\
\hline
& &        \multicolumn{3}{c|}{$ k=3$}   \\
\hline
5	&1/32 &	1792&	 6.28E-04	&	-- \\
6	&1/64 &	4096& 1.22E-04	& 2.36\\
7	& 1/128 &	9216& 1.37E-05	&	3.16\\
8	&1/256 &	20480& 2.02E-06	&	2.76
\\
\hline

\end{tabular}
\label{table:deform_d2}
\end{table}

\begin{figure}[htp]
\begin{center}
\subfigure[]{\includegraphics[width=.42\textwidth]{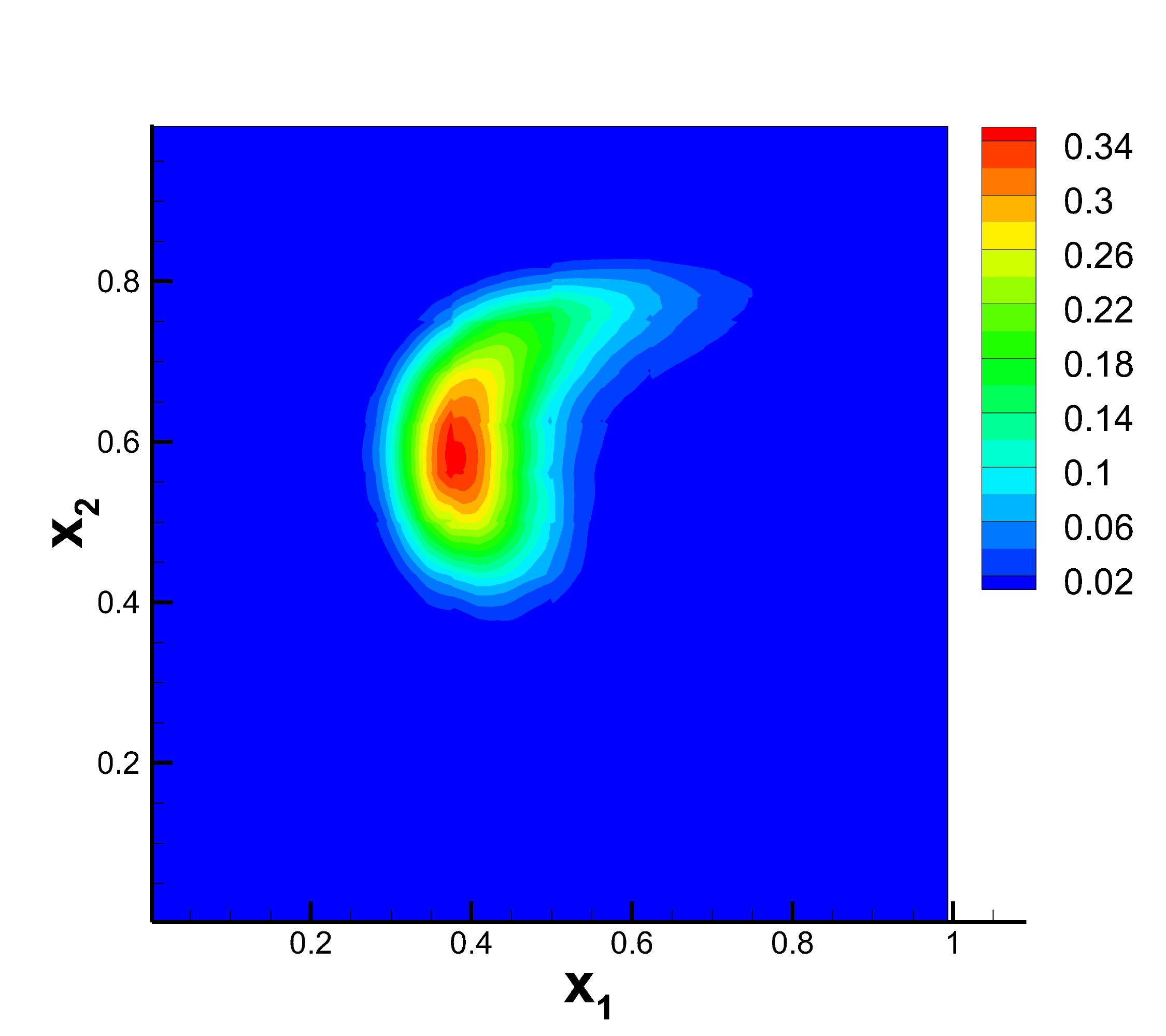}}
\subfigure[]{\includegraphics[width=.42\textwidth]{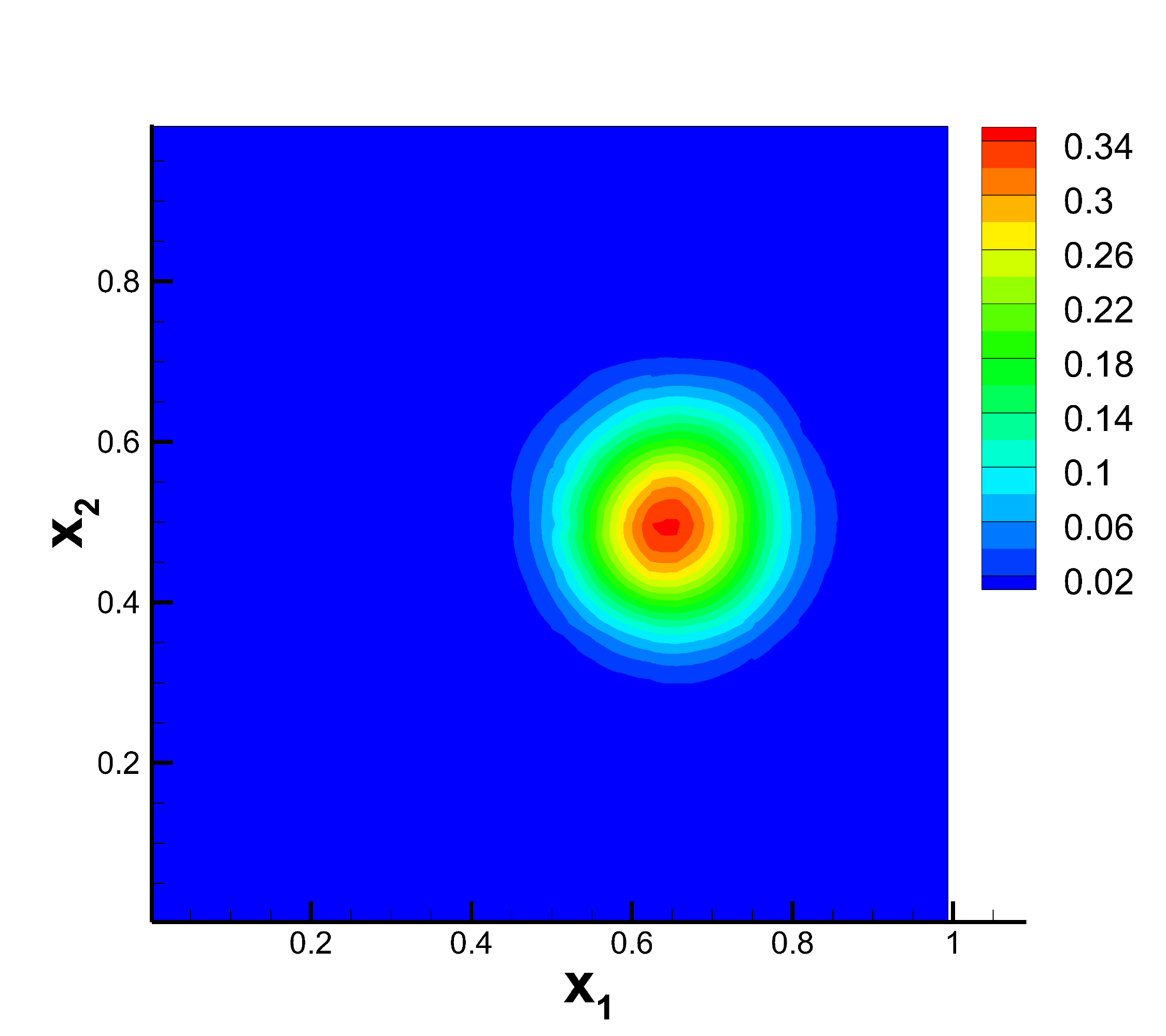}}\\
\subfigure[]{\includegraphics[width=.42\textwidth]{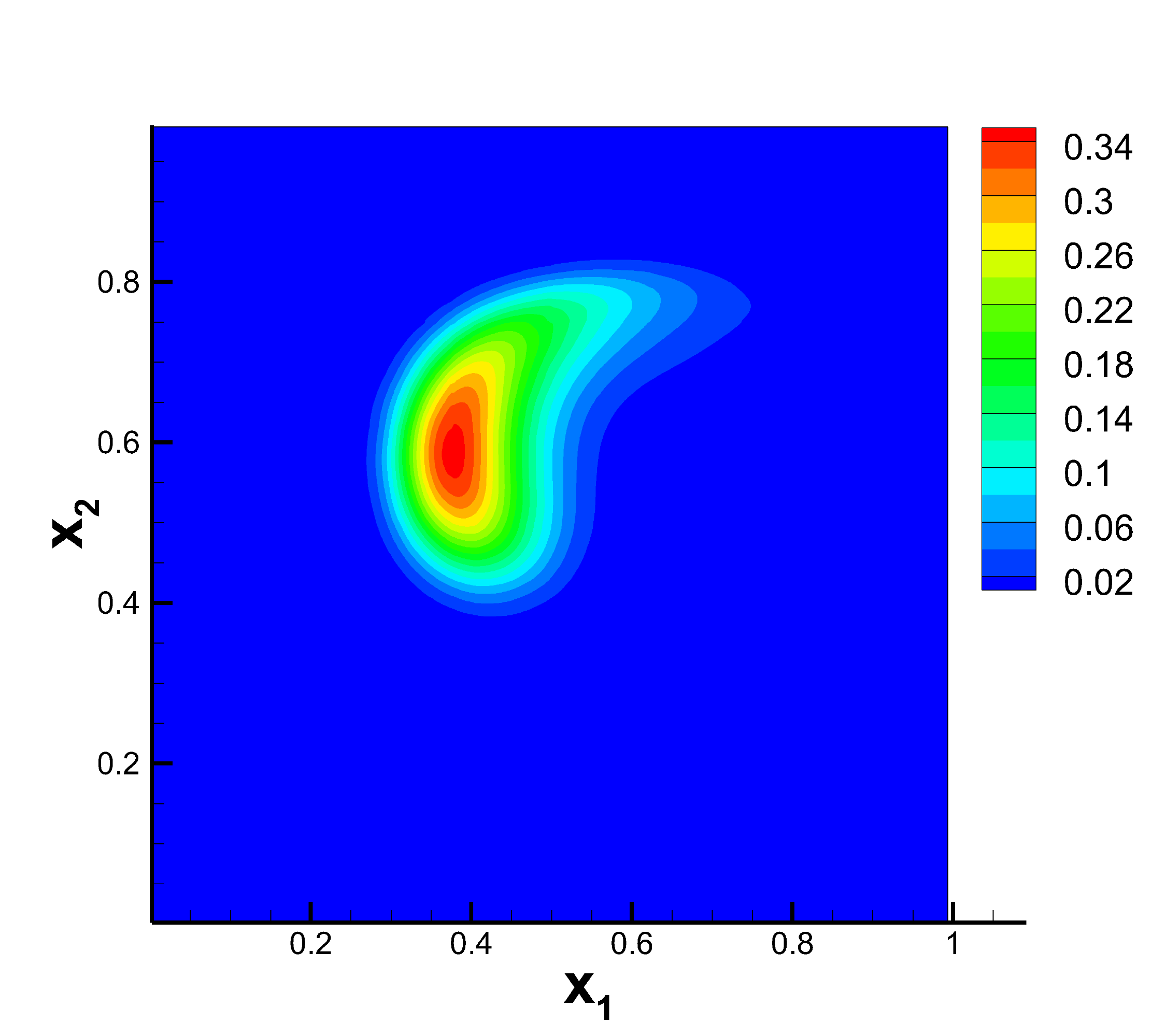}}
\subfigure[]{\includegraphics[width=.42\textwidth]{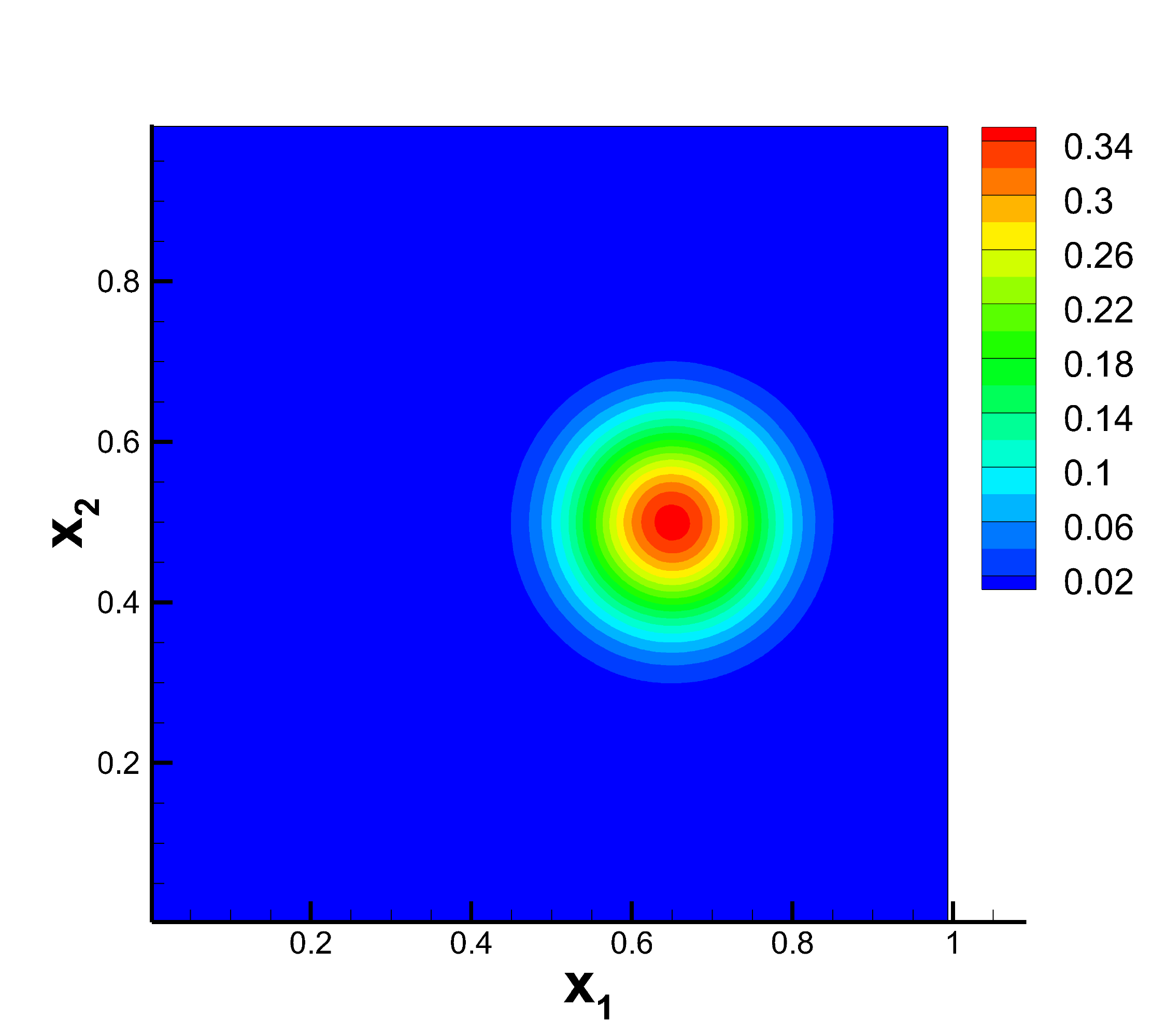}}\\
\subfigure[]{\includegraphics[width=.42\textwidth]{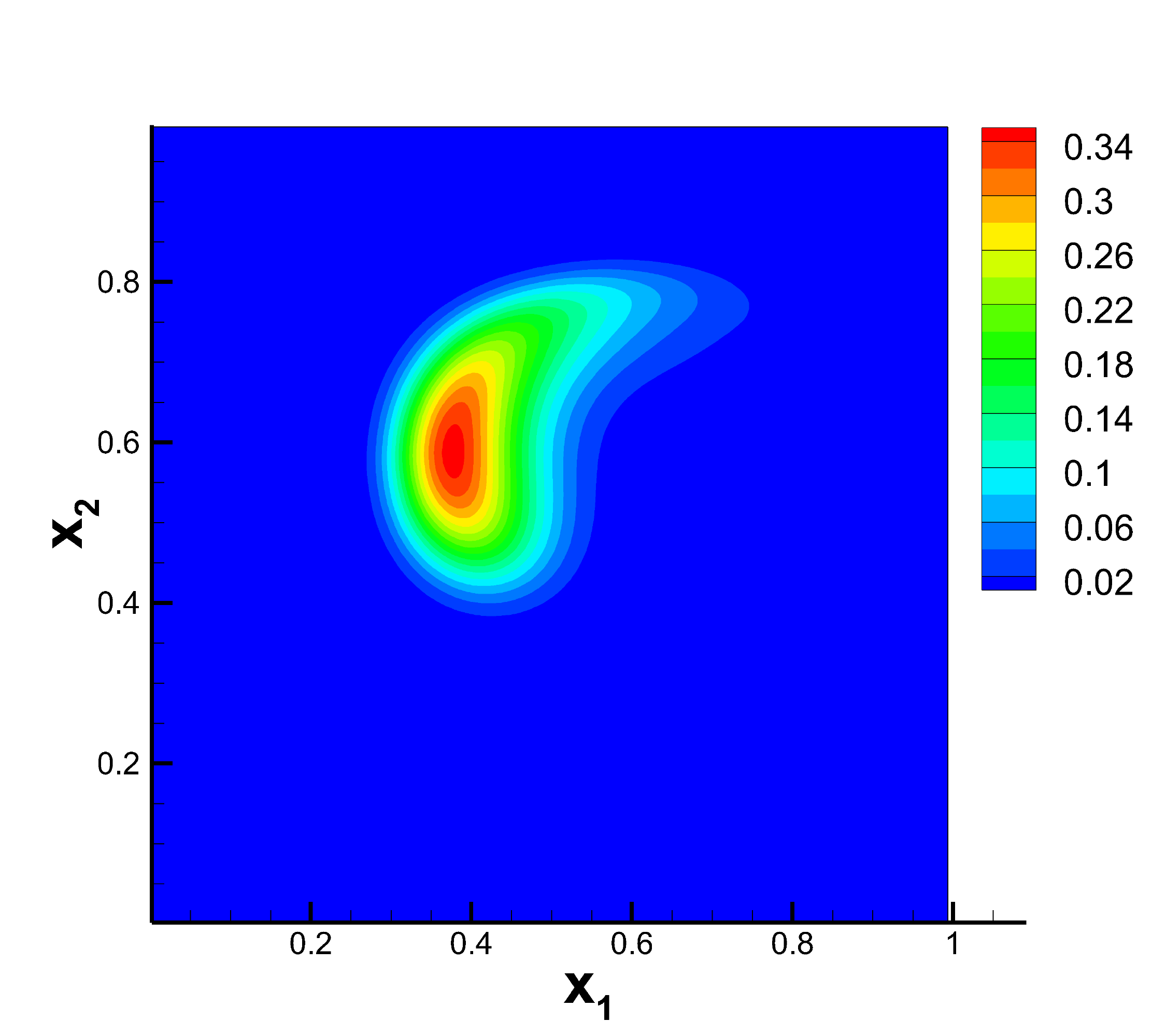}}
\subfigure[]{\includegraphics[width=.42\textwidth]{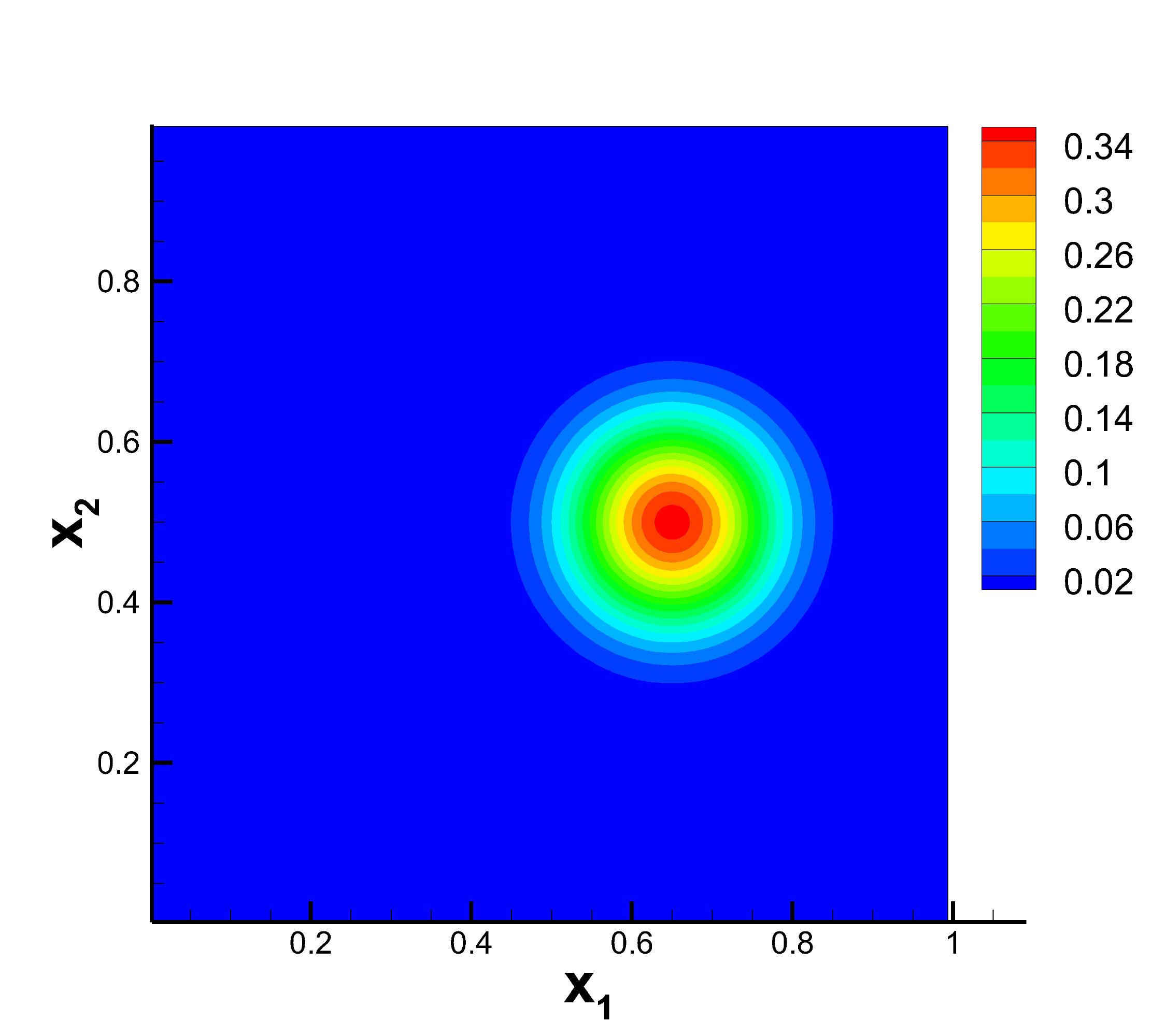}}
\end{center}
  \caption{Example \ref{ex:deformational}. Deformational flow test. The contour plots of the numerical solutions  at $t=T/2$ (a, c, e) and $t=T$ (b, d, f). $k=1$ (a, b), $k=2$ (c, d), and $k=3$ (e, f). $N=7$. }
 \label{fig:defor}
\end{figure}

\section{Kinetic simulations}
\label{sec:kinetic}

In this section, we apply the   sparse grid DG methods to solve kinetic equations. One of the major challenges for deterministic kinetic simulations is the high dimensionality of the underlying system. The unknown probability distribution function depends on the space   and the velocity variables, which translates to numerical computation in six dimensions plus time in real-world settings. It is therefore a good test bed for sparse grid algorithms. In this paper, we   focus on two types of kinetic systems:   the collisionless Vlasov model and the collisional   relaxation model. Both models are well understood in the literature and used as algorithm benchmarks. Our calculation will not involve solutions with discontinuity and numerical tests will be performed focusing on validating accuracy and conservation.


\subsection{The collisionless Vlasov equation} 

In this subsection, we consider the non-dimensionalized single-species nonlinear Vlasov-Amp\`{e}re system (VA) for plasma simulations in the zero-magnetic limit
\begin{eqnarray}
&&f_t + \bv\cdot\nabla_\bx f + \bE(t,\bx)\cdot\nabla_\bv f=0, \label{eq:V}\\
&&\partial_t \bE = -\bJ, \label{eq:A}
\end{eqnarray}
where $f(t,\bx,\bv)$ denotes the probability distribution function of electrons. $\bE(t,\bx)$ is the self-consistent electrostatic field given by Amp\`{e}re's law \eqref{eq:A} and $\bJ(t,\bx)=\int_{\bv} f(t,\bx,\bv) \bv d\bv$ denotes the electron current density. Ions are assumed to form a neutralizing background. 
Note that the discussion here can be easily generalized to other popular Vlasov systems, such as Vlasov-Poisson or Vlasov-Maxwell systems. 
 
The sparse grid DG method for the VA system can be formulated as follows. We solve the Vlasov equation by the algorithm discussed in previous sections,  i.e. we compute $f_h$ according to a similar formulation as in \eqref{eq:DGformulation} treating the Vlasov equation as a variable coefficient transport problem, while the Amp\`{e}re equation can be solved exactly, since $\bJ_h$ is a piecewise polynomial supported on the discrete mesh in $\bx$ and can be obtained exactly. In fact, the calculation of electron current $\bJ_h$ is straightforward since $\bv$ is orthogonal to all bases except those residing on level $0.$ 
The third order TVD-RK scheme is  then used as the time discretization.


A focus of this work is to verify the conservation properties of the scheme. It is well known that the VA system preserves many physical invariants, including the particle number $\int_\bx\int_\bv f(t,\bx,\bv)\,d\bx d\bv,$ total energy $\frac12\int_\bx\int_\bv f(t,\bx,\bv)|\bv|^2\,d\bx d\bv +\frac12 \int_\bx |\bE(t,\bx)|^2\,d\bx,$ and enstrophy $\int_\bx\int_\bv |f(t,\bx,\bv)|^2\,d\bx d\bv.$ Though it is difficult to design a numerical scheme which is able to preserve all the above invariants on the discrete level, tracking these quantities in time provides a good measurement for the performance of numerical schemes.  Traditional DG methods were developed to solve the VA system  \cite{cheng2014energy} and superior performance in conservation of particle number and energy has been observed. As for the sparse grid DG method, we can easily verify similar properties for the semi-discrete scheme. The main reason is that the test function $1$ and $|\bv|^2$ used as the key steps in the proof of numerical conservation will still belong to the space $\hat{\bV}^k_N$ as long as $k \ge 2.$ In fact, they are functions that belong to the coarsest level of mesh.

\begin{thm}[Conservation properties]
	Without the boundary effect, the semi-discrete sparse grid DG scheme   for solving the VA system  conserves the particle number. If $k\ge2$, the scheme also conserves the total energy. The scheme is $L^2$ stable, i.e. $\frac{d}{dt}\int_\bx\int_\bv |f_h(t,\bx,\bv)|^2\,d\bx d\bv  \leq 0.$
\end{thm}  

The proof is  similar to \cite{cheng_vm, cheng2014energy} and is omitted. \hfill\ensuremath{\blacksquare}

Inspired by the theorem above, if one wants to design a DG scheme with particle number and energy conservation, it is enough to choose a basis set that includes $1$ and $|\bv|^2$ on level $0$, while on other levels the bases can be chosen freely according to accuracy consideration. Such adaptivity is of interest to our future studies. In this paper, we consider the following two benchmark test cases in a 1D1V setting.



\begin{itemize}
  \item Landau damping:
  \begin{equation}
  f(0,x,v) = f_M(v)(1+A\cos(kx)),\quad x\in[0,L],\,v\in[-V_c,V_c],
  \end{equation}
  where $A=0.5$, $k=0.5$, $L=4\pi$, $V_c=2\pi$, and $f_M(v)=\frac{1}{\sqrt{2\pi}}e^{-v^2/2}$.
  \item Two-stream instability:
  \begin{equation}
  f(0,x,v) = f_{TS}(v)(1+A\cos(kx)),\quad x\in[0,L],\,v\in[-V_c,V_c],
  \end{equation}
  where $A=0.05$, $k=0.5$, $L=4\pi$, $V_c=2\pi$, and $f_{TS}(v)=\frac{1}{\sqrt{2\pi}}v^2e^{-v^2/2}$.
\end{itemize}
The periodic boundary condition is imposed in $x$-space. As a standard practice, the computational domain in $v$ is truncated to $  [-V_c,V_c]$, where $V_c$ is a constant chosen large enough  to impose zero boundary condition in the $v$-direction  $f_h(t,x,\pm V_c) = 0.$


We first perform the accuracy test. Here we utilize the time reversibility of the VA system, i.e. if we let $f(0,x,v)$ be the initial condition and $f(T,x,v)$ be the solution of  the VA system at $t=T$. When we reverse the velocity field of the solution, yielding $f(T,x,-v)$, and evolve the VA system again to $t=2T$, we would recover $f(0,x,-v)$, which is the initial condition with reverse velocity field. In our simulation, we use the sparse grid DG scheme to compute the solution to $T=1$ and then back to $T=2$, and compare it with the initial condition.
The $L^2$ errors and orders of accuracy are reported in Table \ref{table:va}. For both tests, we can observe slight order reduction from the optimal accuracy, which is similar to previous examples. The loss of accuracy for  $\hat{V}_9^3$ for two-stream instability is due to the domain cut-off in the  $v$-space, which causes   local truncation of about $10^{-8}$. This error can be reduced by taking a larger $V_c$. 

\begin{table}
	\caption{$L^2$ errors and orders of accuracy for the Vlasov-Amp\`{e}re system.
		 $N$ is the number of mesh levels, $h_N$ is the size of the smallest mesh in each direction, $k$ is the polynomial order, $d$ is the dimension. DOF denotes the degrees of freedom of the   sparse approximation space $\hat{V}^k_N$.  $L^2$ order is calculated with respect to  $h_N$. $d=2$. \label{table:va}}
	\centering
	\begin{tabular}{|c|c |c c c|   c c c|  }
		\hline
		$N$ & $h_N$& DOF&$L^2$ error & order&  DOF&$L^2$ error & order \\
		\hline
		& &        \multicolumn{3}{c|}{ Landau damping}          &     \multicolumn{3}{c|}{Two-stream instability}  \\
		\hline
		& &        \multicolumn{3}{c|}{$ k=1$}          &     \multicolumn{3}{c|}{$ k=1$}  \\
		\hline
		5   &   $4\pi$/32 &  448     & 1.44E-01	&	--  &  	448& 2.77E-02	&	--  \\	
		6	&   $4\pi$/64 &	1024&   5.71E-02	&	1.34   &1024&	7.37E-03	&	1.91   \\
		7	&	$4\pi$/128 &  2304& 1.17E-02	&	2.28  & 2304   &  2.12E-03	&	1.80 
		\\
		8	&	$4\pi$/256 &  5120   & 3.07E-03	&	1.94 & 	5120& 5.89E-04	&	1.85
		\\
		9	&	$4\pi$/512 & 11264& 8.01E-04	&	1.94 & 11264&1.52E-04	&	1.96
		\\	
		\hline 
		& &        \multicolumn{3}{c|}{$ k=2$ }          &     \multicolumn{3}{c|}{$ k=2$ }  \\
		\hline
		5   &  $4\pi$/32 & 1008 &  1.03E-02	&	--       &	1008 &   2.58E-03	&	--
		\\ 
		6	& $4\pi$/64 &	2304&	  3.07E-03	&	1.75 &	2304& 3.89E-04	&	2.73
		\\
		7	& $4\pi$/128 &	5184&	  4.62E-04	&	2.73 &	5184 & 6.13E-05	&	2.67
		\\
		8	& $4\pi$/256 &	11520&   1.09E-04	&	2.08 &	11520 &9.66E-06	&	2.67
		\\
		9	& $4\pi$/512 &	25344&   1.86E-05	&	2.55 &25344 & 1.59E-06 &		2.60
		
		\\
		\hline 
		& &        \multicolumn{3}{c|}{$ k=3 $}          &     \multicolumn{3}{c|}{$ k=3$ }  \\
		\hline
		5	& $4\pi$/32 &	1792&	1.95E-03	&	-- & 1792& 1.52E-04	&	--
		\\
		6	& $4\pi$/64 &	4096&   4.26E-04	&	2.19 &4096&1.15E-05	&	3.72
		\\
		7	& $4\pi$/128 &	9216&	3.54E-05	&	3.59 & 9216&	8.82E-07	&	3.71
		\\
		8	& $4\pi$/256 &	20480 & 4.44E-06	&	3.00  & 	20480 &5.89E-08	&	3.90
		\\	
		9	& $4\pi$/512 &	45056 & 2.65E-07	&	4.07  & 45056 &  3.56E-08	&	0.73
		\\
		\hline
		
	\end{tabular}
\end{table}


In Figures \ref{fig:evo_lan}-\ref{fig:evo_two}, we plot evolution of the relative error in the total particle number, total energy, and enstrophy, and evolution of error in momentum for the proposed method with sparse approximation space $\hat{\bV}_8^3$. For both cases, it is observed that the momentum  are conserved up to machine error. The total particle number ($f_{\mathbf{1},\bzer}^\bzer$ in the code) is conserved up to the truncation errors at the boundary. Also note that even though the total energy is not  conserved, the relative error is still on a quite small scale ($10^{-7}$). This is due to the boundary effects and the errors from the Runge-Kutta schemes. The visible decay in enstrophy is expected because of  the dissipative upwind flux. 
\begin{figure}[htp]
\begin{center}
\subfigure[]{\includegraphics[width=.42\textwidth]{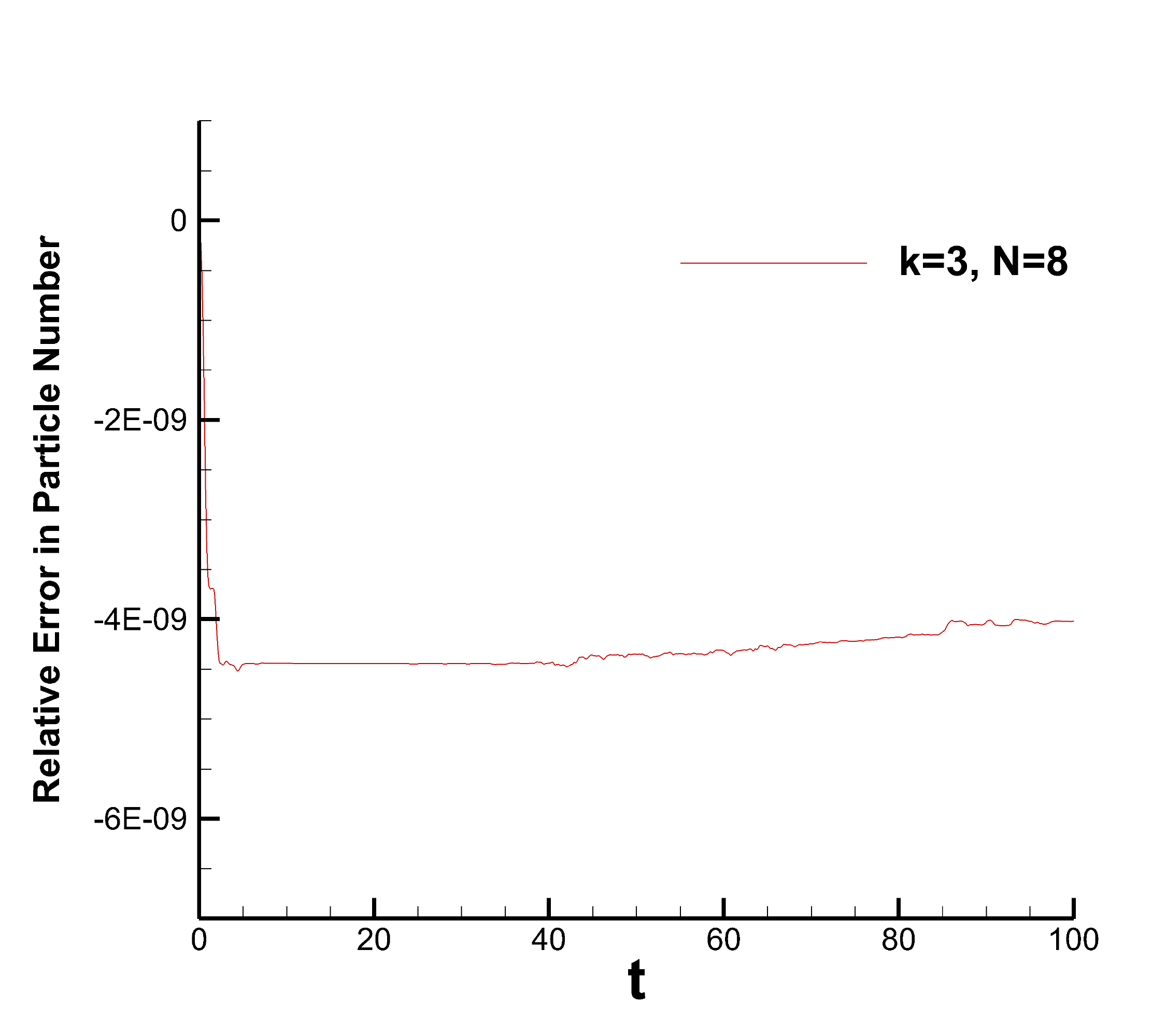}}
\subfigure[]{\includegraphics[width=.42\textwidth]{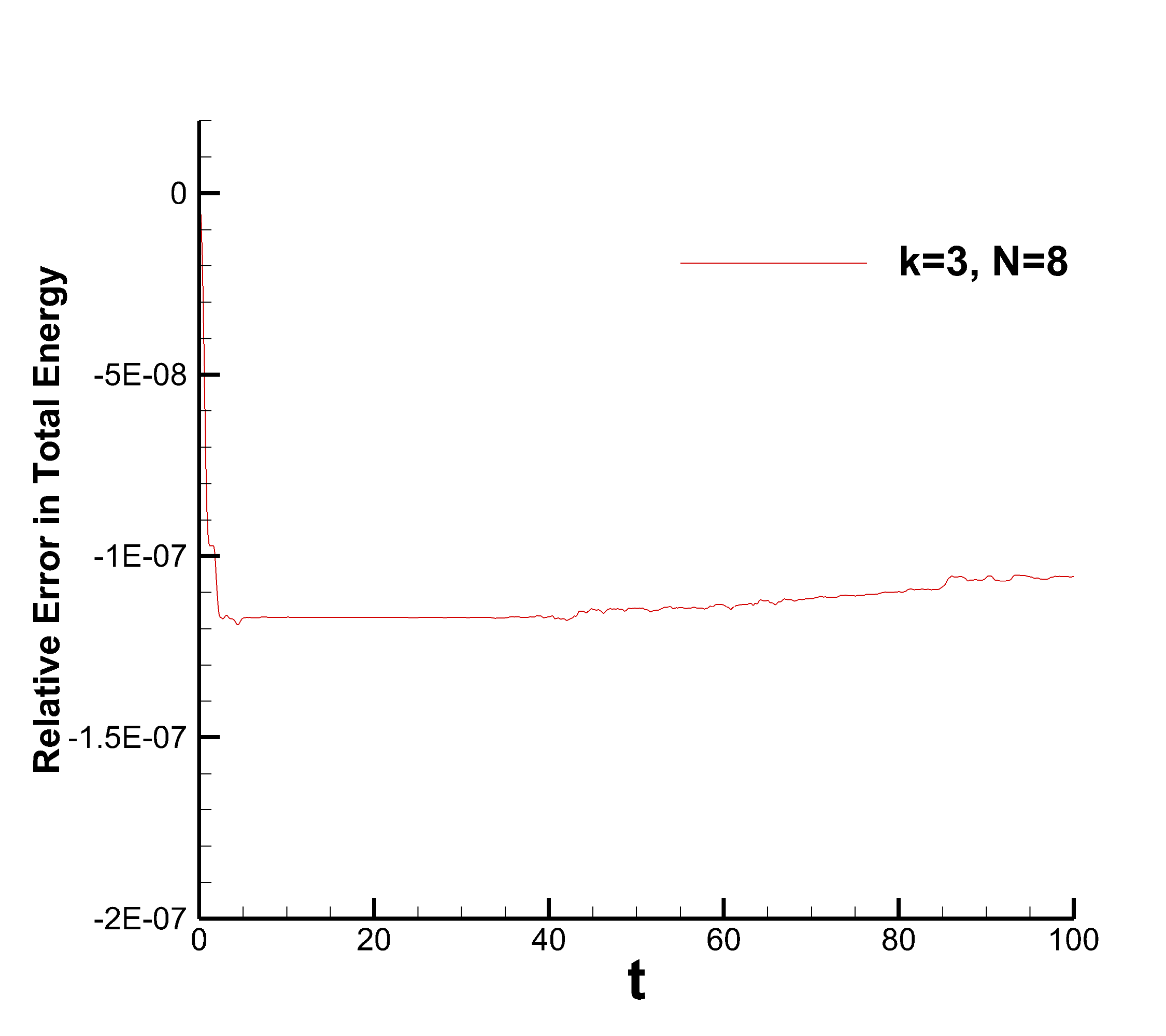}}\\
\subfigure[]{\includegraphics[width=.42\textwidth]{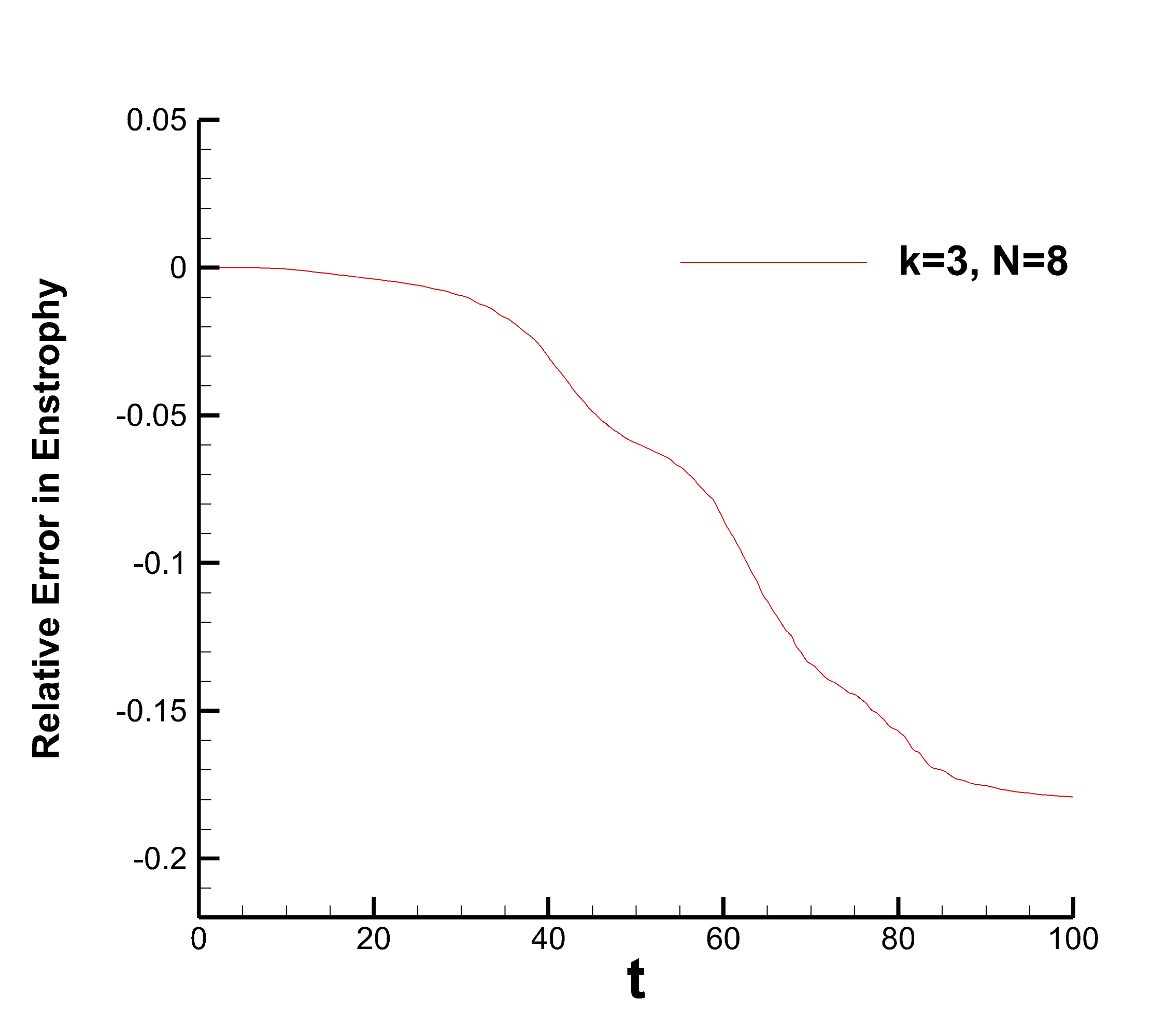}}
\subfigure[]{\includegraphics[width=.42\textwidth]{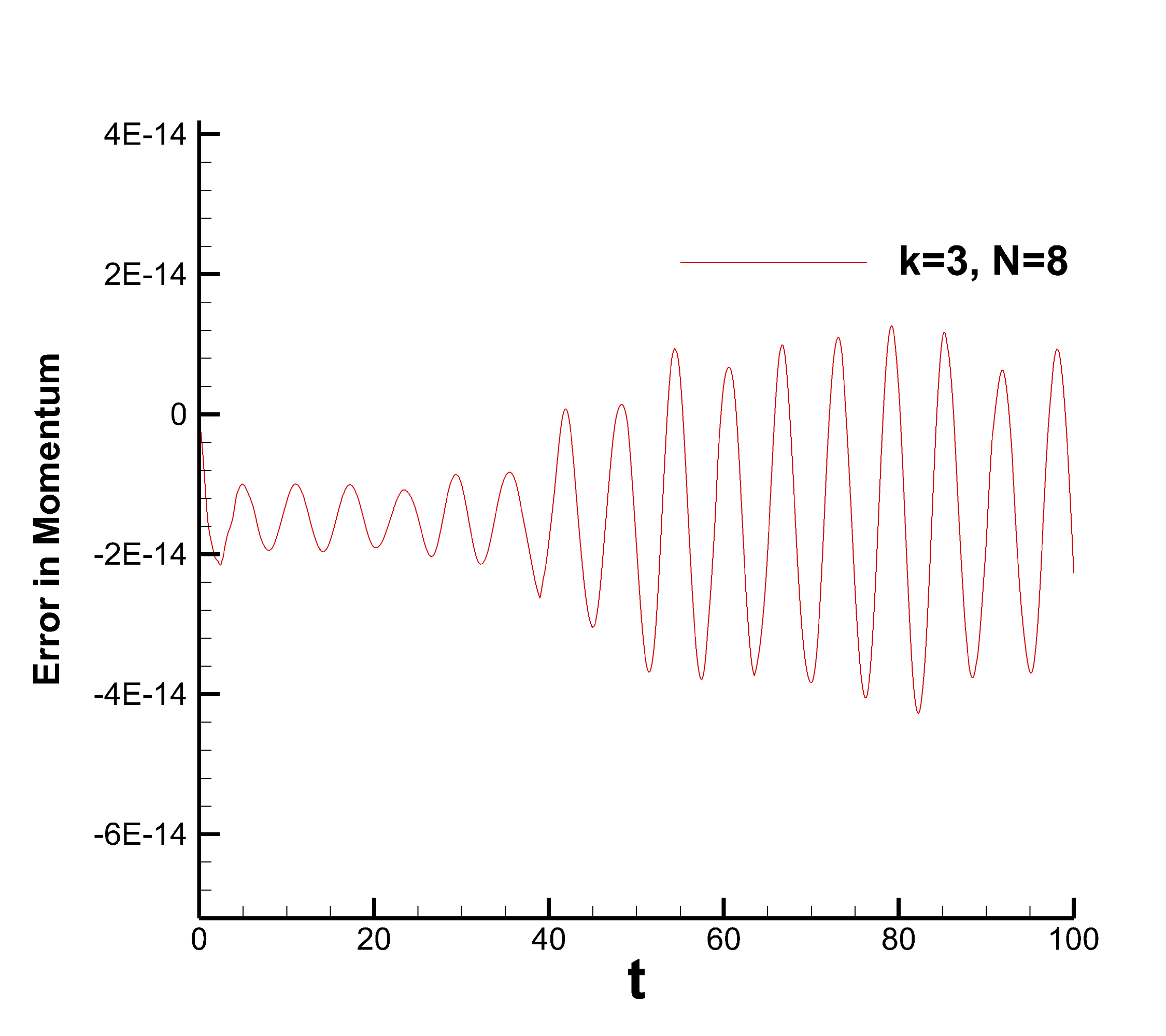}}\\
\end{center}
  \caption{Landau damping. Evolution of the relative errors in total particle number (a), total energy (b), entrophy (c), and evolution of error in momentum (d). $k=3$, $N=8$. }
 \label{fig:evo_lan}
\end{figure}

\begin{figure}[htp]
\begin{center}
\subfigure[]{\includegraphics[width=.42\textwidth]{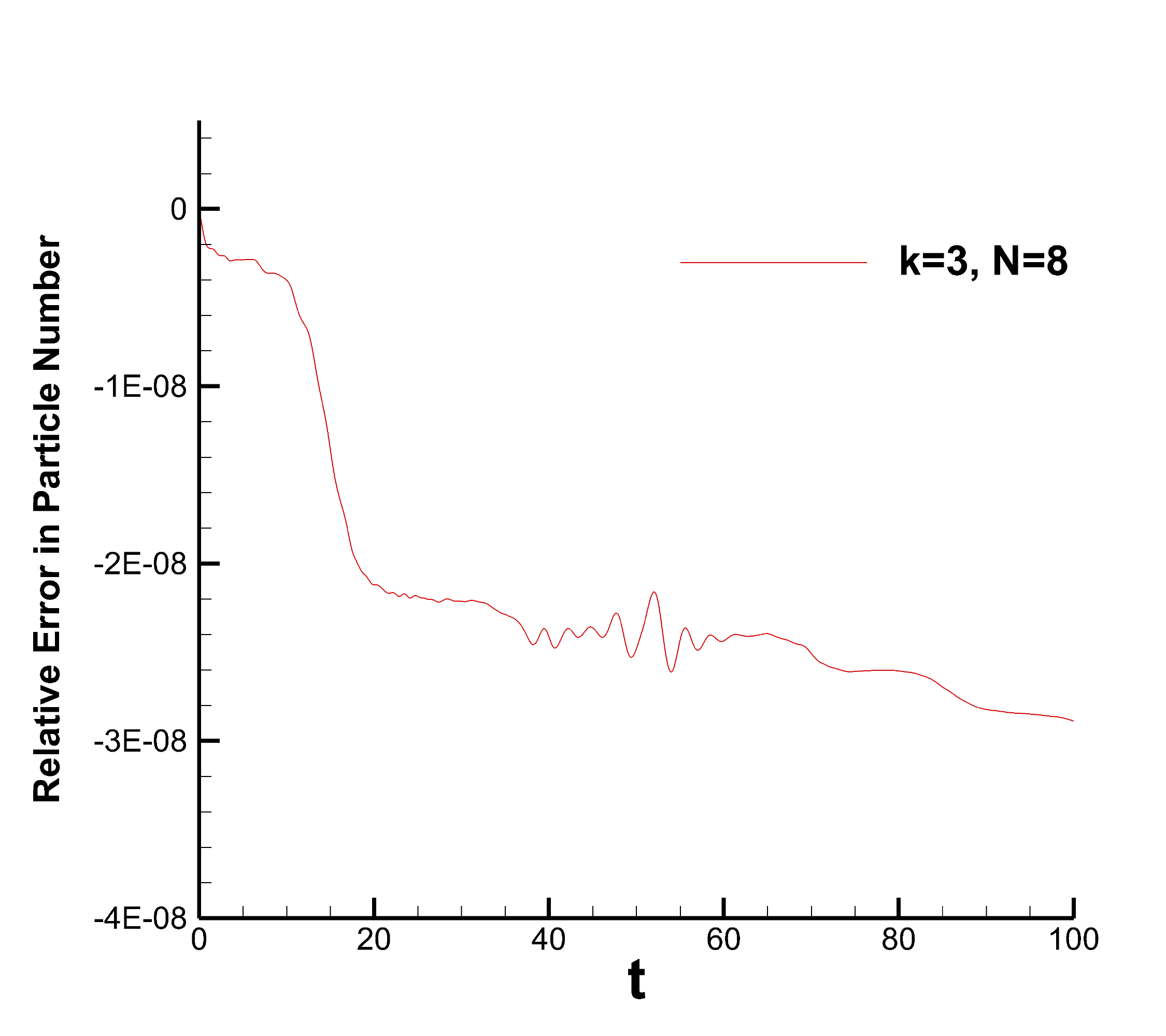}}
\subfigure[]{\includegraphics[width=.42\textwidth]{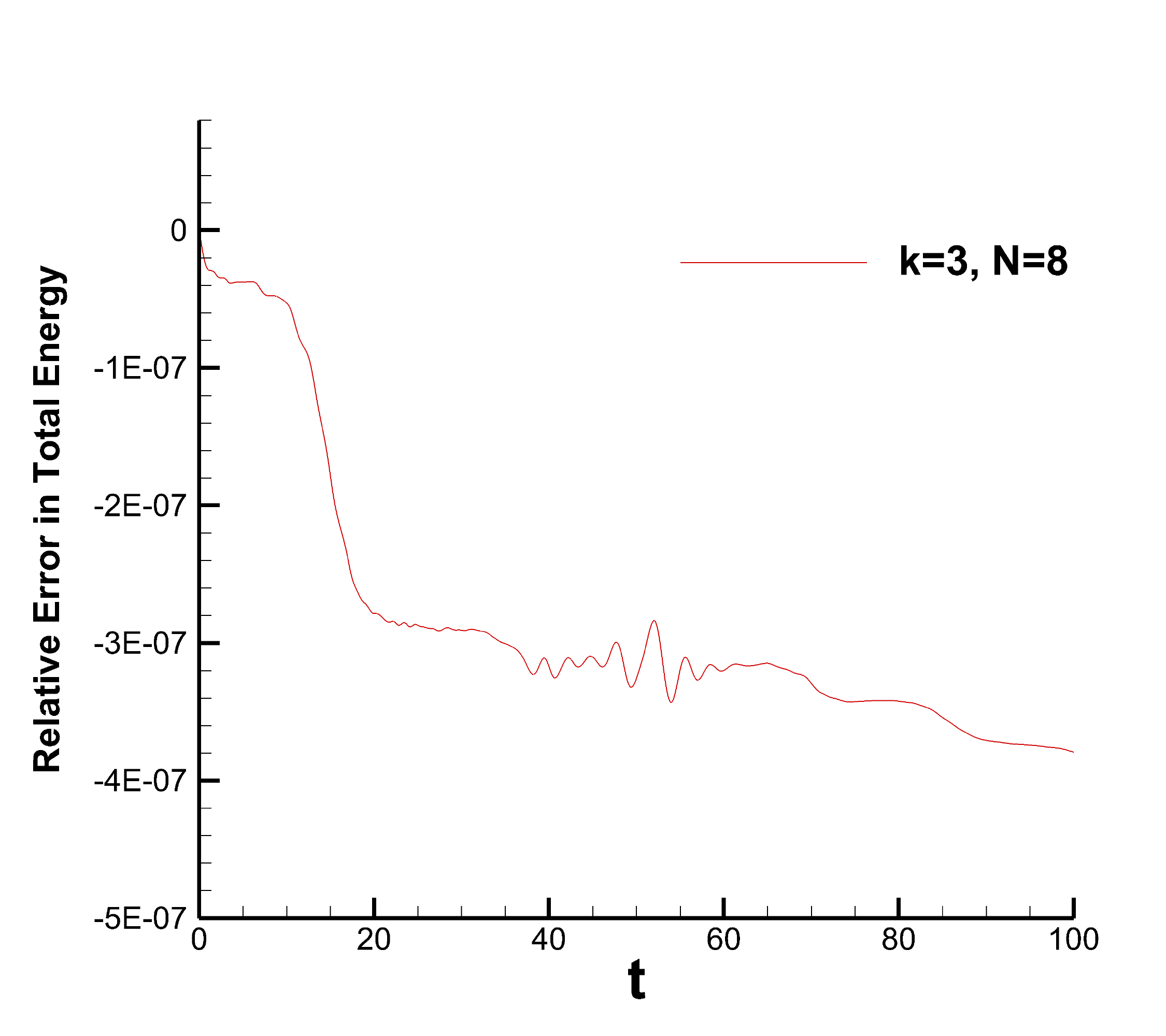}}\\
\subfigure[]{\includegraphics[width=.42\textwidth]{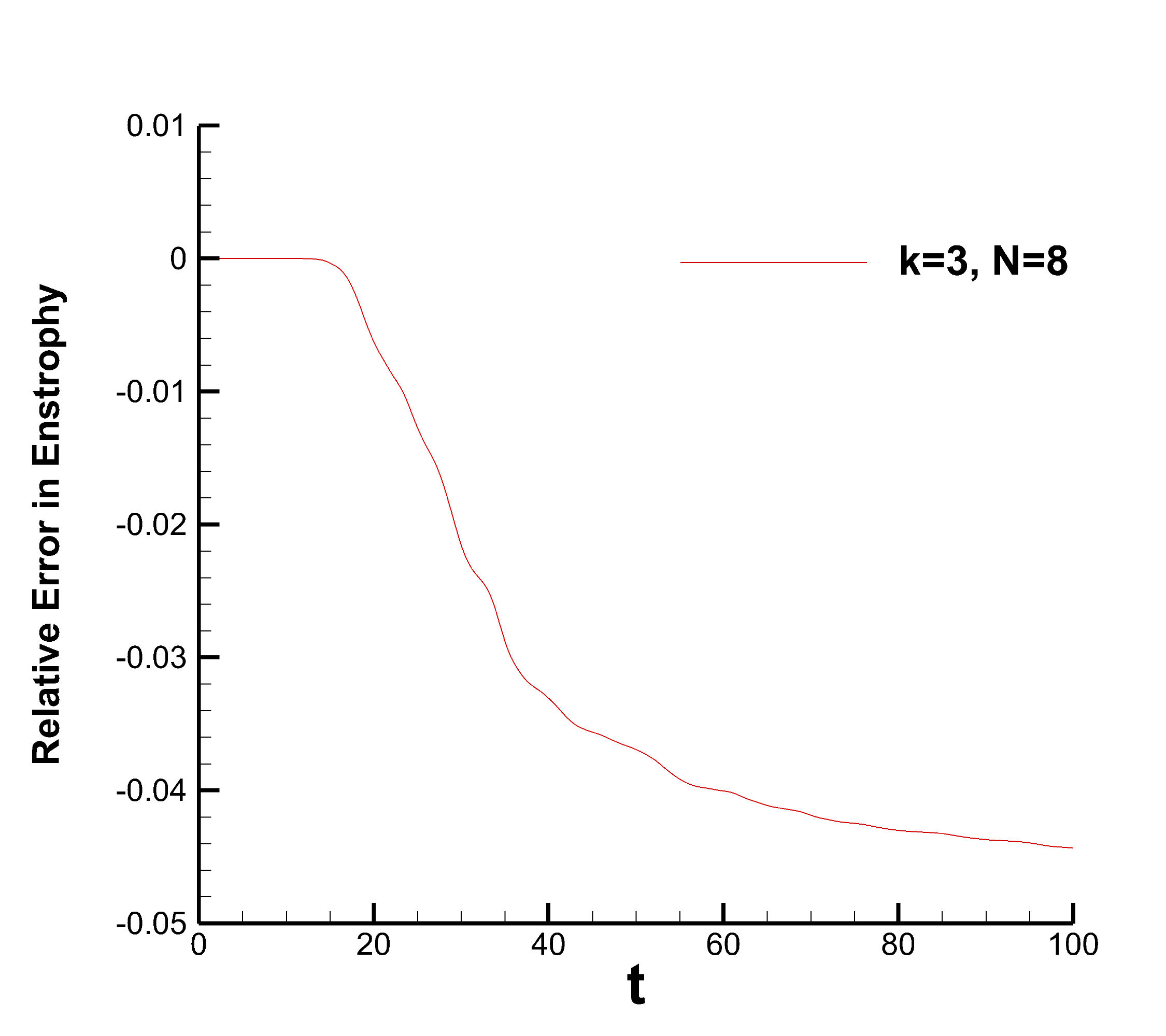}}
\subfigure[]{\includegraphics[width=.42\textwidth]{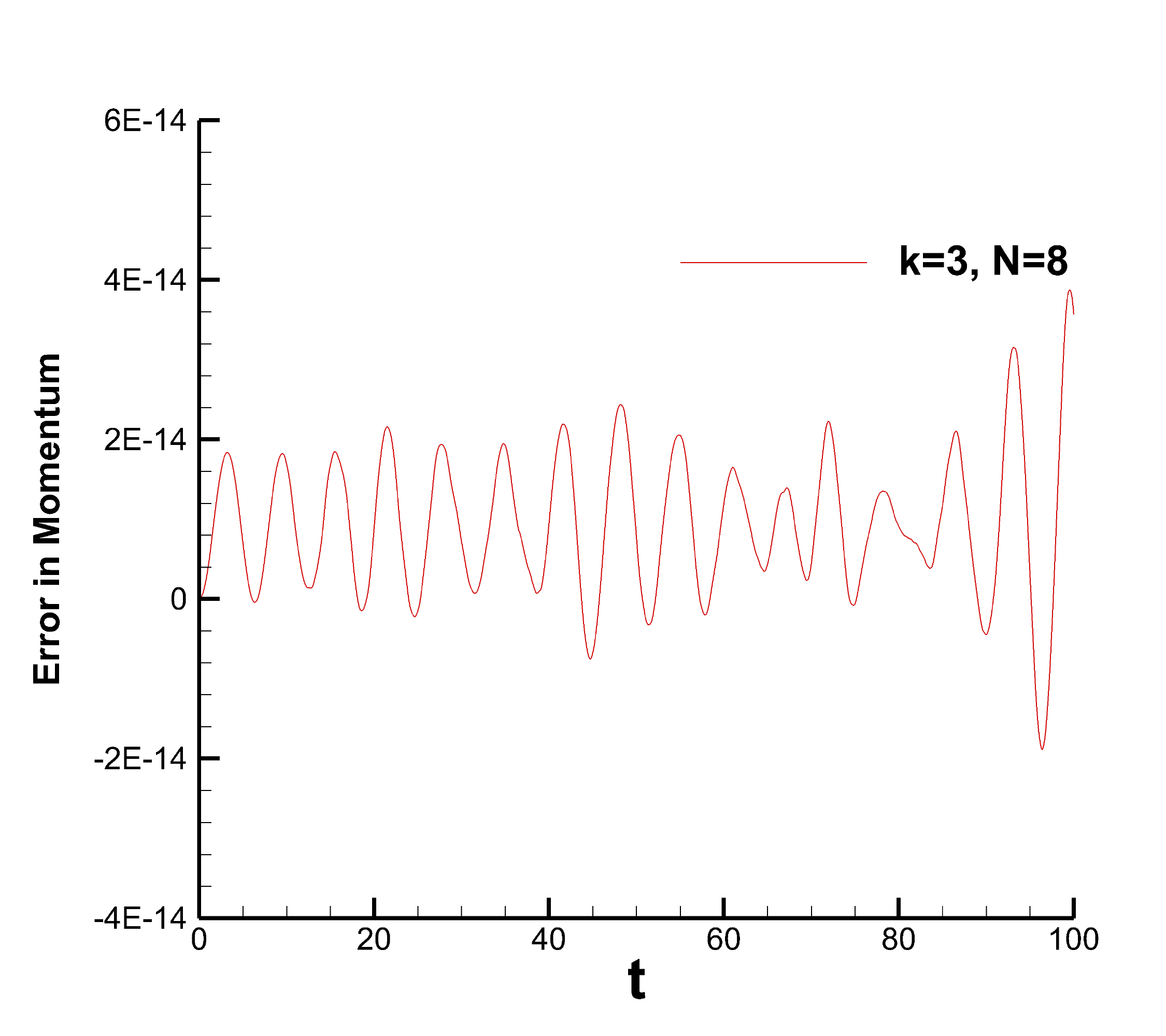}}\\
\end{center}
  \caption{Two-stream instability. Evolution of the relative errors in total particle number (a), total energy (b), entrophy (c), and evolution of error in momentum (d). $k=3$, $N=8$. }
 \label{fig:evo_two}
\end{figure}

We further present some numerical data to benchmark the proposed scheme. We first consider the Log Fourier modes of the electric field $E(t,x)$ as  functions of time, which are defined as 
$$\log FM_n(t)=\log_{10}\left(\frac1L\sqrt{\left| \int_0^L E(t,x)\,\sin(knx)\,dx\right|^2+\left| \int_0^L E(t,x)\,\cos(knx)\,dx\right|^2}\right).$$ 
In Figures \ref{fig:fou_lan}-\ref{fig:fou_two}, we show the time evolution of the first four Log Fourier modes when simulating Landau damping and two-stream instability, respectively.  The sparse approximation space $\hat{\bV}_8^3$ is used in our computation.
The results   agree with other calculations in the literature. 
In Figures \ref{fig:con_lan}-\ref{fig:con_two}, we present the phase space contour plots at several instances of time for Landau damping and two stream instability  computed   with space $\hat{\bV}_8^3$. Note that the number of degrees of freedom of used sparse approximation space is relatively small (which is 20480), yet those numerical results agree with the benchmarks in the literature. 


\begin{figure}[htp]
\begin{center}
\subfigure[]{\includegraphics[width=.42\textwidth]{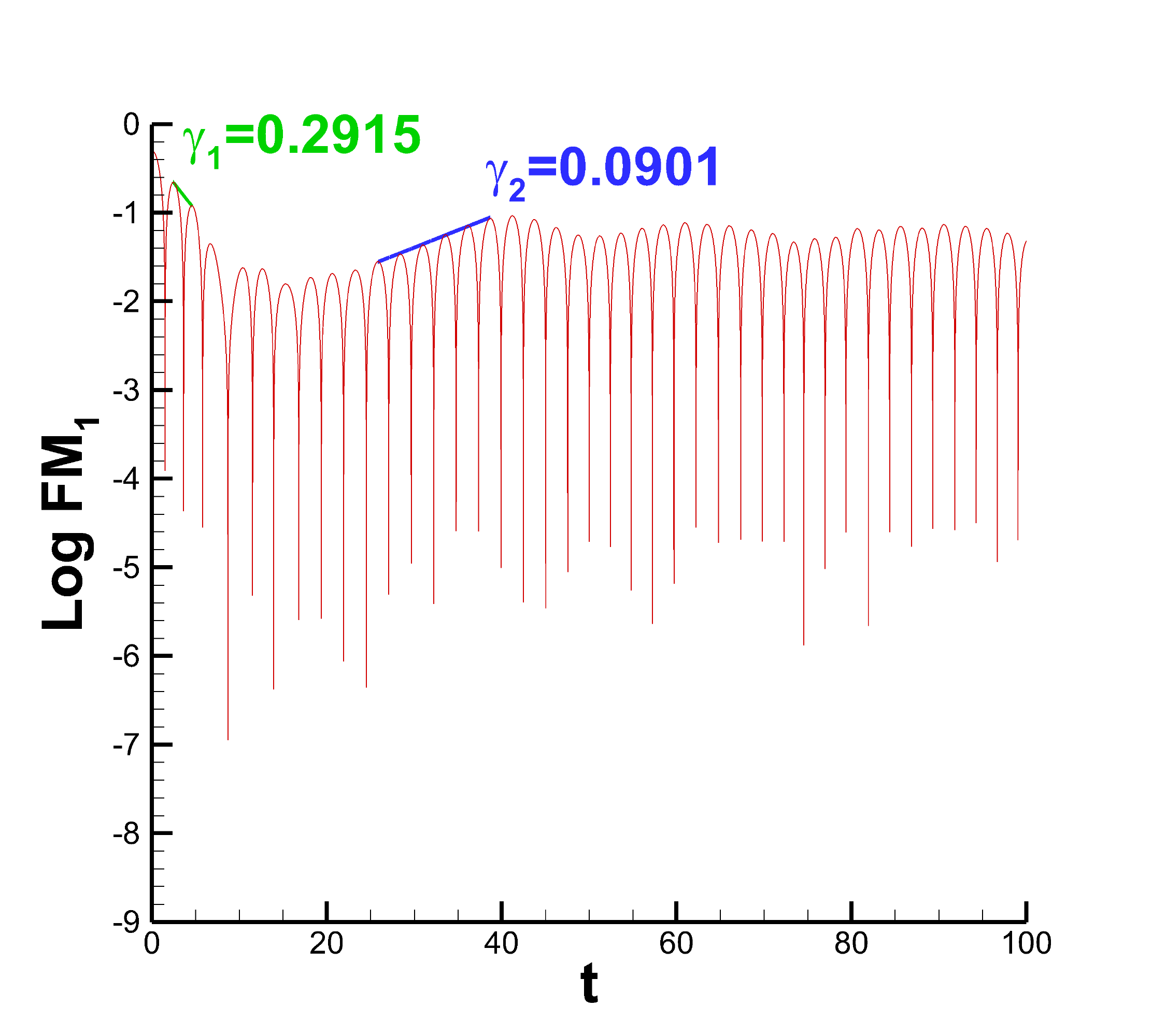}}
\subfigure[]{\includegraphics[width=.42\textwidth]{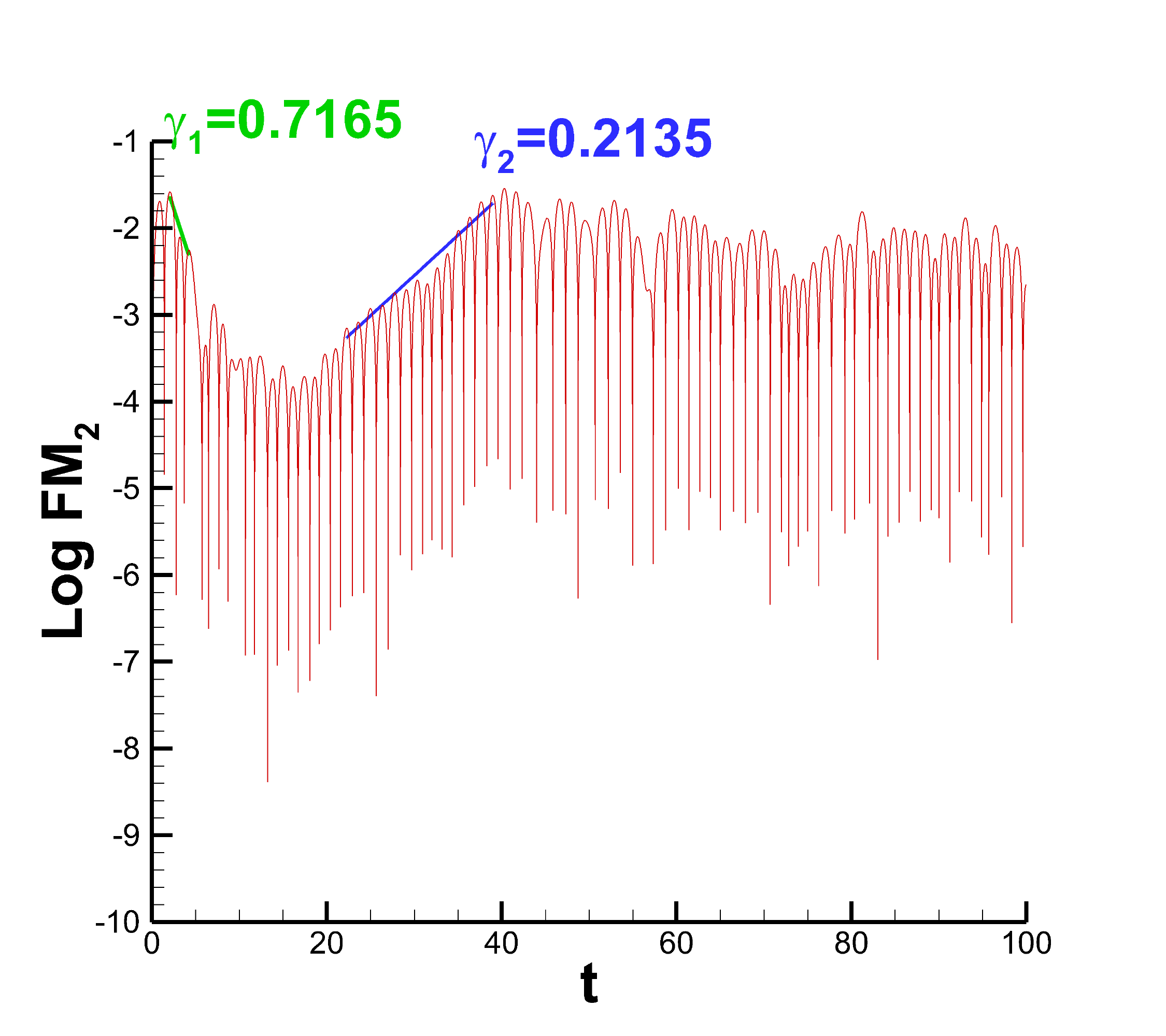}}\\
\subfigure[]{\includegraphics[width=.42\textwidth]{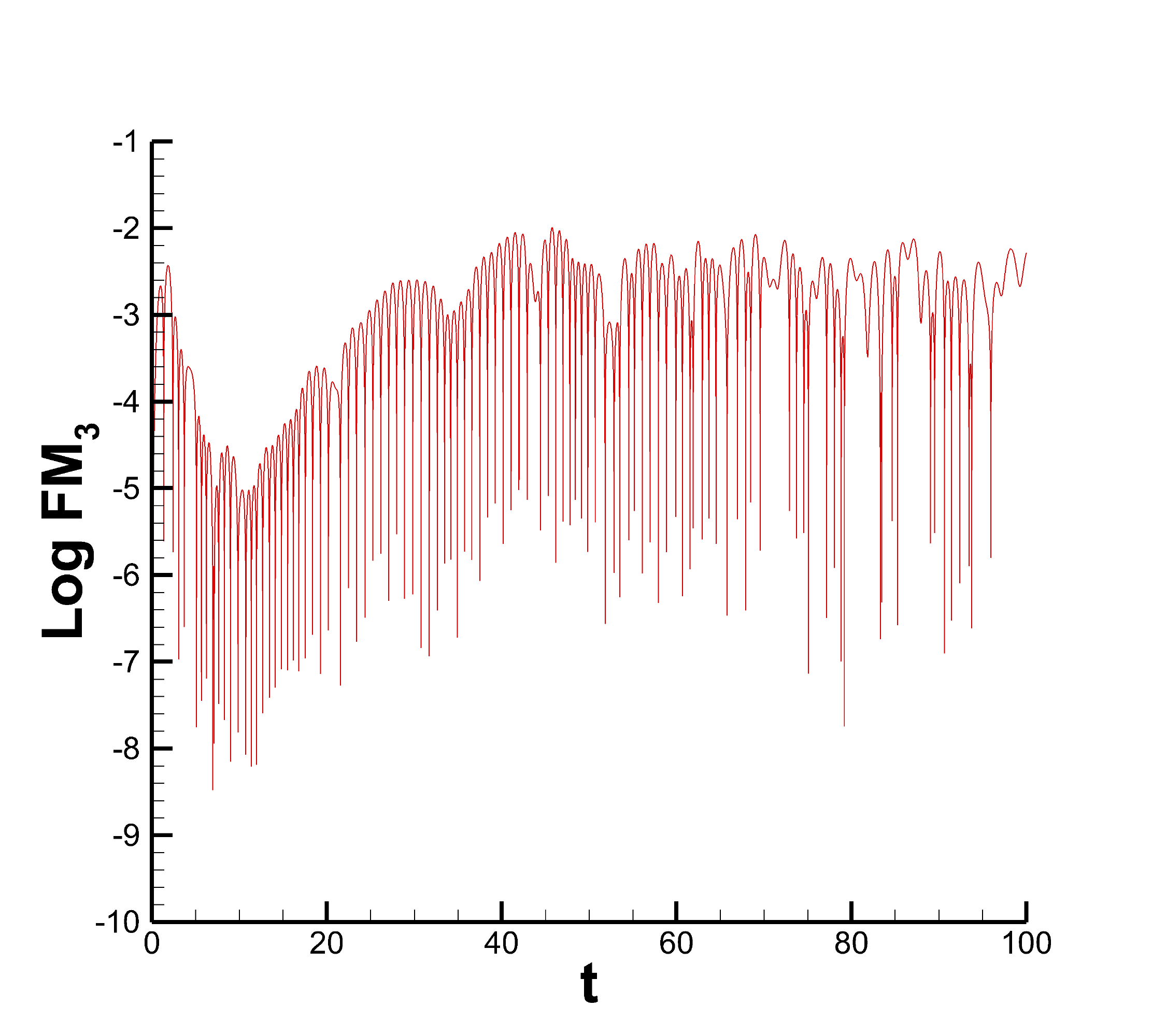}}
\subfigure[]{\includegraphics[width=.42\textwidth]{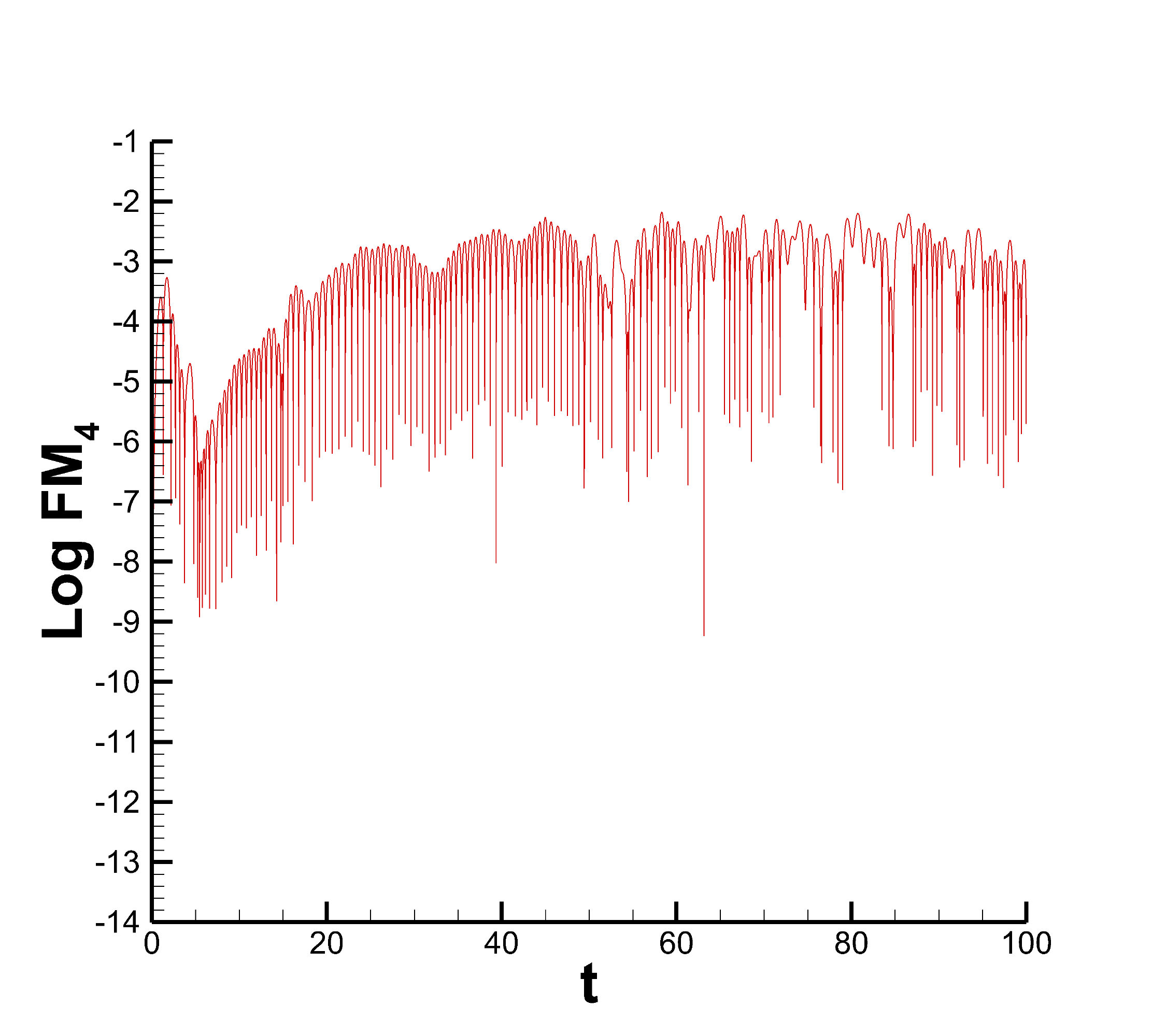}}\\
\end{center}
  \caption{The first four log Fourier modes of Landau damping. $k=3$, $N=8$.}
 \label{fig:fou_lan}
\end{figure}

\begin{figure}[htp]
\begin{center}
\includegraphics[width=.5\textwidth]{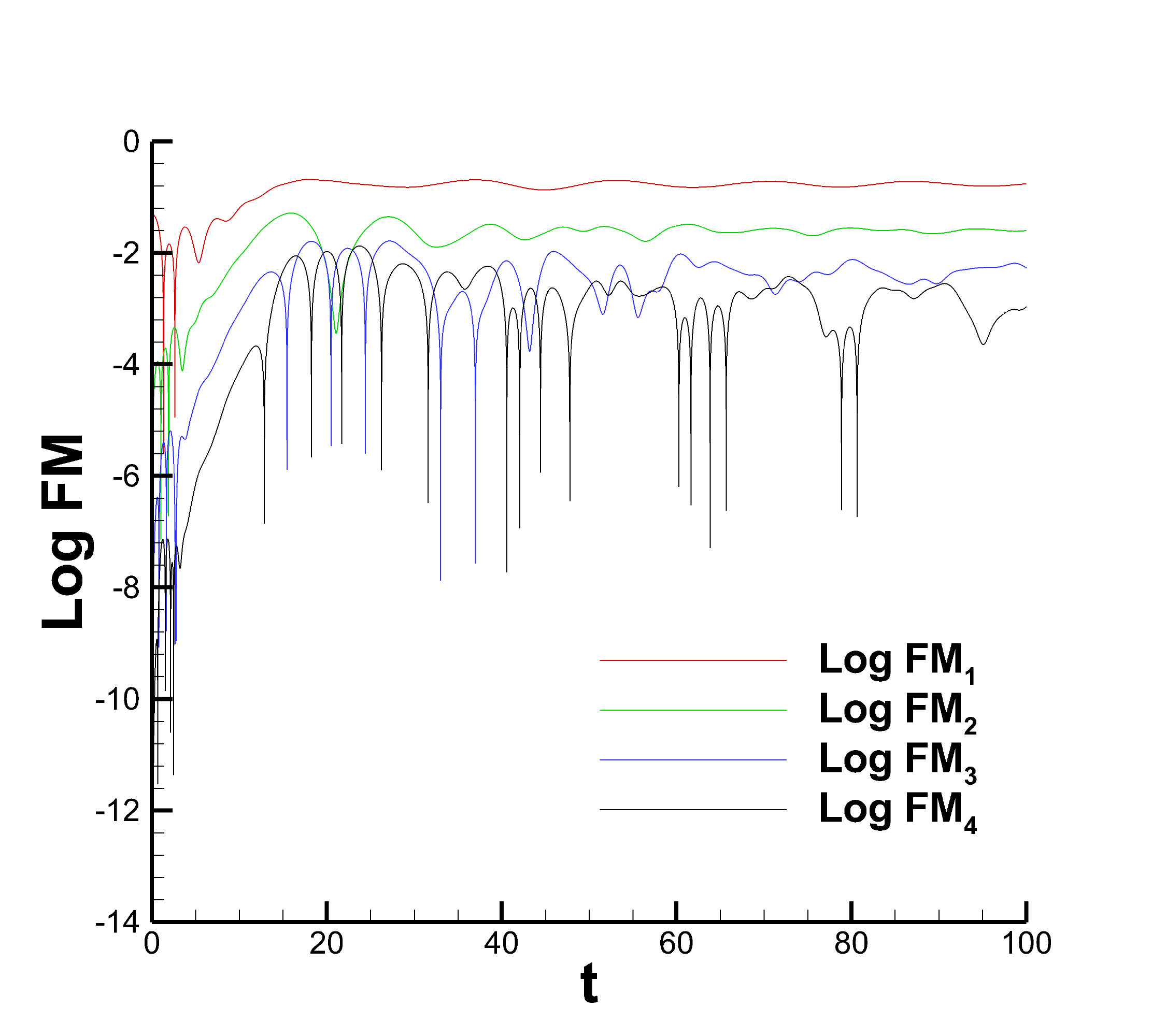}
\end{center}
  \caption{The first four log Fourier modes of two-stream instability. $k=3$, $N=8$.}
 \label{fig:fou_two}
\end{figure}

\begin{figure}[htp]
	\begin{center}
		\subfigure[]{\includegraphics[width=.42\textwidth]{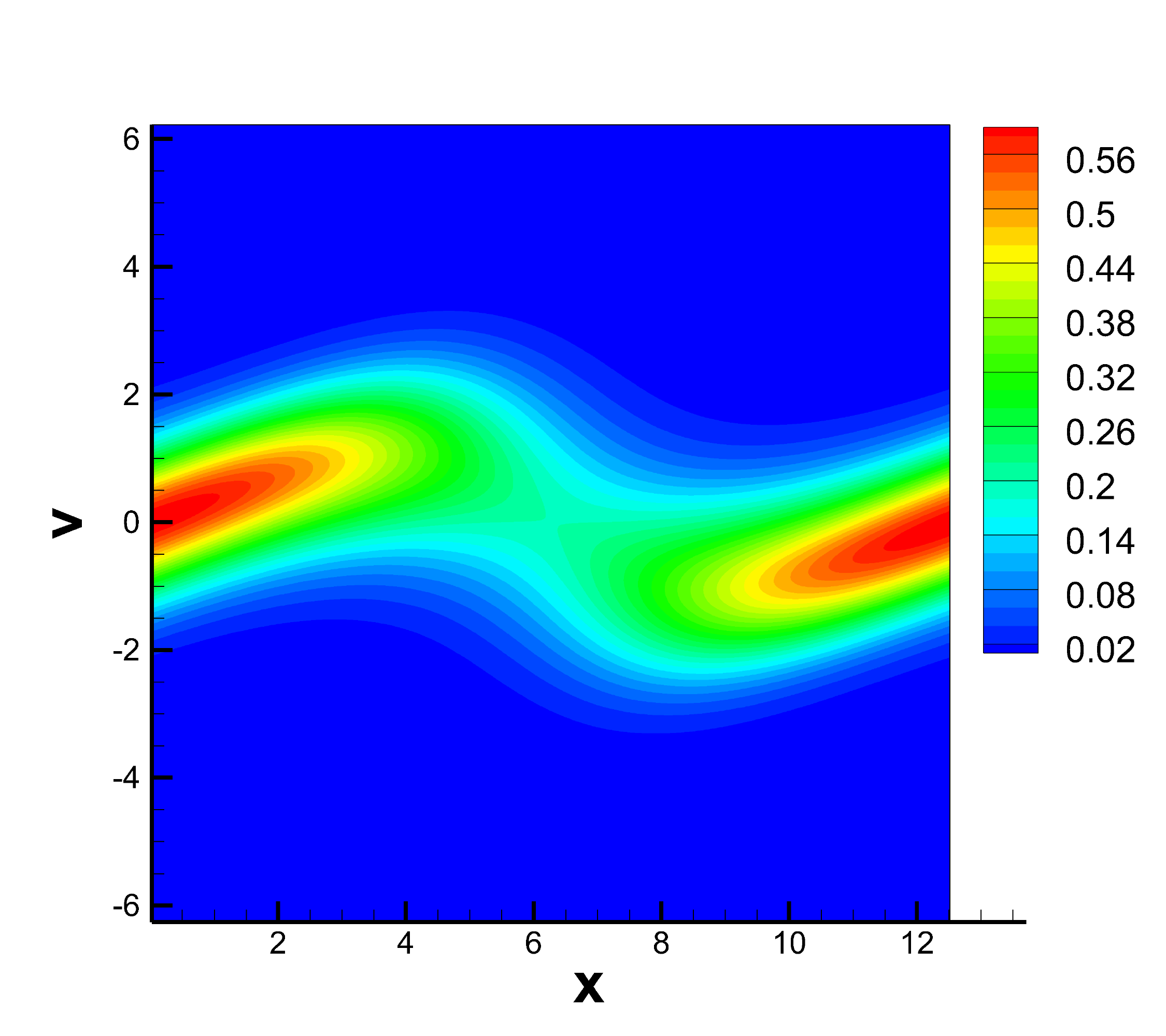}}
		\subfigure[]{\includegraphics[width=.42\textwidth]{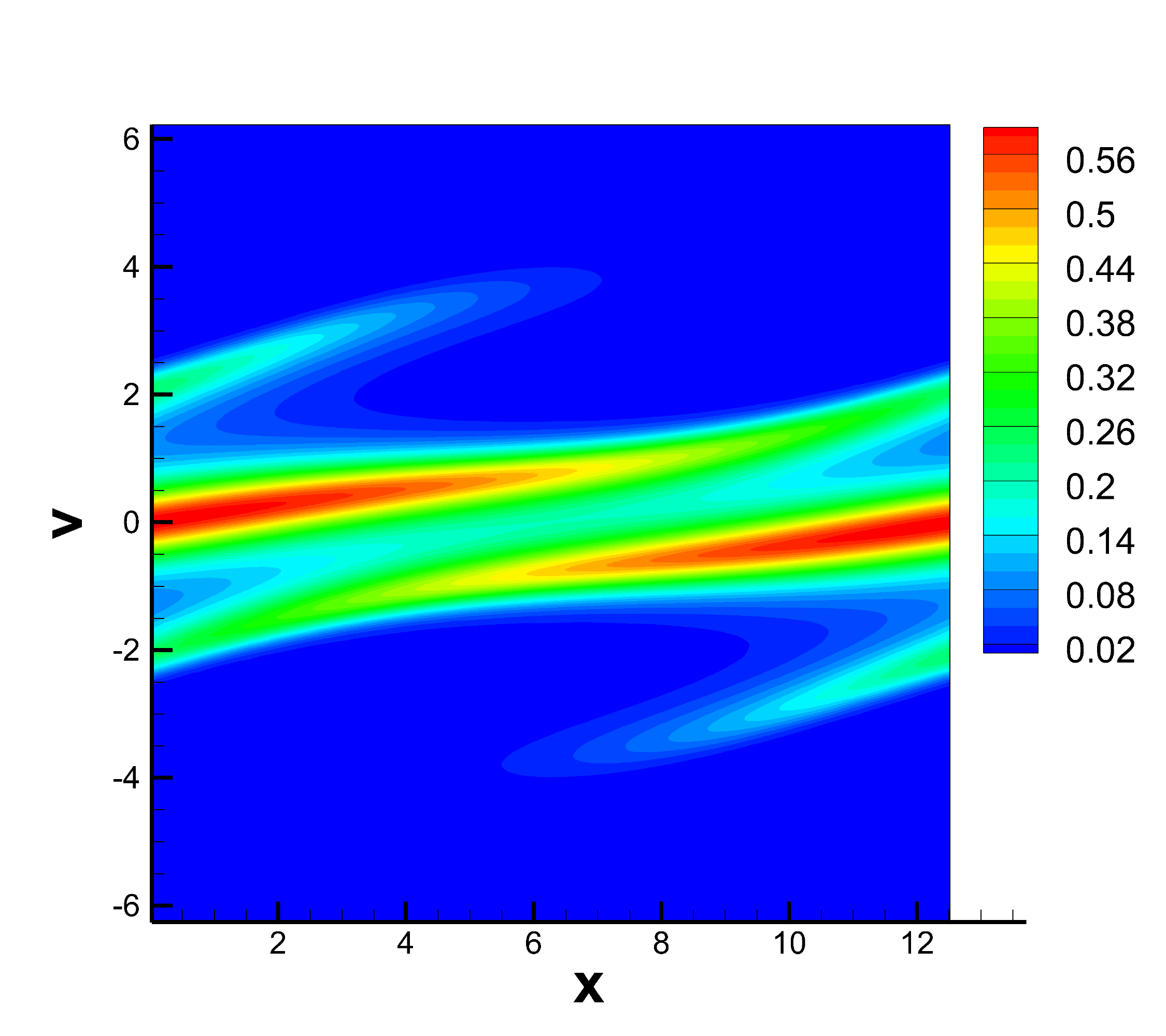}}\\
		\subfigure[]{\includegraphics[width=.42\textwidth]{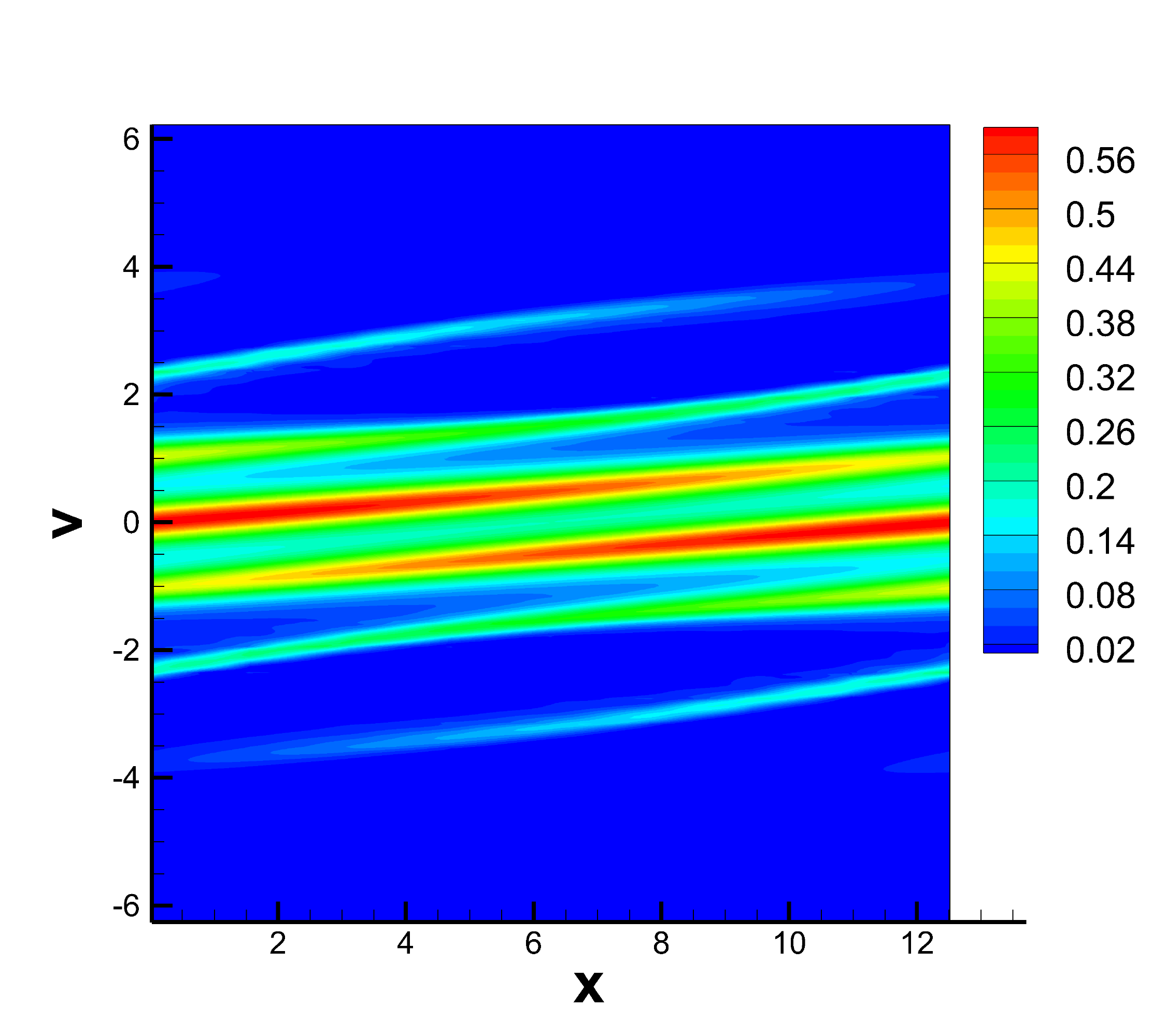}}
		\subfigure[]{\includegraphics[width=.42\textwidth]{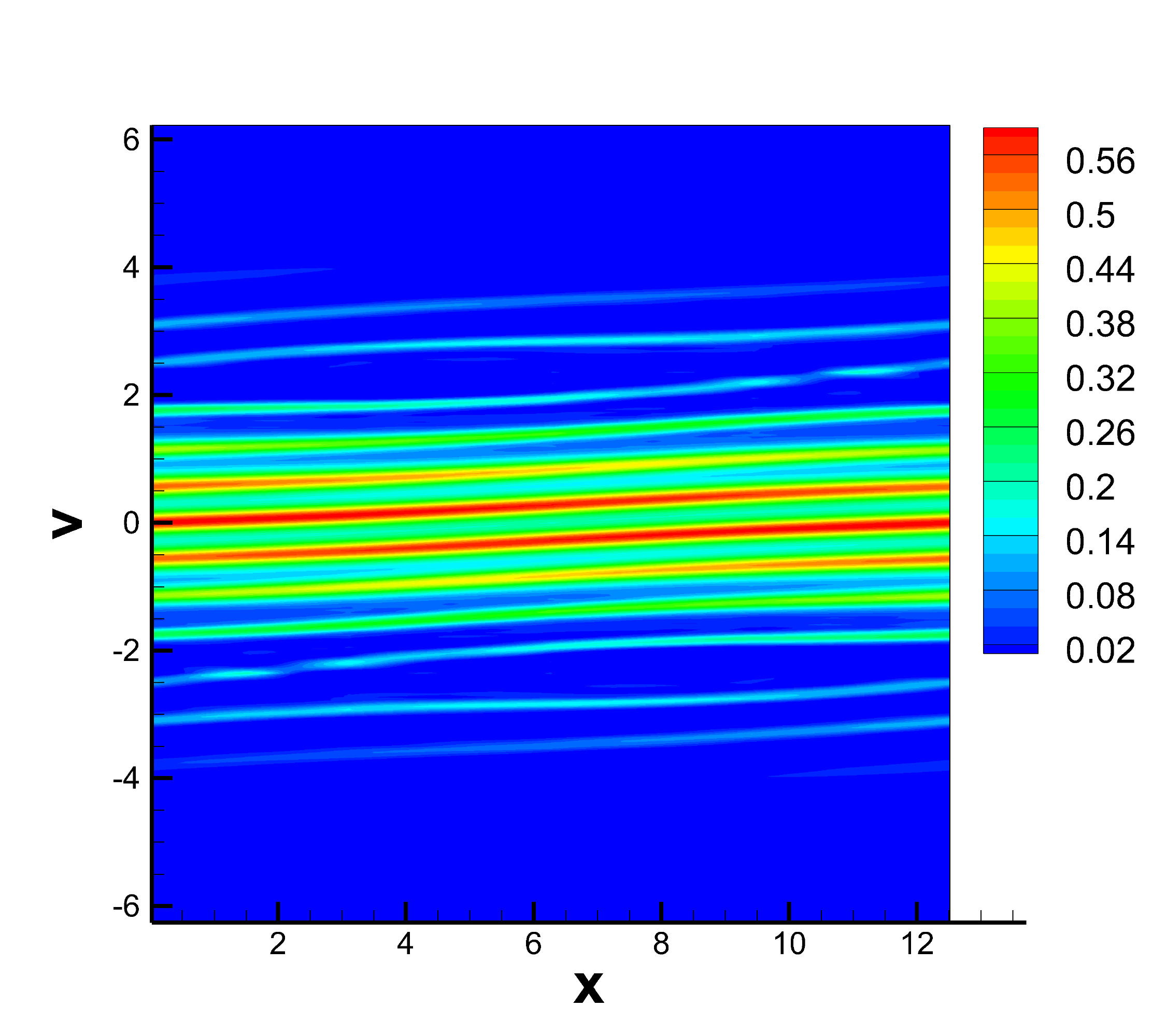}}\\
	\end{center}
	\caption{Landau damping. Phase space contour plots at   $t=1$ (a),  $t=5$ (b), $t=10$ (c), and $t=20$ (d). $k=3$, $N=8$.}
	\label{fig:con_lan}
\end{figure}

\begin{figure}[htp]
\begin{center}
\subfigure[]{\includegraphics[width=.42\textwidth]{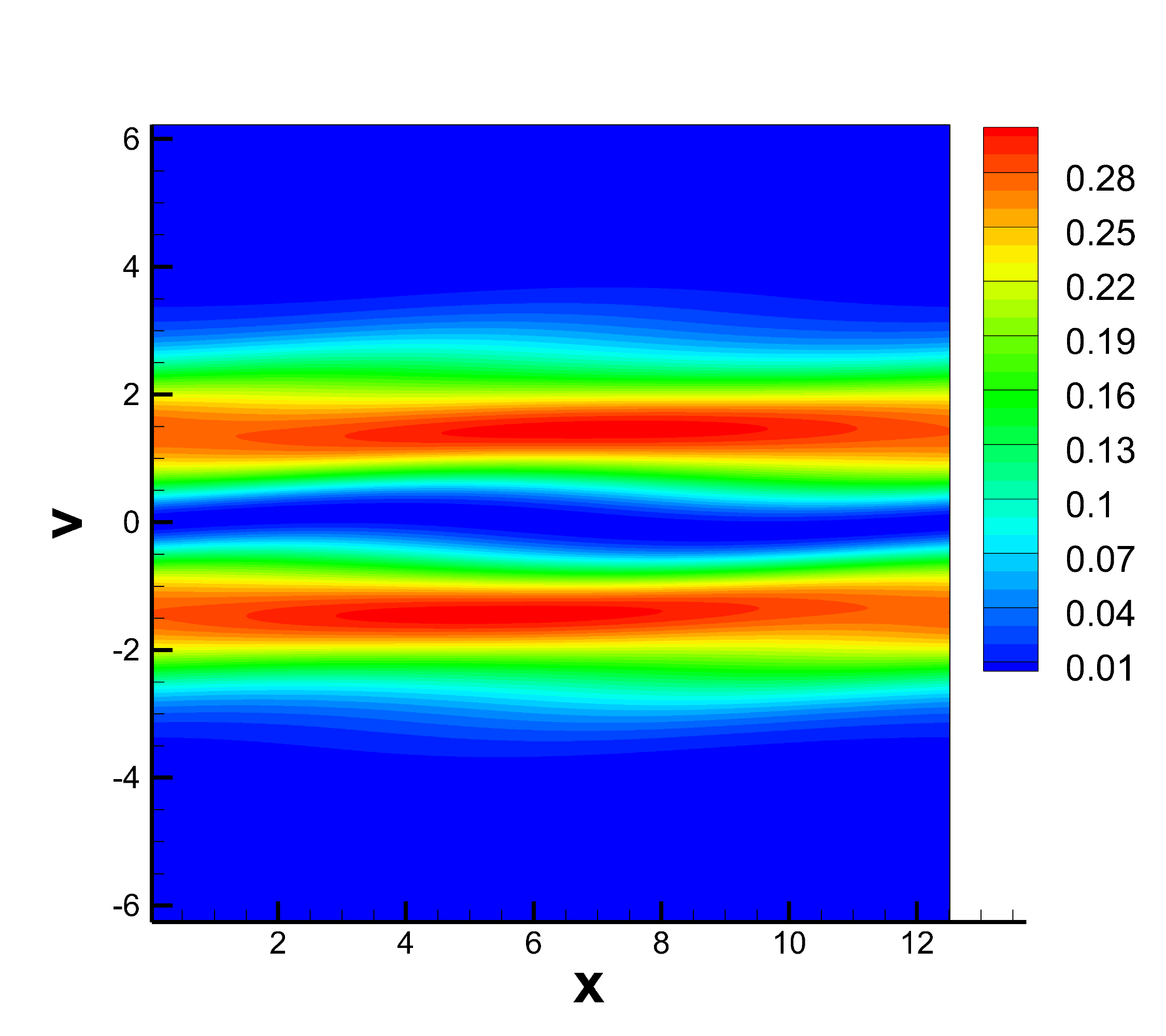}}
\subfigure[]{\includegraphics[width=.42\textwidth]{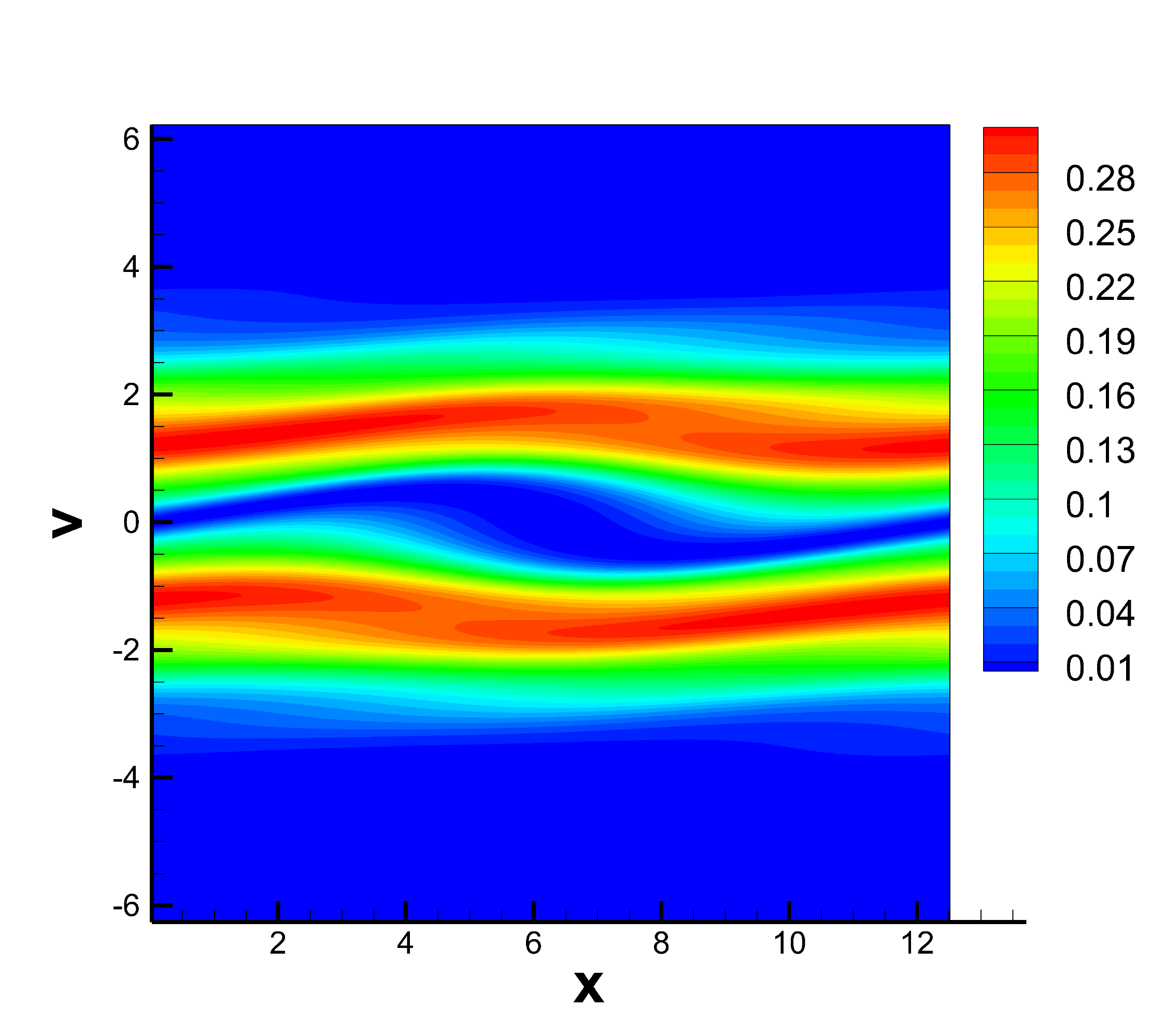}}\\
\subfigure[]{\includegraphics[width=.42\textwidth]{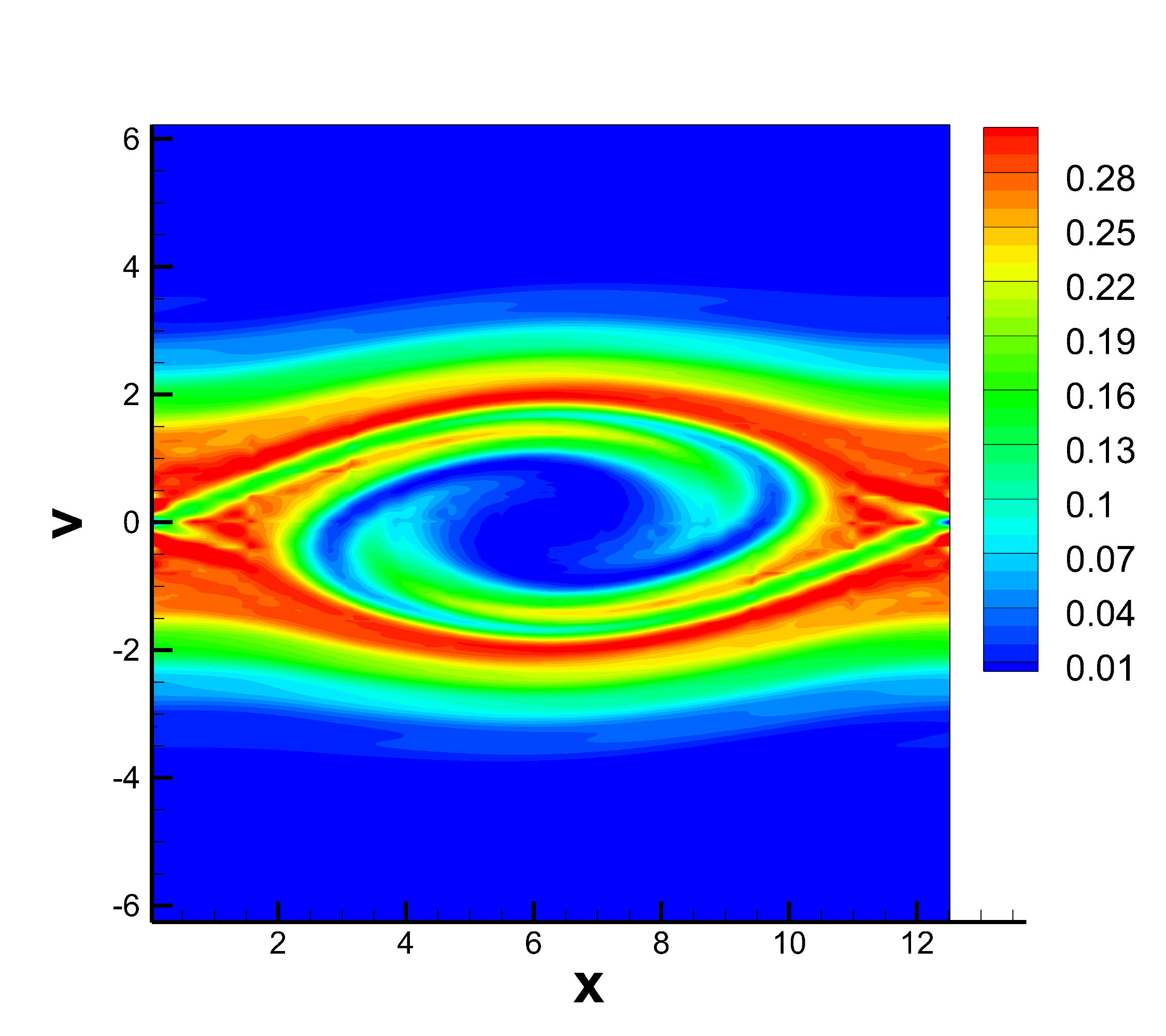}}
\subfigure[]{\includegraphics[width=.42\textwidth]{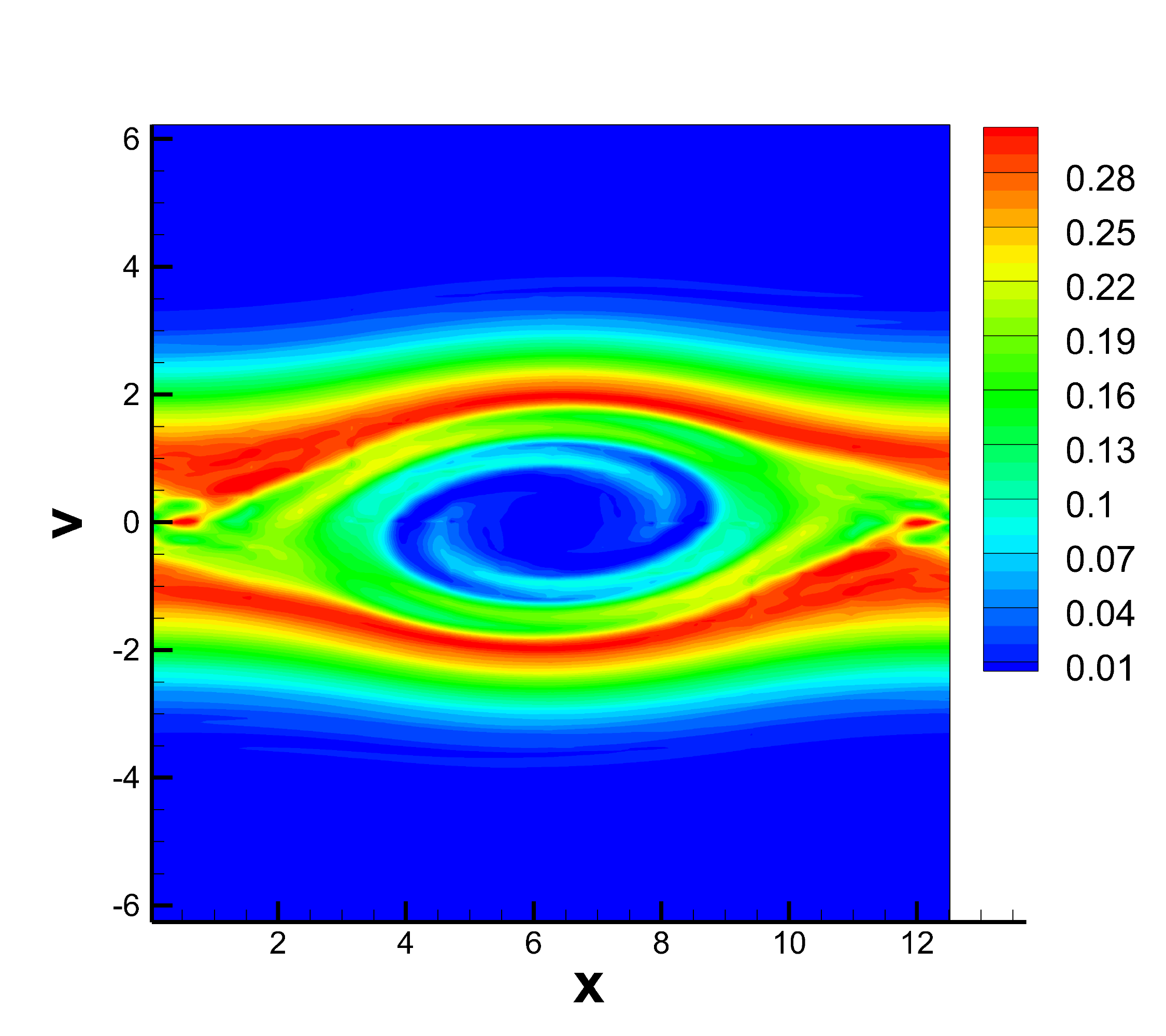}}\\
\end{center}
  \caption{Two-stream instability. Phase space contour plots at   $t=5$ (a),  $t=10$ (b), $t=20$ (c), and $t=40$ (d). $k=3$, $N=8$.}
 \label{fig:con_two}
\end{figure}

\subsection{The relaxation model}

The last example we present in this paper is the linear Vlasov-Boltzmann transport equation
\begin{equation}
\label{eq:vb}
f_t + \bv\cdot\nabla_\bx f + \bE(t,\bx)\cdot\nabla_\bv f=L(f),
\end{equation}
where $L(f)$ denotes the linear relaxation operator
$$L(f)=\frac{\mu_\infty(\bv)\rho(t,\bx)-f(t,\bx,\bv)}{\tau},$$
and $\mu_\infty(\bv)$ is an absolute Maxwellian distribution defined as 
$$\mu_\infty(\bv) = \frac{\exp(-\frac{|\bv|^2}{2\theta})}{(2\pi\theta)^{d/2}},$$
and $$\rho(t,\bx)=\int_{\bv}f(t,\bx,\bv)\,d\bv$$
denotes the macroscopic density. The external electric field $\bE(t,\bx)$  is given by a known  electrostatic potential 
$$ \bE(\bx)= -\nabla_\bx \Phi(\bx)\quad \text{with}\quad \Phi(x) = \frac{|\bx|^2}{2}.$$
In this case, the unique stationary state solution $\mathcal{M}(\bx,\bv)$ of equation \eqref{eq:vb} is the global Maxwellian distribution $\mu_\infty(\bv)$ multiplied by the stationary macroscopic density  given by the spatial Maxwellian
$$\rho_\infty(\bx)=\frac{\exp(-\Phi(\bx)/\theta)}{\int_{\bx}\exp(-\Phi(\bx)/\theta)\,d\bx},$$
i.e.,
$$\mathcal{M}(\bx,\bv)=\rho_\infty(\bx)\mu_\infty(\bv) = \frac{\exp\left(-\left(\frac{|\bv|^2}{2}+\Phi(\bx)\right)/\theta\right)}{(2\pi\theta)^{d/2}\int_{\bx}\exp(-\Phi(\bx)/\theta)\,d\bx}.$$

A DG scheme with traditional piecewise polynomial space for \eqref{eq:vb} has been developed and analyzed in \cite{cheng_gamba_proft_2010}. Here, we use the sparse DG scheme to solve the problem  \eqref{eq:vb}  in a cut-off domain $\Omega=[-L,L]^d\times[-V_{c},V_{c}]^{d}$ with $d=1,\,2$ (two dimensional and four dimensional calculations, resp.).
The following initial conditions are used:
\begin{align*}
f(0,\bx,\bv) &= \frac{1}{s_1}\sin(x^2/2)^2\exp(-(x^2+v^2)/2),\quad\text{when}\quad d=1,\\
f(0,\bx,\bv) &= \frac{1}{s_2}\sin(x_1^2/2)^2\cos(x_2^2/2)^2\exp(-(x_1^2+x_2^2+v_1^2+v_2^2)/2),\quad\text{when}\quad d=2,
\end{align*}
where $s_1,\,s_2$ are normalization constants such that $\int_{\Omega}f(0,\bx,\bv)d\bx d\bv=1$.
The zero boundary conditions are imposed in both $\bx-$ and $\bv-$spaces.
We take $V_{c}=5$, and $\tau=\theta=1$. 

For the case of $d=1$, we investigate the decay rate of the initial state to equilibrium by tracking time evolution of the following two entropy functionals
$$\mathcal{H}_{log}(t)=\int_{\Omega}H\log(H)\mathcal{M}\,dxdv,\quad\mathcal{H}_2(t)=\int_{\Omega}H^2\mathcal{M}\,dxdv,$$
where $H(t,x,v)=f_h/\mathcal{M}$ is the global relative entropy function. $\mathcal{H}_{log}$ and $\mathcal{H}_2$ should relax to $0$ and $1$ indicating convergence to equilibrium. In the simulation, we use the sparse DG scheme with $\hat{\bV}^3_7$ and plot the decay rate with respect to both entropy functionals in Figure \ref{fig:relax_vb}. Both plots verify the convergence of numerical solution to $\mathcal{M}.$ 
For the case of $d=2$, the sparse DG method with approximation space $\hat{\bV}^3_8$ is used.  In Figures  \ref{fig:con_vb4_x0v0}-\ref{fig:con_vb4_x0x1}, we report the  evolution of the two-dimensional cuts in $x_1-v_1$ plane at $x_2=0,\,v_2=0$ and in $x_1-x_2$ plane  at $v_1=0,\,v_2=0$. Similar to the case of $d=1$,  
 $f_h$ in both two-dimensional cuts is observed to relax towards the equilibrium distribution $\mathcal{M}$. 

\begin{figure}[htp]
	\begin{center}		\subfigure[]{\includegraphics[width=.42\textwidth]{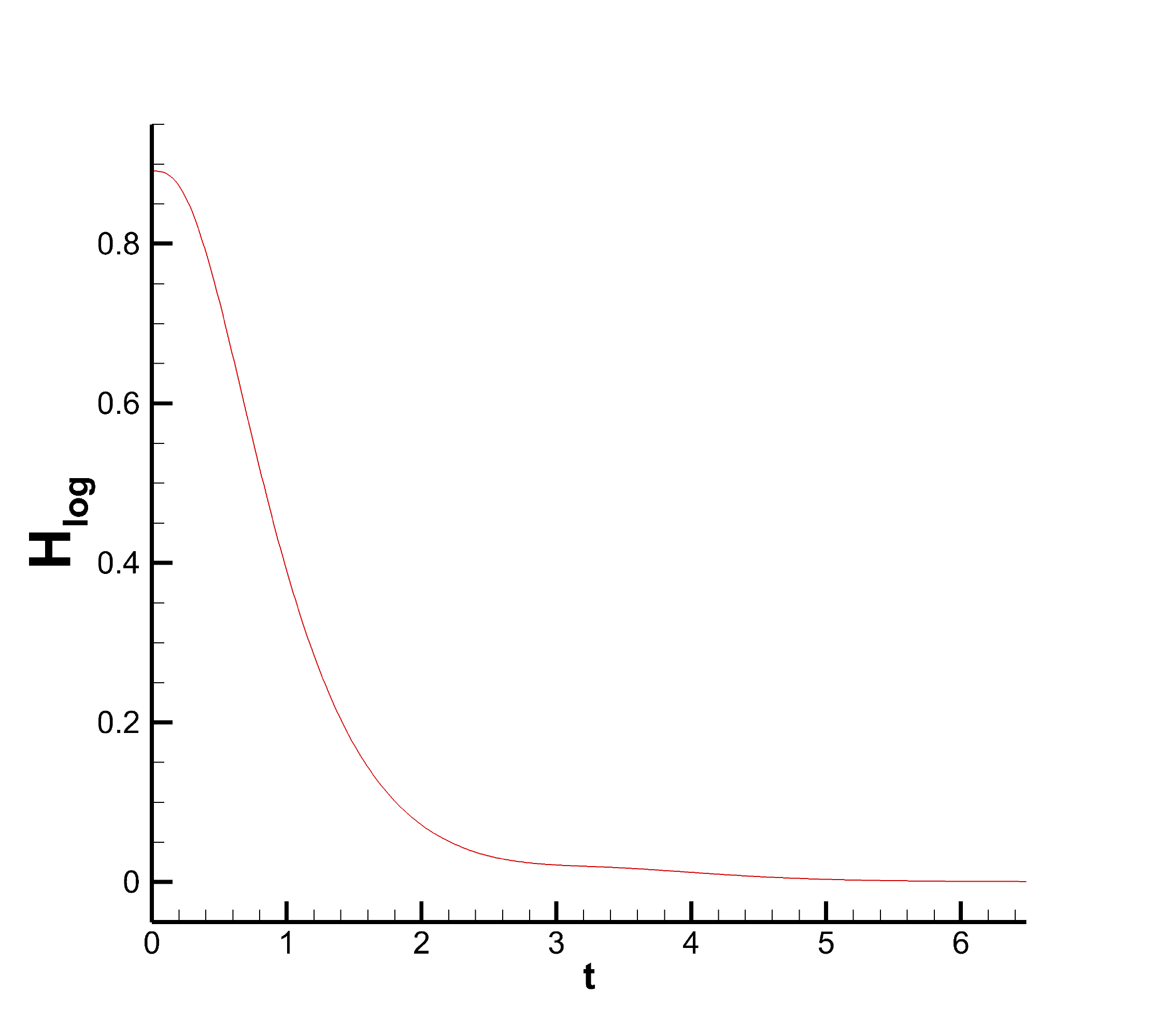}}
	\subfigure[]{\includegraphics[width=.42\textwidth]{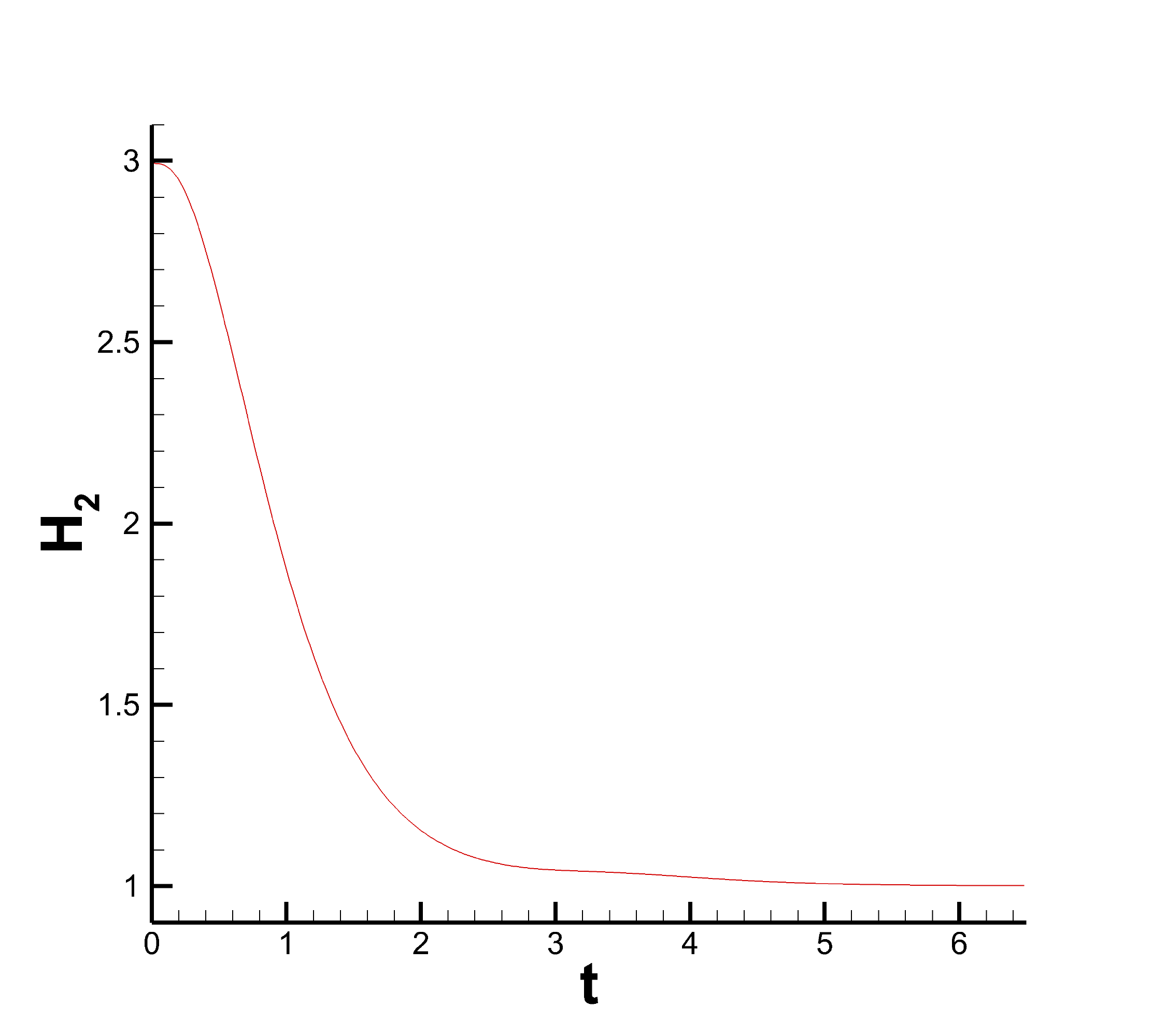}}
	\end{center}
	\caption{Linear Vlasov-Boltzmann equation. Decay rate for entropy functional $\mathcal{H}_{log}$ (a) and $\mathcal{H}_{2}$ (b).  $k=3$, $N=7$, $d=1$.}
	\label{fig:relax_vb}
\end{figure}


\begin{figure}[htp]
	\begin{center}
		\subfigure[]{\includegraphics[width=.42\textwidth]{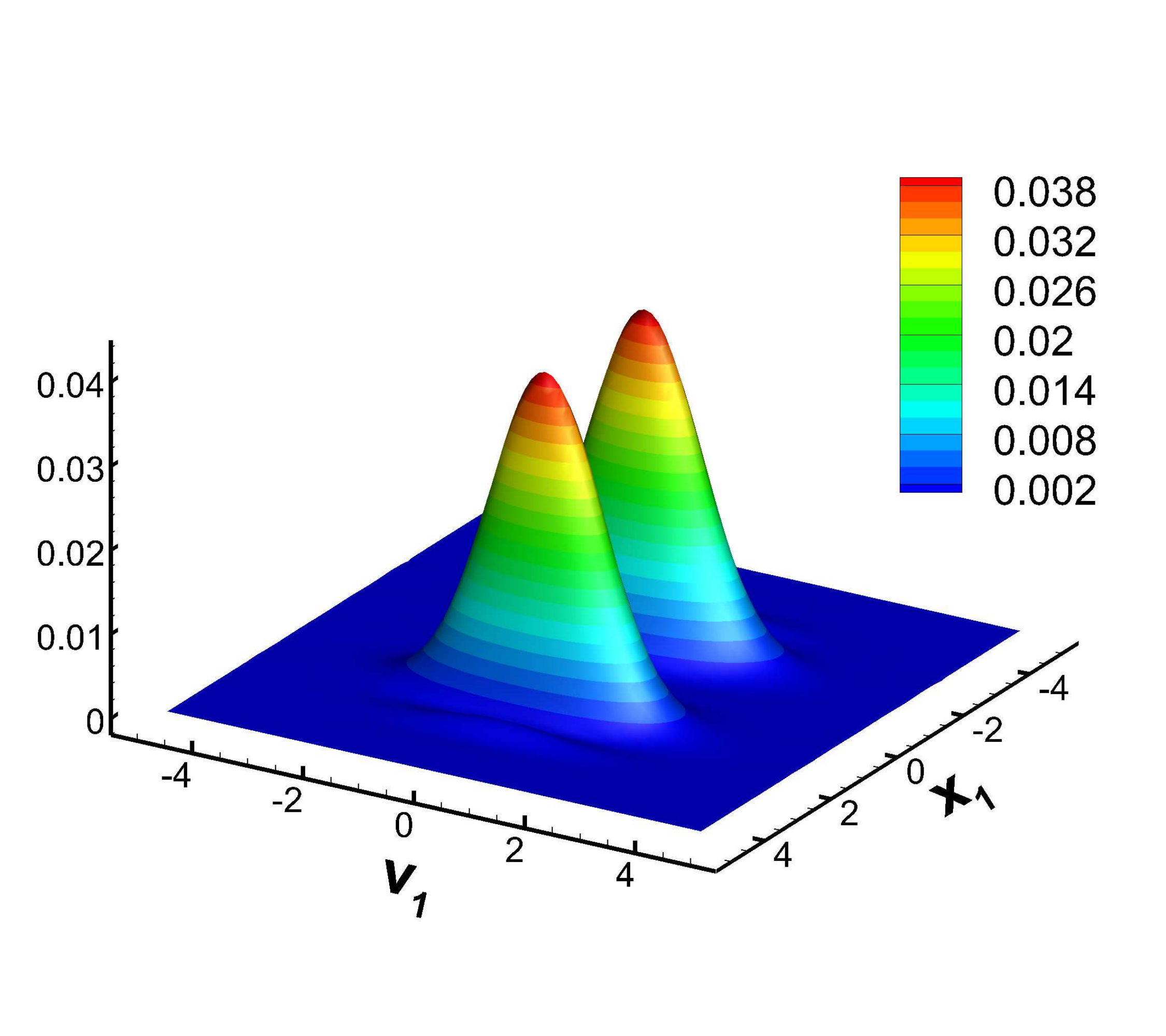}}
		\subfigure[]{\includegraphics[width=.42\textwidth]{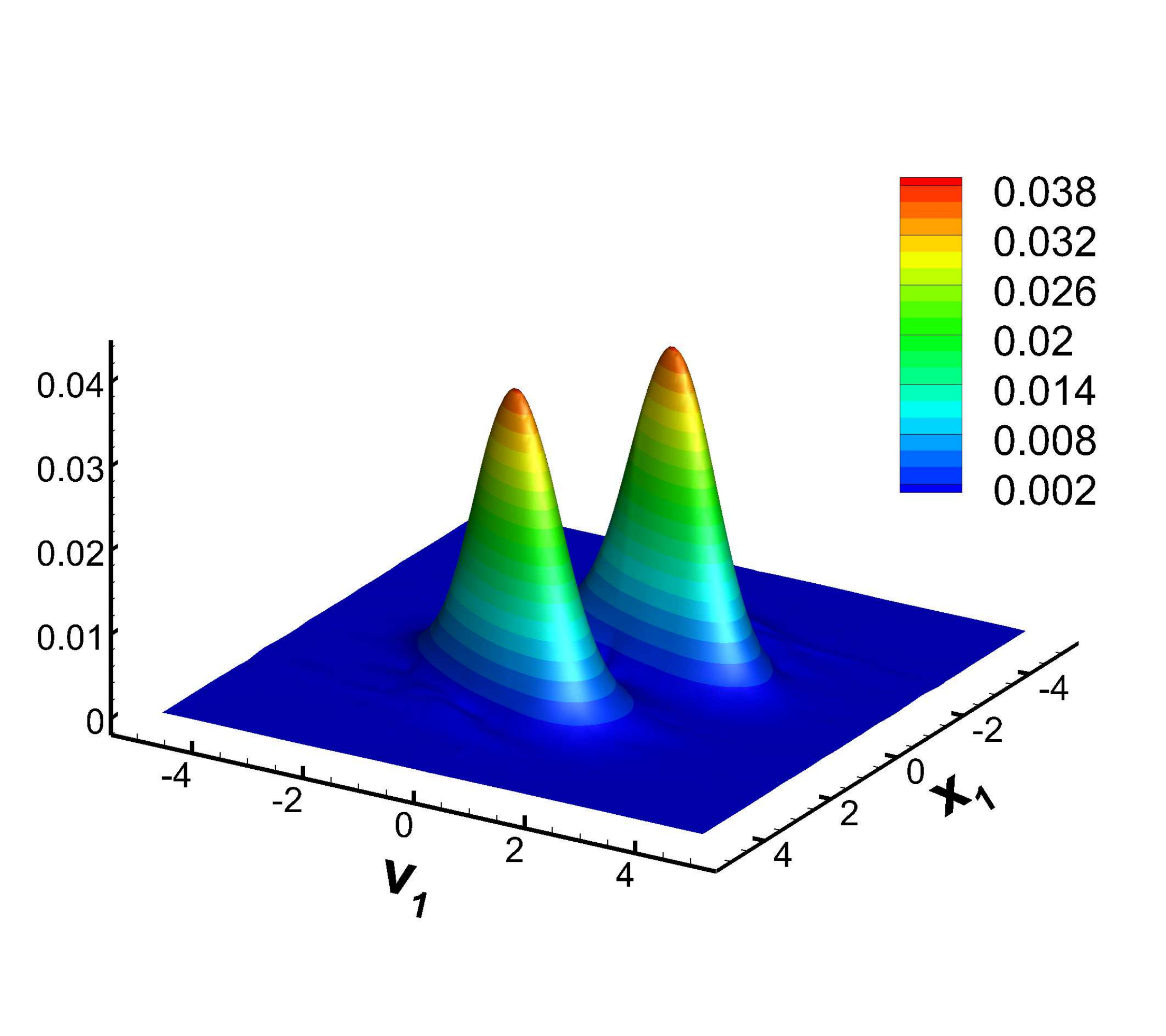}}\\
		\subfigure[]{\includegraphics[width=.42\textwidth]{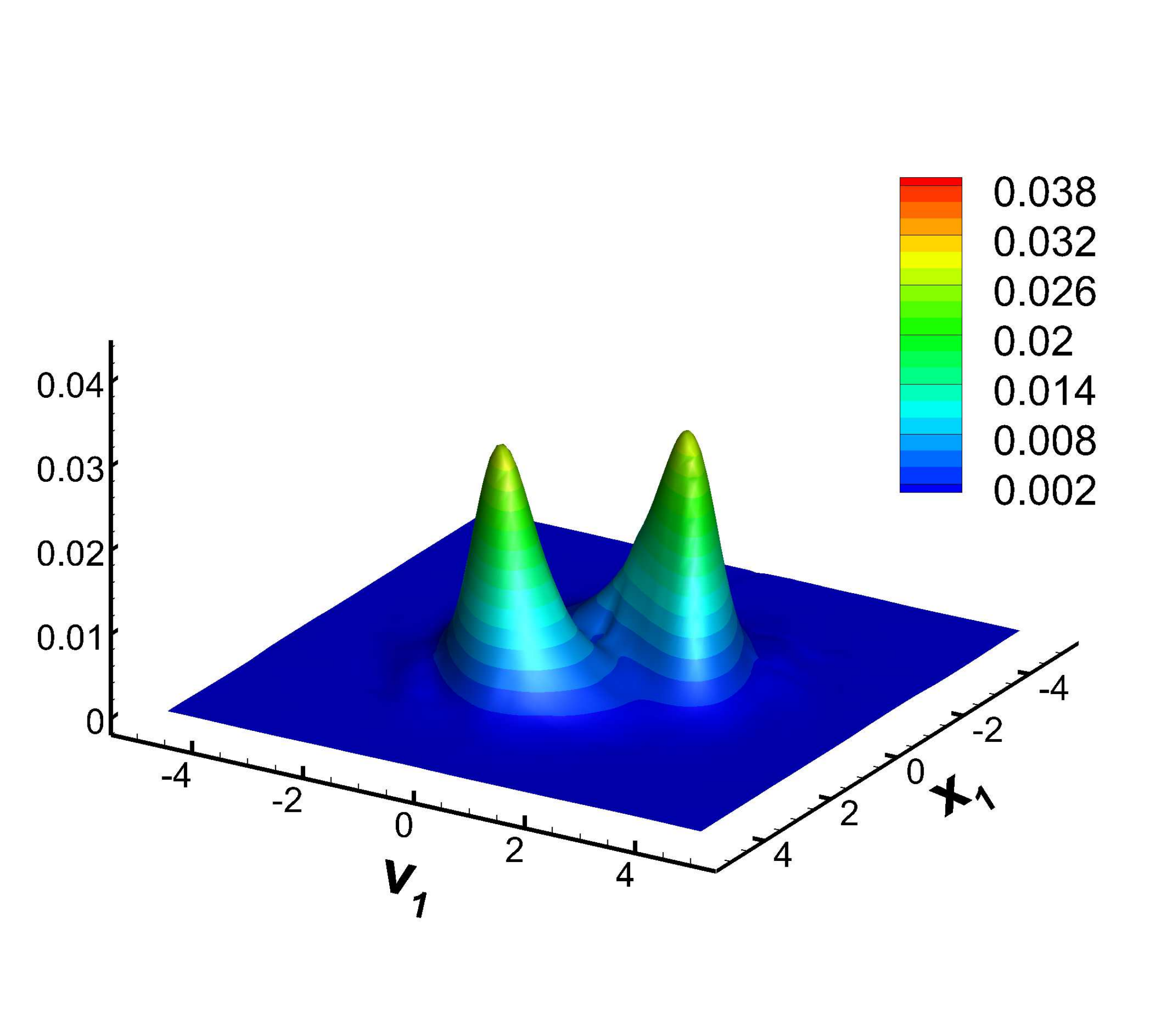}}
		\subfigure[]{\includegraphics[width=.42\textwidth]{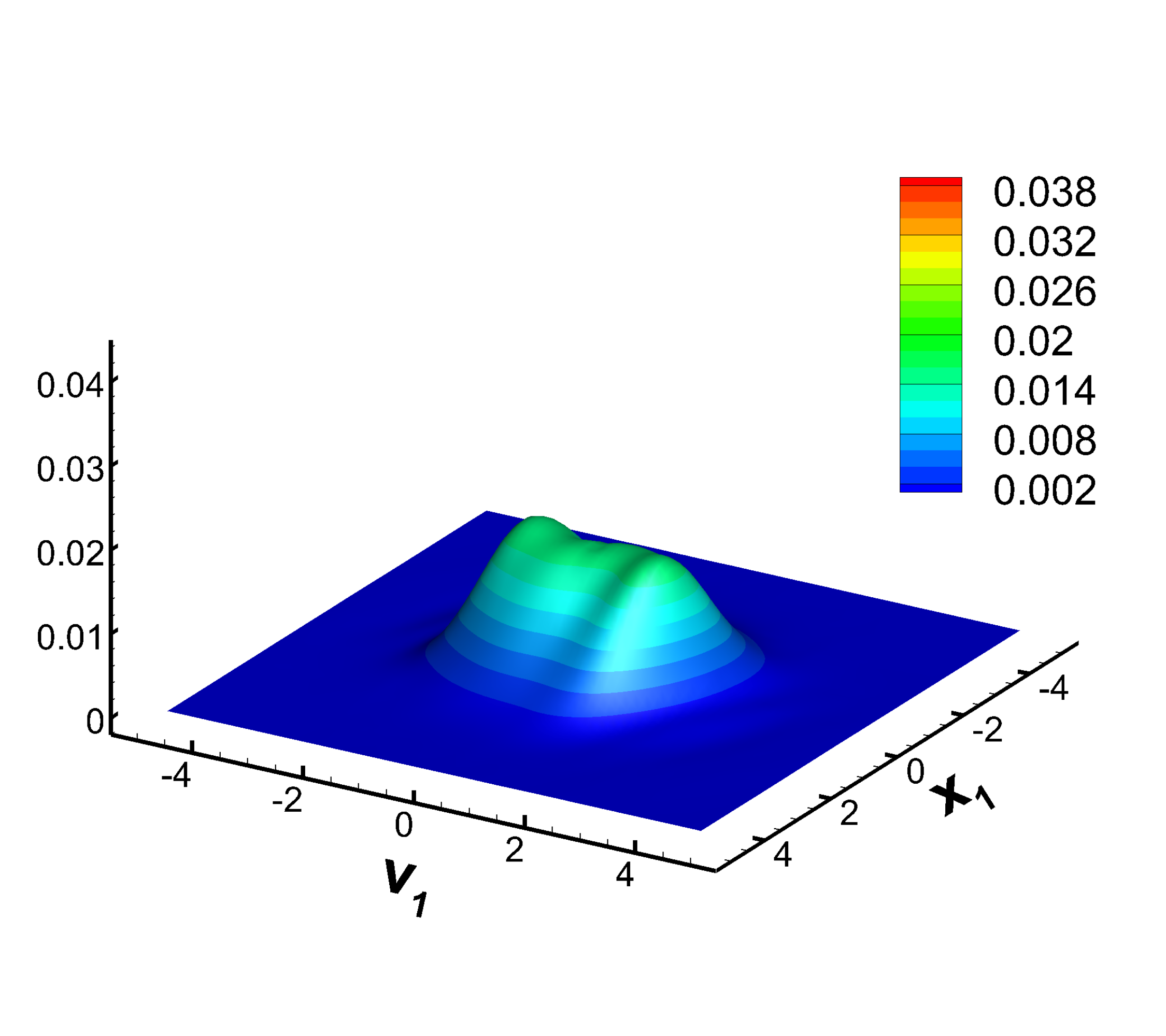}}\\
		\subfigure[]{\includegraphics[width=.42\textwidth]{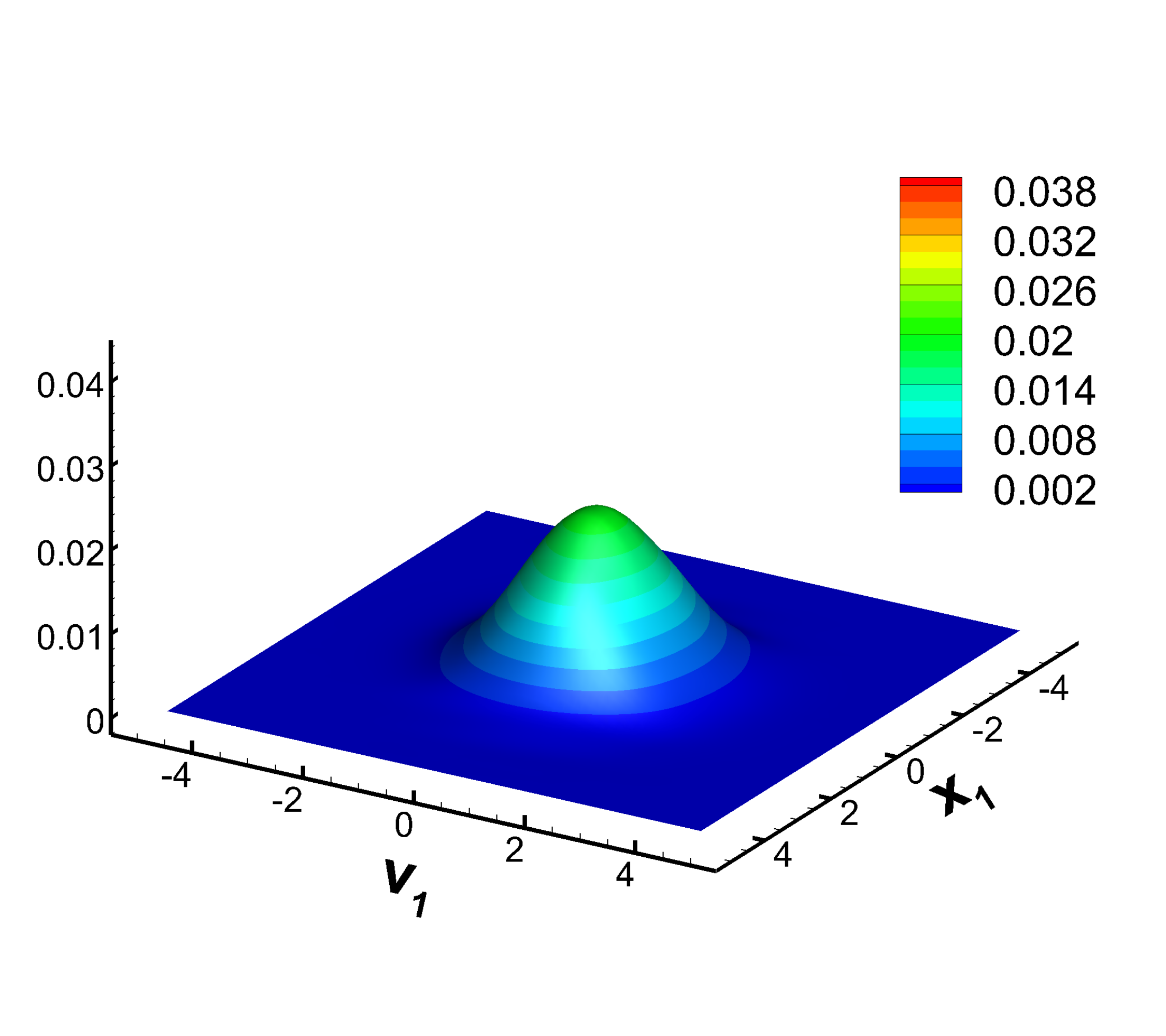}}	
		\subfigure[]{\includegraphics[width=.42\textwidth]{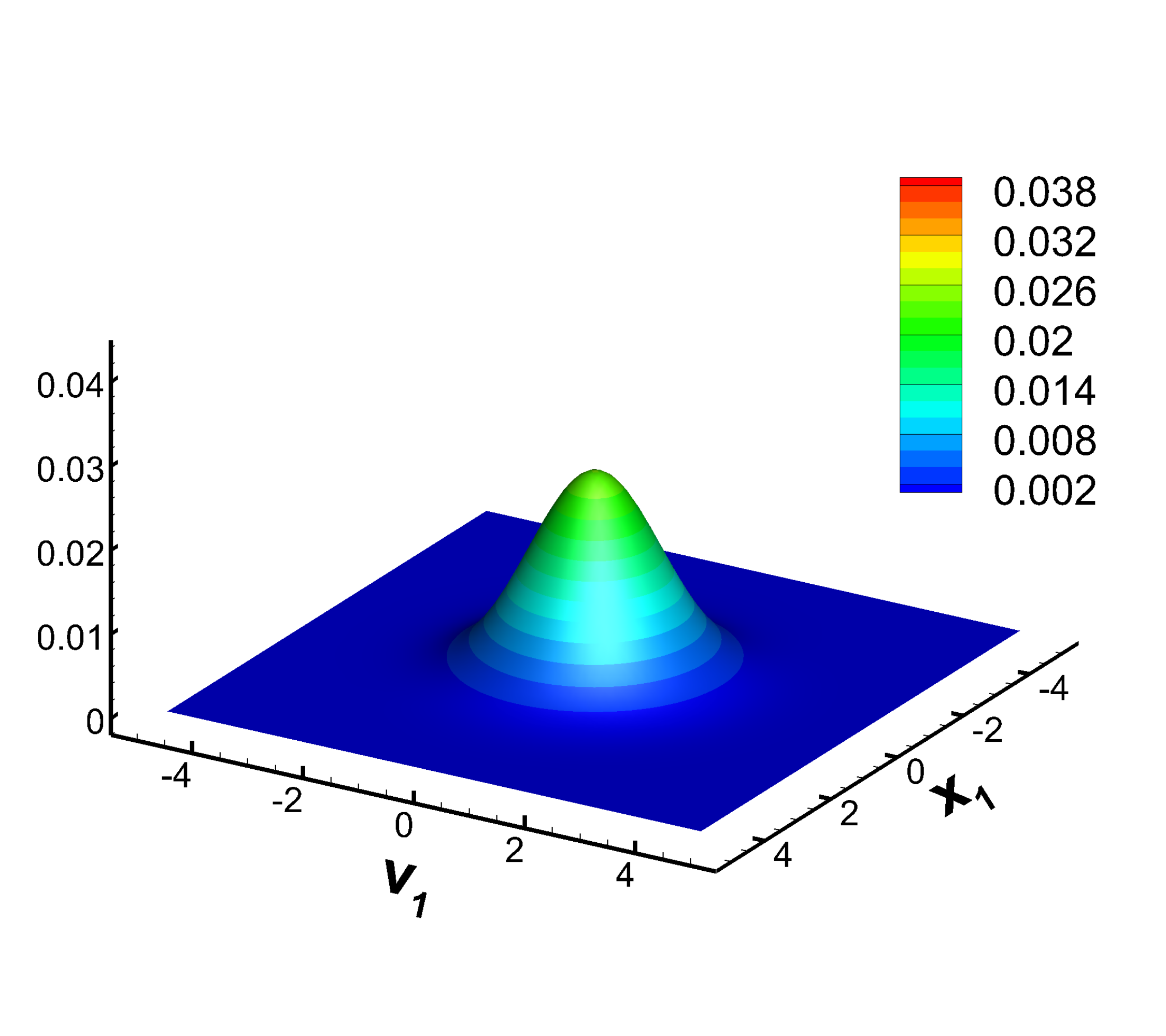}}
	\end{center}
	\caption{Linear Vlasov-Boltzmann equation. The two-dimensional cuts in $x_1-v_1$ plane of the evolution of $f_h$ towards equilibrium at $x_2=0$ and $v_2=0$.  $t=0$ (a), $t=0.5$ (b),  $t=1$ (c), $t=2$ (d), $t=3$ (e), and $t=6$ (f). $k=3$, $N=8$, $d=2$.}
	\label{fig:con_vb4_x0v0}
\end{figure}

\begin{figure}[htp]
	\begin{center}
		\subfigure[]{\includegraphics[width=.42\textwidth]{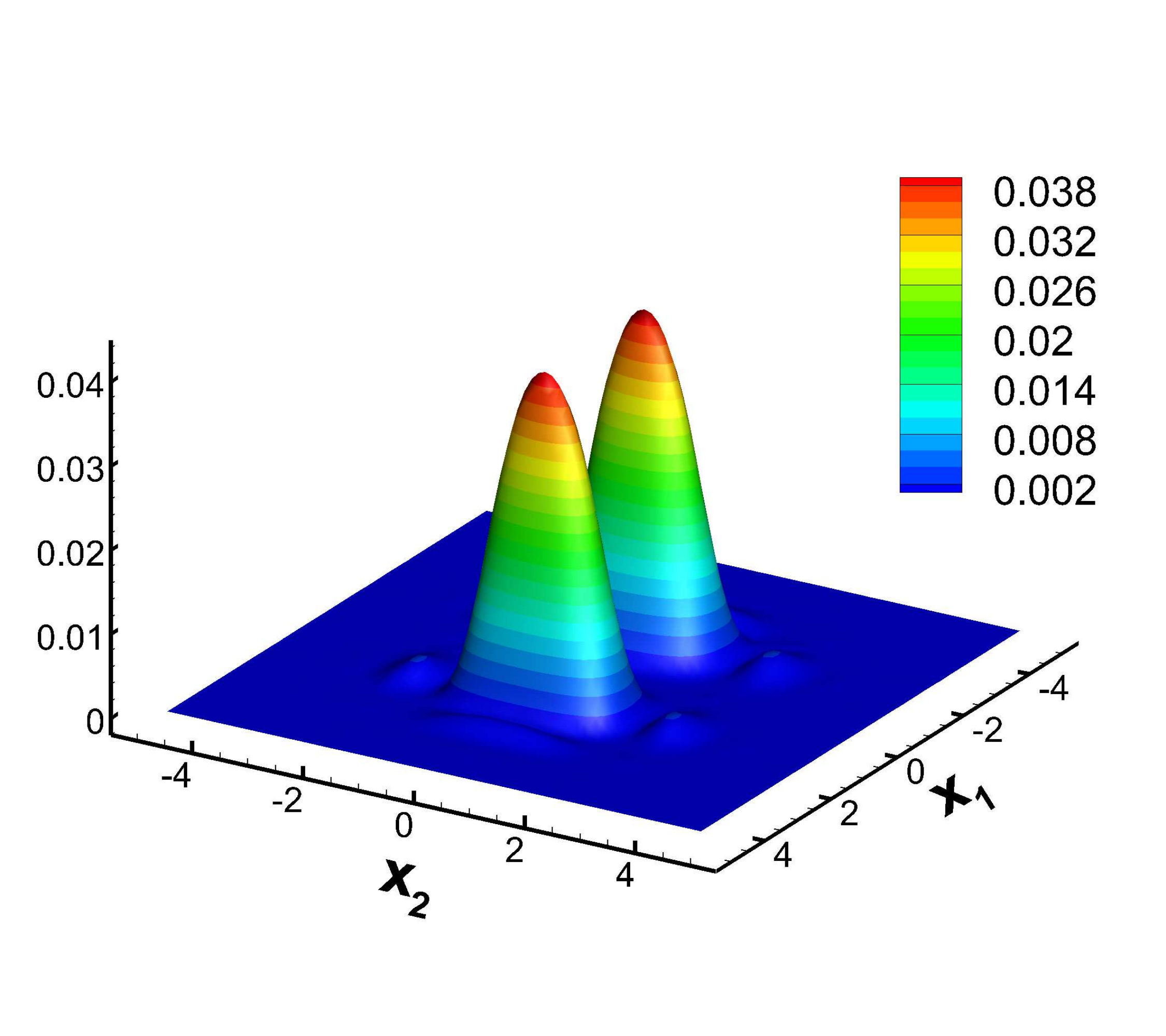}}
		\subfigure[]{\includegraphics[width=.42\textwidth]{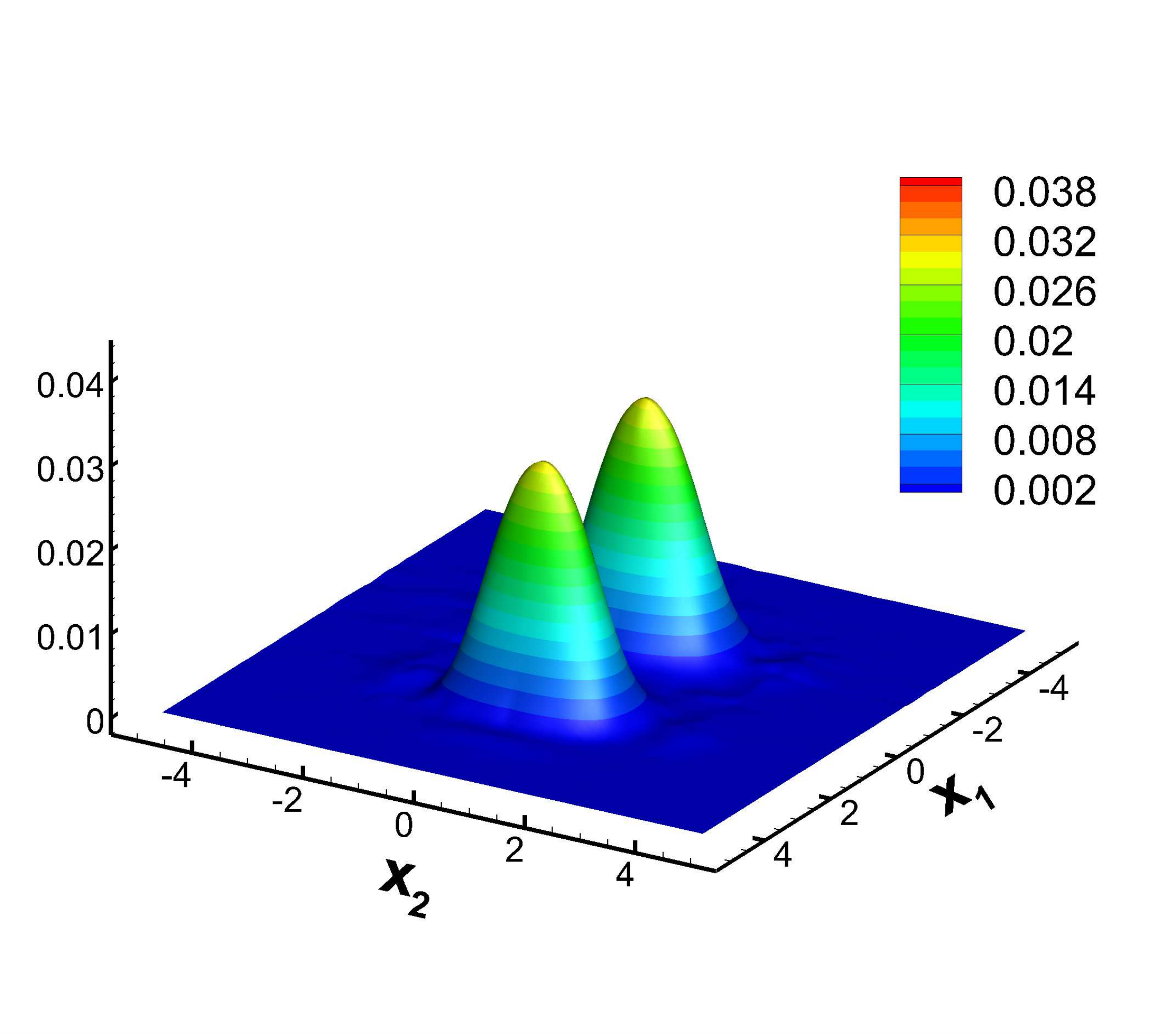}}\\
		\subfigure[]{\includegraphics[width=.42\textwidth]{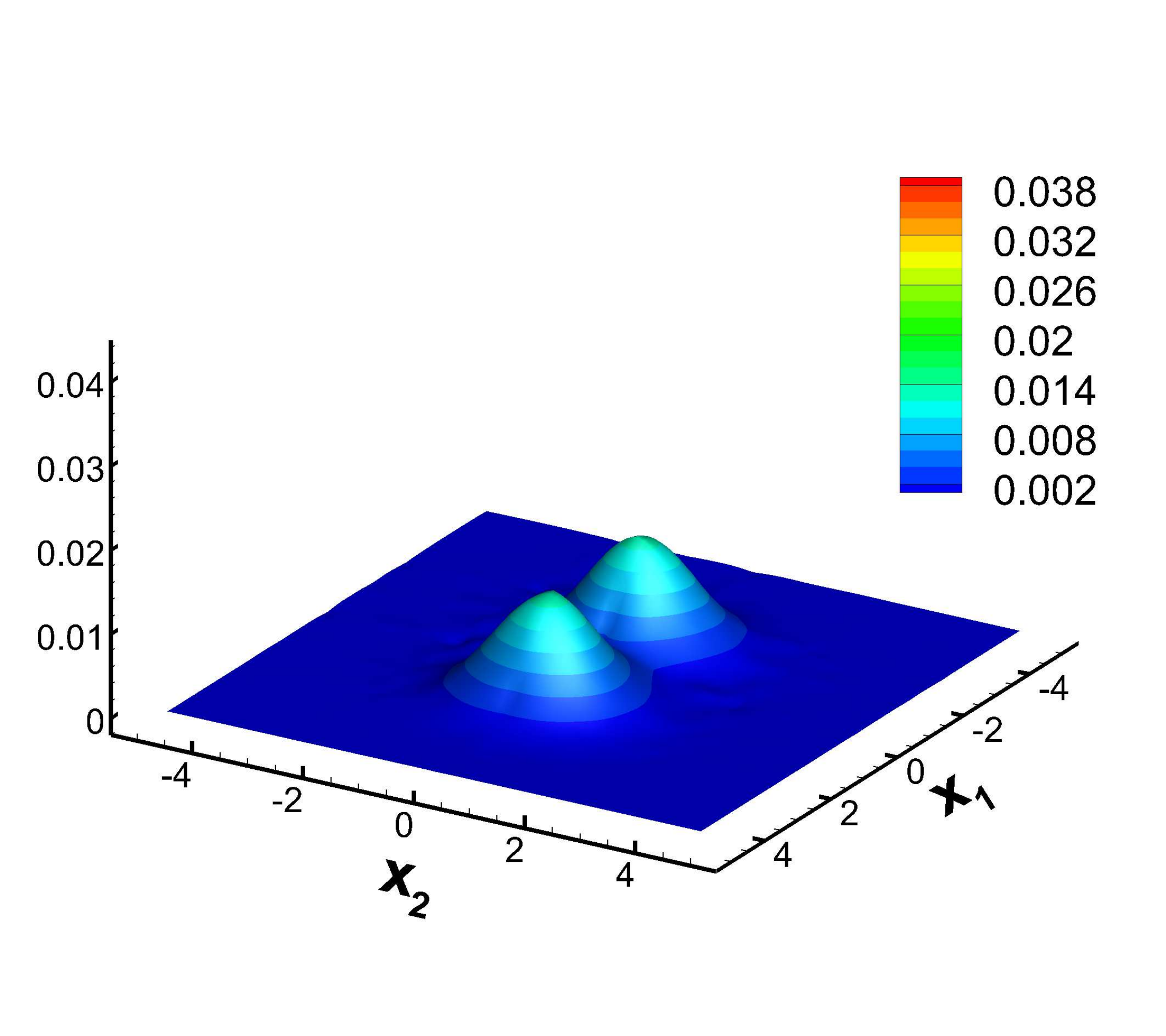}}
		\subfigure[]{\includegraphics[width=.42\textwidth]{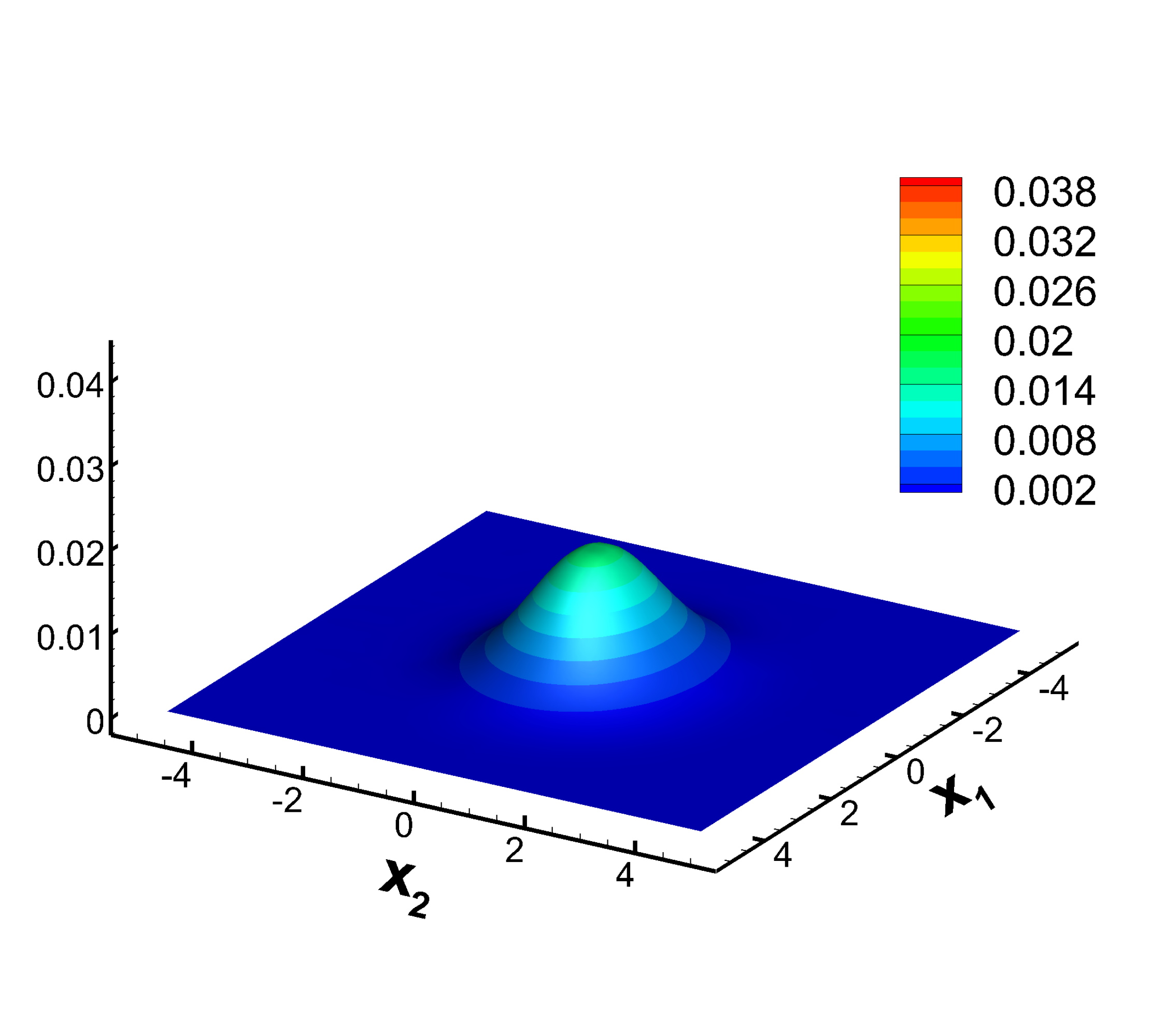}}\\
		\subfigure[]{\includegraphics[width=.42\textwidth]{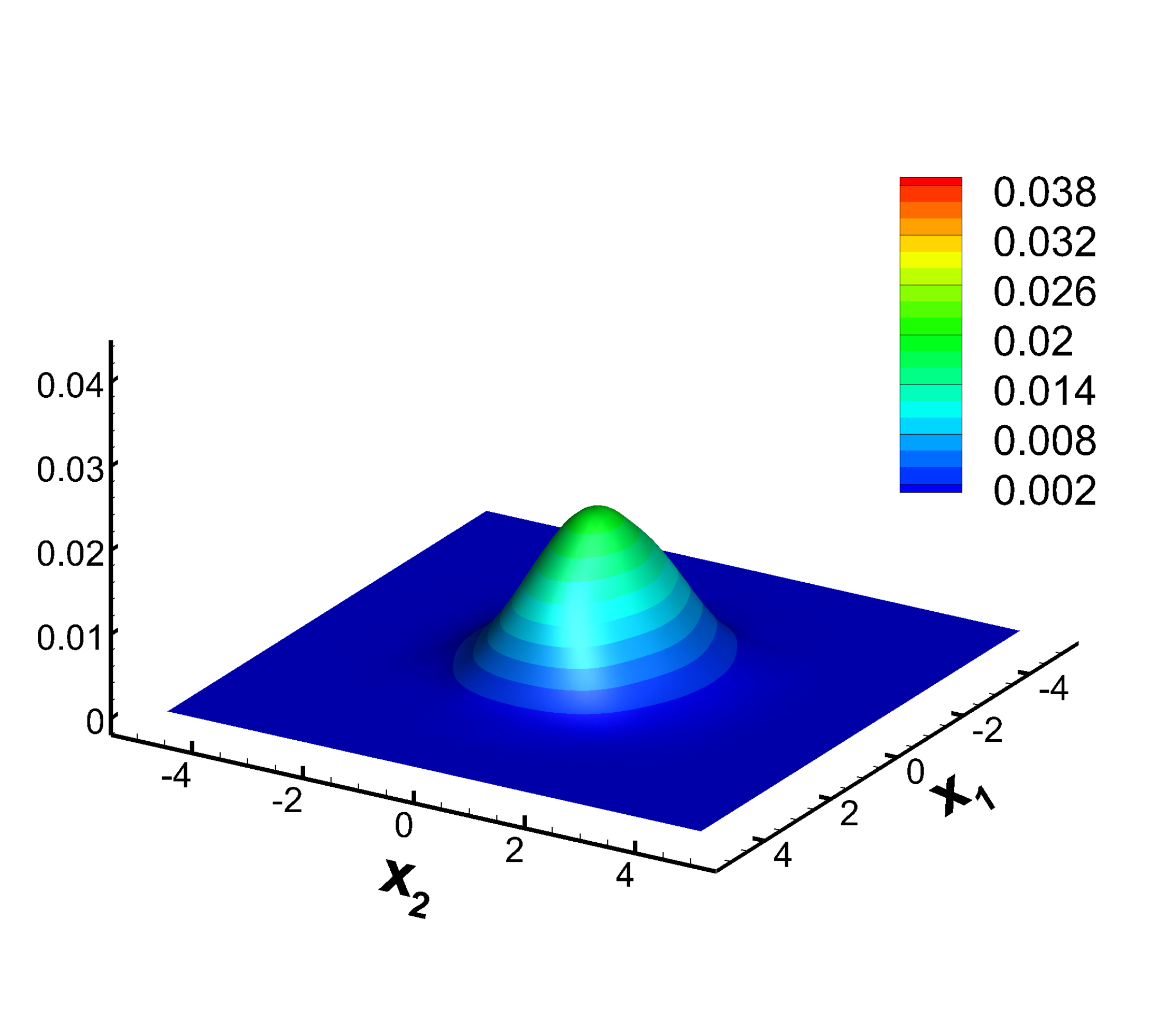}}	
		\subfigure[]{\includegraphics[width=.42\textwidth]{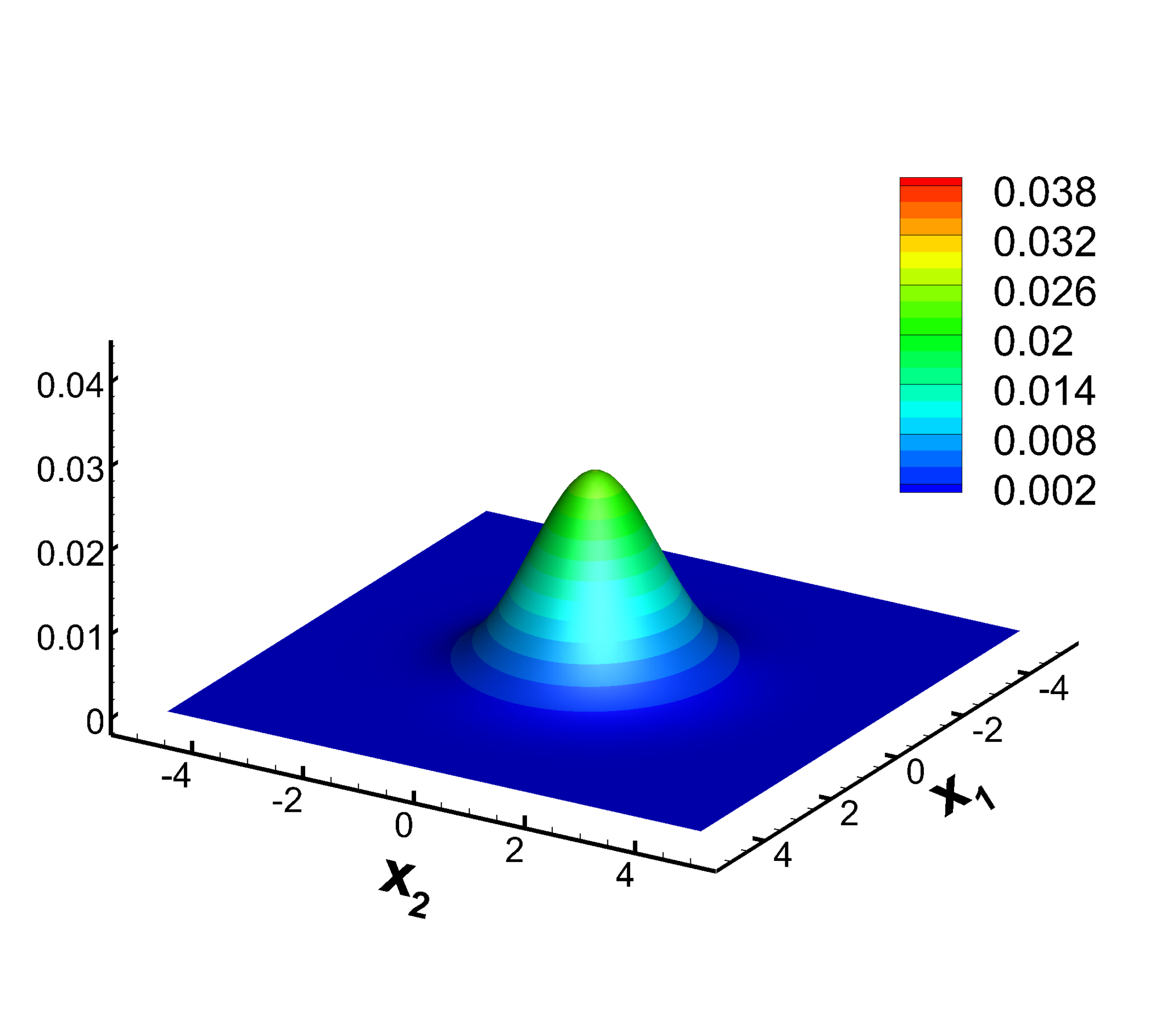}}
	\end{center}
	\caption{Linear Vlasov-Boltzmann equation. The two-dimensional cuts in $x_1-x_2$ plane at $v_1=0$ and $v_2=0$ of the evolution of $f_h$ towards equilibrium.  $t=0$ (a), $t=0.5$ (b),  $t=1$ (c), $t=2$ (d), $t=3$ (e), and $t=6$ (f). $k=3$, $N=8$, $d=2$.}
	\label{fig:con_vb4_x0x1}
\end{figure}

\section{Conclusions and future work}
\label{sec:Con}

In this paper, we developed a sparse grid DG scheme for variable-coefficient transport equations. The method uses a sparse finite element space based on a hierarchical construction of basis functions and are shown to reduce the computational degrees of freedom significantly in high dimensions. Weak formulations of traditional DG methods are incorporated ensuring many nice properties such as conservation and stability. For constant coefficient equations, we established  semi-discrete  $L^2$ stability and error estimate of order $O( (\log h)^d h^{k+1/2})$  for smooth enough solutions where $h$ is the size of the finest mesh in each dimension. The method is   applied to kinetic simulations of Vlasov and Boltzmann equations. Good performance in accuracy and conservation is observed. Future work includes further investigation of error estimates in  various norms, and development of the scheme for solutions with less regularities.



\appendix

\section{Proof of Lemma \ref{thm:appx}}

In this appendix, we prove the approximation results of  Lemma \ref{thm:appx} in three steps. Many discussions are closely related to  \cite{schwab2008sparse}, where the $C^0$ sparse finite element space was considered. The main difference lies in the splitting of the errors into two parts to fine tune the estimates for discontinuous piecewise polynomials in $ \hat{\bV}_N^k$.

\subsection{An alternative representation of the $L^2$ projection}
To facilitate the discussion,  we seek an alternative representation of the $L^2$ projection $\mathbf{P}$ onto $ \hat{\bV}_N^k$ following a similar   construction in \cite{schwab2008sparse}. We   denote $P_n^k$ as the standard $L^2$ projection operator from  $L^2(0,1)$ to $V_n^k$, and the induced
 increment projector
\[Q_n^k: =
 \left\{
  \begin{array}{ll}
P_n^k - P_{n-1}^k,  & \textrm{if} \, n \ge1,\\
P_0^k, & \textrm{if} \, n = 0, 
  \end{array} \right.\]
and further denote
\begin{equation}
\label{eqn:multiproj}
\hat{\bP}_N^k:=\sum_{\substack{ |\bl|_1 \leq N\\\bl \in \mathbb{N}_0^d}} {Q}^k_{l_1, x_1} \otimes \cdots \otimes {Q}^k_{l_d, x_d},
\end{equation}
where the second subindex of ${Q}^k_{l_i, x_i}$ indicates that the increment operator is defined in $x_i$-direction.
We can easily show that $\bP=\hat{\bP}_N^k.$ In fact, for any $v$, it's clear that $\hat{\bP}_N^k v\in  \hat{\bV}_N^k$. Therefore, we only need
\begin{equation}
\label{eq:portho}
\int_\Omega (\hat{\bP}_N^k v -v) w \, d\bx=0, \qquad \forall \, w \in  \hat{\bV}_N^k.
\end{equation}
It suffices to show \eqref{eq:portho} for $v \in C^\infty(\Omega)$ which is a dense subset of $L^2(\Omega)$. In fact, we have 
$$
v=\sum_{\bl \in \mathbb{N}_0^d}  {Q}^k_{l_1, x_1} \otimes \cdots \otimes {Q}^k_{l_d, x_d} v.
$$
Therefore,
\begin{equation}
\label{eq:portho1}
\int_\Omega (\hat{\bP}_N^k v -v) w \, d\bx= \int_\Omega   (\sum_{\substack{ |\bl|_1>N\\\bl \in \mathbb{N}_0^d}}  {Q}^k_{l_1, x_1} \otimes \cdots \otimes {Q}^k_{l_d, x_d} v  ) \,w \, d\bx.
\end{equation}
Since the one dimensional projectors satisfy
$$
\int_{[0,1]} Q_n^k \phi \, \varphi \,dx=0, \quad   \varphi \in V_{n-1}^k,
$$
for any $ n\ge1, \phi \in L^2(0,1).$ \eqref{eq:portho1} is immediate from the definition of the tensor-product operators and we are done.

%



\subsection{Properties of the one-dimensional and tensor product projections}
Here, we review some classical approximation results about the one-dimensional projection $P_n^k$  and  their tensor-product constructions.

\begin{prop}[Convergence property of the 1D projections \cite{Ciarlet_1975_FEM_Elliptic}]
\label{prop2}

For a function $v \in  H^{p+1}(0,1)$,    we have the convergence property of the $L^2$ projection $P_n^k$ as follows: for any integer $t$ with $1 \leq q \leq \min\{p, k\}$, $s=0,\,1$,
\begin{equation}
\label{eq:proj1}
 | P_n^k v- v|_{H^s(I^j_n)}\le c_{k,s, q} 2^{-n(q+1-s)} |v|_{H^{q+1}(I^j_n)},\quad j=0,\cdots,2^{n+1}-1,
 \end{equation}
where $c_{k,s, q}$ is a constant that depends on $k, s, q$, but not on $n$, and the case of $s=0$ refers to the $L^2$ norm. 
\end{prop}

From this property, using basic algebra, we can deduce that for $n \geq 1$,
\begin{equation}
 | Q_n^k v  |_{H^s(I^j_n)}\le \tilde{c}_{k,s, q} 2^{-n(q+1-s)} |v|_{H^{q+1}(I^j_n)},
\end{equation}
with
\begin{equation}
\label{eq:proj2}
\tilde{c}_{k,s, q}=c_{k,s, q} \,(1+2^{-(q+1-s)}).
 \end{equation}

When $n=0$, first note that the $L^2$ projector $Q_0^k=P_0^k$ preserves the $L^2$ norm, i.e.,
$$\| Q_0^k v\|_{L^2(0,1)}= \|v\|_{L^2(0,1)}.$$
Moreover, by the inverse inequality, see, e.g., \cite{schwab1998p}, we have
$$| Q_0^k v|_{H^1(0,1)}\le c^k_{inv} \|Q_0^k v\|_{L^2(0,1)}=c^k_{inv} \| v\|_{L^2(0,1)},$$
where $c^k_{inv} = \sqrt{12}k^2.$
Therefore, we have obtained the following estimation
 \begin{equation}
\label{eq:proj3}
 | Q_0^k v |_{H^s(0,1)}\le \hat{c}_{k,s} \|v\|_{L^2(0,1)},\quad s=0,\,1,
 \end{equation}
where $\hat{c}_{k,0}=1$ and $\hat{c}_{k,1}=c_{inv}^k$.

For multi-dimensions, if we consider the  $L^2$ projection $\bP^k_{N}={P}^k_{N, x_1} \otimes \cdots \otimes {P}^k_{N, x_d}$ onto the standard piecewise polynomial space $\bV_N^k$, we have the following approximation results.

\begin{prop}[Convergence property of the multi-dimensional tensor-product projectors \cite{Ciarlet_1975_FEM_Elliptic}] 
\label{prop:tensor}

For a function $v \in  H^{p+1}(\Omega)$,    we have the convergence property of the $L^2$ projection $\bP^k_{N}$ onto $\bV_N^k$ as follows: for any integer $t$ with $1 \leq q \leq \min\{p, k\}$, $s=0,\,1$,
\begin{equation}
\label{eq:tensorproj}
 | \bP_N^k v- v|_{H^s(\Omega_N)}\le \bar{\bar{c}}_{k,s, q} 2^{-N(q+1-s)} |v|_{H^{q+1}(\Omega)}, 
 \end{equation}
where $\bar{\bar{c}}_{k,s, q}$ is a constant that depends on $k, s, q$, but not on $N$. 
\end{prop}

\subsection{Proof of Lemma \ref{thm:appx}} In this subsection, we will prove Lemma \ref{thm:appx}.   For any function $v \in L^2(\Omega)$, we split the  error into two parts as follows,
$$
v-\bP v= v-{\hat{\bP}}_N^k v= v -\bP^k_{N} v+\bP^k_{N} v -{\hat{\bP}}_N^k v.
$$
The term $ v-\bP_N^k v$ can be estimated by Property \ref{prop:tensor}. Therefore we only need to bound
\begin{align*}
\bP^k_{N} v -{\hat{\bP}}_N^k v&=\sum_{\substack{ |\bl|_\infty \leq N\\\bl \in \mathbb{N}_0^d}} {Q}^k_{l_1, x_1} \otimes \cdots \otimes {Q}^k_{l_d, x_d} v -\sum_{\substack{ |\bl|_1 \leq N\\\bl \in \mathbb{N}_0^d}} {Q}^k_{l_1, x_1} \otimes \cdots \otimes {Q}^k_{l_d, x_d} v \\
&=\sum_{\substack{ |\bl|_\infty \leq N, |\bl|_1 > N\\\bl \in \mathbb{N}_0^d}} {Q}^k_{l_1, x_1} \otimes \cdots \otimes {Q}^k_{l_d, x_d} v.
\end{align*}

In what follows, we will estimate   the tensor-product construction of increment projections when $ |\bl|_\infty \leq N, |\bl|_1>N$.  
 Assume $u\in\mathcal{H}^{p+1}(\Omega)$, and $q$ is an integer with $1\leq q\leq \min\{p,k\}. $ For a multi-index $\bl$ with $ |\bl|_\infty \leq N, |\bl|_1>N$,
let $L = supp(\bl):=\{i_1,i_2,\cdots,i_r\}\subset\{1,2,\cdots,d\}$, i.e.,  $l_{i_s}\neq 0$ iff $i_s\in L$, and $r = |L|$. Combining the approximation property of the one-dimensional increment operator $Q_n^k$ in \eqref{eq:proj1} and \eqref{eq:proj3} with Proposition 5.1 in  \cite{schwab2008sparse} gives that, on each elementary cell $I_\bl^\bj \in \Omega_\bl$ with $L=supp(\bl)=\{i_1,i_2,\cdots,i_r\}$
\begin{align}
\left\|Q^k_{l_1,x_1}\otimes\cdots\otimes Q^k_{l_d,x_d} v\right\|^2_{L^2(I^\bj_\bl)}
\leq& \tilde{c}^{2r}_{k,0,q}\hat{c}^{2(d-r)}_{k,0}4^{-(q+1)|\bl|_1}\left\|\left(\frac{\partial^{q+1}}{\partial x_{i_1}^{q+1}}\ldots \frac{\partial^{q+1}}{\partial x_{i_r}^{q+1}}\right)v\right\|^2_{L^2(I^\bj_\bl)}\notag\\
=& \tilde{c}^{2r}_{k,0,q}\hat{c}^{2(d-r)}_{k,0}4^{-(q+1)|\bl|_1}|v|^2_{H^{q+1,L}(I_\bl^\bj)}, \label{eq:proj_l2}
\end{align}
\begin{align}
&\left|Q^k_{l_1,x_1}\otimes\cdots\otimes Q^k_{l_d,x_d} v\right|^2_{H^1(I^\bj_\bl)}\notag\\=&\sum_{m=1}^d
\left\|\frac{\partial}{\partial x_m}\left(Q^k_{l_1,x_1}\otimes\cdots\otimes Q^k_{l_d,x_d} \right)v\right\|^2_{L^2(I^\bj_\bl)} \notag \\
=&\sum_{m\in L}\left\|\frac{\partial}{\partial x_m}\left(Q^k_{l_1,x_1}\otimes\cdots\otimes Q^k_{l_d,x_d} \right)v\right\|^2_{L^2(I^\bj_\bl)} + \sum_{m\not\in L}\left\|\frac{\partial}{\partial x_m}\left(Q^k_{l_1,x_1}\otimes\cdots\otimes Q^k_{l_d,x_d} \right)v\right\|^2_{L^2(I^\bj_\bl)}\notag\\
\leq&\sum_{m\in L} \tilde{c}^{2(r-1)}_{k,0,q}\tilde{c}^{2}_{k,1,q}\hat{c}^{2(d-r)}_{k,0}4^{-(q+1)|\bl|_1+l_m}\left\|\left(\frac{\partial^{q+1}}{\partial x_{i_1}^{q+1}}\ldots \frac{\partial^{q+1}}{\partial x_{i_r}^{q+1}}\right)v\right\|^2_{L^2(I^\bj_\bl)}\notag\\
&+\sum_{m\not\in L} \tilde{c}^{2r}_{k,0,q}\hat{c}^{2}_{k,1}\hat{c}^{2(d-r-1)}_{k,0}4^{-(q+1)|\bl|_1}\left\|\left(\frac{\partial^{q+1}}{\partial x_{i_1}^{q+1}}\ldots \frac{\partial^{q+1}}{\partial x_{i_r}^{q+1}}\right)v\right\|^2_{L^2(I^\bj_\bl)}\notag\\
= &\sum_{m\in L} \tilde{c}^{2(r-1)}_{k,0,q}\tilde{c}^{2}_{k,1,q}\hat{c}^{2(d-r)}_{k,0}4^{-(q+1)|\bl|_1+l_m}\left|v\right|^2_{H^{q+1,L}(I^\bj_\bl)}+\sum_{m\not\in L} \tilde{c}^{2r}_{k,0,q}\hat{c}^{2}_{k,1}\hat{c}^{2(d-r-1)}_{k,0}4^{-(q+1)|\bl|_1}\left|v\right|^2_{H^{q+1,L}(I^\bj_\bl)}\notag\\
=& \left( \tilde{c}^2_{k,1,q}\hat{c}^{2}_{k,0}\sum_{m\in L}4^{l_m} + (d-r)\tilde{c}^2_{k,0,q}\hat{c}^2_{k,1}\right)\tilde{c}^{2(r-1)}_{k,0,q}\hat{c}^{2(d-r-1)}_{k,0}4^{-(q+1)|\bl|_1}\left|v\right|^2_{H^{q+1,L}(I^\bj_\bl)}\notag\\
\leq& d \bar{c}_{k,q}\tilde{c}^{2(r-1)}_{k,0,q}\hat{c}^{2(d-r-1)}_{k,0}4^{-(q+1)|\bl|_1+|\bl|_\infty}\left|v\right|^2_{H^{q+1,L}(I^\bj_\bl)},
 \label{eq:proj_h1}
\end{align}
where $\bar{c}_{k,q}=\max\left( \tilde{c}^2_{k,1,q}\hat{c}^{2}_{k,0}, \tilde{c}^2_{k,0,q}\hat{c}^2_{k,1}\right)$.

Hence, summing up on all elements,
\begin{align*}
&\left\|Q^k_{l_1,x_1}\otimes\cdots\otimes Q^k_{l_d,x_d}  v \right\|^2_{L^2(\Omega_N)}
=\left\|Q^k_{l_1,x_1}\otimes\cdots\otimes Q^k_{l_d,x_d}  v \right\|^2_{L^2(\Omega_\bl)}\\
=&\sum_{\bzer\leq \bj\leq 2^\bl-\mathbf{1}} \|Q^k_{l_1,x_1}\otimes\cdots\otimes Q^k_{l_d,x_d} v\|^2_{L^2(I_\bl^\bj)}\\
\leq&\sum_{\bzer\leq \bj\leq 2^\bl-\mathbf{1}}\tilde{c}^{2r}_{k,0,q}\hat{c}^{2(d-r)}_{k,0}4^{-(q+1)|\bl|_1}|v|^2_{H^{q+1,L}(I_\bl^\bj)}
=\tilde{c}^{2r}_{k,0,q}\hat{c}^{2(d-r)}_{k,0}4^{-(q+1)|\bl|_1}|v|^2_{H^{q+1,L}(\Omega)},
 \end{align*}
and
\begin{align*}
&\left|Q^k_{l_1,x_1}\otimes\cdots\otimes Q^k_{l_d,x_d} v\right|^2_{H^1(\Omega_N)}=\left|Q^k_{l_1,x_1}\otimes\cdots\otimes Q^k_{l_d,x_d} v\right|^2_{H^1(\Omega_\bl)}\\
 =&\sum_{\bzer\leq \bj\leq 2^\bl-\mathbf{1}} \left|Q^k_{l_1,x_1}\otimes\cdots\otimes Q^k_{l_d,x_d} v\right|^2_{H^1(I_\bl^\bj)}\\
\leq&\sum_{\bzer\leq \bj\leq 2^\bl-\mathbf{1}}d \bar{c}_{k,q}\tilde{c}^{2(r-1)}_{k,0,q}\hat{c}^{2(d-r-1)}_{k,0}4^{-(q+1)|\bl|_1+|\bl|_\infty}\left|v\right|^2_{H^{q+1,L}(I^\bj_\bl)}\\
=&d \bar{c}_{k,q}\tilde{c}^{2(r-1)}_{k,0,q}\hat{c}^{2(d-r-1)}_{k,0}4^{-(q+1)|\bl|_1+|\bl|_\infty}|v|^2_{H^{q+1,L}(\Omega)}.
 \end{align*}
The above estimate is valid since $|\bl|_\infty \leq N$, which implies   the function $Q^k_{l_1,x_1}\otimes\cdots\otimes Q^k_{l_d,x_d}  v$ is supported on grid $\Omega_N,$ (i.e. $Q^k_{l_1,x_1}\otimes\cdots\otimes Q^k_{l_d,x_d}  v$ is a polynomial in each element of  $\Omega_N$).

Now, gathering the results, 
for the $L^2$ norm, we have
\begin{align*}
&\| \bP^k_{N} v -{\hat{\bP}}_N^k v \|_{L^2(\Omega_N)}
\le  \sum_{\substack{ |\bl|_\infty \leq N, |\bl|_1 > N\\\bl \in \mathbb{N}_0^d}} \|Q^k_{l_1,x_1}\otimes\cdots\otimes Q^k_{l_d,x_d} v\|_{L^2(\Omega_N)}\\
& \leq \sum_{\substack{ |\bl|_\infty \leq N, |\bl|_1 > N\\\bl \in \mathbb{N}_0^d}} \tilde{c}^{r}_{k,0,q}\hat{c}^{(d-r)}_{k,0}2^{-(q+1)|\bl|_1}|v|_{H^{q+1,L}(\Omega)}, \qquad \textrm{with} \quad L=supp(\bl), r=|L|\\
& \leq \sum_{\substack{ |\bl|_\infty \leq N, |\bl|_1 > N\\\bl \in \mathbb{N}_0^d}} C_{k,q}^d2^{-(q+1)(N+1)} |v|_{\mathcal{H}^{q+1}(\Omega)}=   C_{k,q}^d2^{-(q+1)(N+1)} |v|_{\mathcal{H}^{q+1}(\Omega)} \sum_{\substack{ |\bl|_\infty \leq N, |\bl|_1 > N\\\bl \in \mathbb{N}_0^d}}  1,
\end{align*}
where $C_{k,q}=\max(\tilde{c}_{k,0,q}, \hat{c}_{k,0}).$ Now, we estimate
\begin{align}
&\sum_{\substack{ |\bl|_\infty \leq N, |\bl|_1 > N\\\bl \in \mathbb{N}_0^d}}  1=\sum_{\substack{ |\bl|_\infty \leq N\\\bl \in \mathbb{N}_0^d}}  1-\sum_{\substack{ |\bl|_1 \le N\\\bl \in \mathbb{N}_0^d}}  1 =(N+1)^d-\sum_{s=0}^N \sum_{\substack{ |\bl|_1 =s\\\bl \in \mathbb{N}_0^d}}  1\notag\\
&= (N+1)^d-\sum_{s=0}^N \binom{s+d-1}{d-1}=(N+1)^d- \binom{N+d}{d} \le (N+1)^d (1-\frac{1}{d!}). \label{eq:countl}
\end{align}
Therefore,
\begin{align*}
&\| \bP^k_{N} v -{\hat{\bP}}_N^k v \|_{L^2(\Omega_N)} \le C_{k,q}^d2^{-(q+1)(N+1)} (N+1)^d (1-\frac{1}{d!}) |v|_{\mathcal{H}^{q+1}(\Omega)}\\
&  \le C_{k,q}^d2^{-(q+1)(N+1)} (N+1)^d  |v|_{\mathcal{H}^{q+1}(\Omega)} = A_q \kappa(k,q,N)^d 2^{-N(q+1)}   |v|_{\mathcal{H}^{q+1}(\Omega)},
\end{align*}
with $\kappa(k,q,N)=C_{k,q}(N+1), \,A_q=2^{-(q+1)}.$ Combining with \eqref{eq:tensorproj}, we arrive at the estimate
\begin{align}
&\| \bP v - v \|_{L^2(\Omega_N)} \le \bar{\bar{c}}_{k,0, q} 2^{-N(q+1)} |v|_{H^{q+1}(\Omega)} +A_q \kappa(k,q,N)^d 2^{-N(q+1)}   |v|_{\mathcal{H}^{q+1}(\Omega)}\notag \\
& \le \left(\bar{\bar{c}}_{k,0, q}  +A_q \kappa(k,q,N)^d \right)2^{-N(q+1)}   |v|_{\mathcal{H}^{q+1}(\Omega)}. \label{eq:estl2}
\end{align}

For the $H^1$ broken semi-norm, similarly we have
\begin{align*}
&| \bP^k_{N} v -{\hat{\bP}}_N^k v |_{H^1(\Omega_N)}
\le  \sum_{\substack{ |\bl|_\infty \leq N, |\bl|_1 > N\\\bl \in \mathbb{N}_0^d}} |Q^k_{l_1,x_1}\otimes\cdots\otimes Q^k_{l_d,x_d} v |_{H^1(\Omega_N)}\\
& \leq  \sum_{\substack{ |\bl|_\infty \leq N, |\bl|_1 > N\\\bl \in \mathbb{N}_0^d}} \sqrt{d \bar{c}_{k,q}}\tilde{c}^{(r-1)}_{k,0,q}\hat{c}^{(d-r-1)}_{k,0}2^{-(q+1)|\bl|_1+|\bl|_\infty}|v|_{H^{q+1,L}(\Omega)}, \qquad \textrm{with} \quad L=supp(\bl), r=|L|\\
& \leq \sqrt{d \bar{c}_{k,q}}C_{k,q}^{d-2} |v|_{\mathcal{H}^{q+1}(\Omega)}  \sum_{\substack{ |\bl|_\infty \leq N, |\bl|_1 > N\\\bl \in \mathbb{N}_0^d}} 2^{-(q+1)|\bl|_1+|\bl|_\infty}\\
& \leq \sqrt{d \bar{c}_{k,q}}C_{k,q}^{d-2} |v|_{\mathcal{H}^{q+1}(\Omega)} \sum_{s=N+1}^{Nd}  2^{-(q+1)s} \sum_{\substack{ |\bl|_\infty \leq N, |\bl|_1 =s\\\bl \in \mathbb{N}_0^d}} 2^{|\bl|_\infty}.
\end{align*}
In  \cite{schwab2008sparse}, it was shown that 
$$
\sum_{\substack{|\bl|_1 =s\\\bl \in \mathbb{N}_0^d}} 2^{|\bl|_\infty} \le d \,2^{d-1+s}.
$$
Therefore,
\begin{align*}
&| \bP^k_{N} v -{\hat{\bP}}_N^k v |_{H^1(\Omega_N)} \le \sqrt{d \bar{c}_{k,q}}C_{k,q}^{d-2} |v|_{\mathcal{H}^{q+1}(\Omega)} \sum_{s=N+1}^{Nd}  2^{-(q+1)s}  d \,2^{d-1+s}\\
& \le \sqrt{d \bar{c}_{k,q}}C_{k,q}^{d-2} |v|_{\mathcal{H}^{q+1}(\Omega)} d 2^{d-1}   \sum_{s=N+1}^{Nd}  2^{-qs}\\
& \le \sqrt{d \bar{c}_{k,q}}C_{k,q}^{d-2} |v|_{\mathcal{H}^{q+1}(\Omega)} d 2^{d-1}    2^{-qN}\\
& \leq d^{3/2}  B_{k,q} (2 C_{k,q})^d 2^{-Nq} |v|_{\mathcal{H}^{q+1}(\Omega)},
\end{align*}
where $B_{k,q} =  \sqrt{\bar{c}_{k,q}} C_{k,q}^{-2}/2$. Combining with  \eqref{eq:tensorproj}, we get
\begin{align}
&| \bP v -v |_{H^1(\Omega_N)} \le \bar{\bar{c}}_{k,1, q} 2^{-Nq} |v|_{H^{q+1}(\Omega)}+d^{3/2}  B_{k,q} (2 C_{k,q})^d 2^{-Nq} |v|_{\mathcal{H}^{q+1}(\Omega)}\notag \\
& \le \left (\bar{\bar{c}}_{k,1, q}  +d^{3/2}  B_{k,q} (2 C_{k,q})^d \right ) 2^{-Nq} |v|_{\mathcal{H}^{q+1}(\Omega)}, \label{eq:esth1}
\end{align}
which completes the proof of Lemma  \ref{thm:appx}.

\hfill\ensuremath{\blacksquare}

\bibliographystyle{abbrv}
\bibliography{ref_cheng,ref_cheng_2}

\end{document}